\documentclass[letterpaper,11pt,sort&compress]{elsarticle}

\usepackage[margin=1in]{geometry}

\usepackage{amsmath}
\usepackage{amssymb}
\usepackage{amsfonts}
\usepackage{amsthm}
\usepackage{arydshln}
\usepackage{cancel}
\usepackage{bm}
\usepackage{multirow}
\usepackage{booktabs}
\usepackage{natbib}
\usepackage[labelfont=bf]{caption}
\usepackage{subcaption}
\usepackage[colorlinks]{hyperref}
\usepackage{fancyhdr}

\usepackage{lineno,hyperref}

\hypersetup{pdfauthor={Name}}

\begin{document}

\let\today\relax
\makeatletter
\def\ps@pprintTitle{%
    \let\@oddhead\@empty
    \let\@evenhead\@empty
    }
\makeatother

\begin{frontmatter}

\title{Quadratic Point Estimate Method for Probabilistic Moments Computation}

\author[a]{Minhyeok Ko}

\author[a]{Konstantinos G. Papakonstantinou\corref{cor}}
\ead{kpapakon@psu.edu}

\cortext[cor]{Corresponding author}

\address[a]{Department of Civil and Environmental Engineering, The Pennsylvania State University, University Park, PA 16802, USA}

\begin{abstract}
This paper presents in detail the originally developed Quadratic Point Estimate Method (QPEM), aimed at efficiently and accurately computing the first four output moments of probabilistic distributions, using $2n^2+1$ sample (or sigma) points, with $n$, the number of input random variables. The proposed QPEM particularly offers an effective, superior, and practical alternative to existing sampling and quadrature methods for low- and moderately-high-dimensional problems. Detailed theoretical derivations are provided proving that the proposed method can achieve a fifth or higher-order accuracy for symmetric input distributions. Various numerical examples, from simple polynomial functions to nonlinear finite element analyses with random field representations, support the theoretical findings and further showcase the validity, efficiency, and applicability of the QPEM, from low- to high-dimensional problems.
\end{abstract}


\end{frontmatter}

\section{Introduction}
\label{sec:intro}
In many engineering and scientific applications, uncertainty propagation through a nonlinear system often entails calculating probabilistic moment integrals related to a probability density function (PDF). Computing these integrals is a critical component in determining the output probabilistic characteristics \cite{zhao2000third,xu2004dimension,fan2016adaptive,gottardi2020reliability,zhou2020adaptive,zhao2021structural,liu2023efficient,li2023efficient,yuan2023sample}. Analytical solutions for such integrals, in general cases, do not exist and thus numerical estimations are required. The Monte Carlo (MC) method is the most versatile and widely used numerical sampling technique, utilized in various scientific fields \cite{decker1991monte, amar2006monte,liu2008monte,mordechai2011applications,kroese2013handbook,rubinstein2016simulation}. However, MC generally requires a considerable number of samples, for sufficiently low precision levels, which may not be computationally feasible, particularly in engineering problems involving computationally intensive simulations.

Variance reduction techniques, such as quasi-Monte Carlo methods \cite{caflisch1998monte,sobol1998quasi,jank2005quasi} and Latin Hypercube Sampling \cite{mckay1979lhs,olsson2002latin, shields2016generalization}, can serve as efficient alternatives to the MC method and are characterized by sampling points that fill the space more uniformly. These techniques offer the advantage of faster estimator convergence rates over the MC method under certain conditions and dimensions.

Quadrature techniques have also been extensively utilized to calculate the moment integrals and, overall, work well in low-dimensional spaces or deterministic cases. However, probabilistic problems are usually described by numerous random variables, particularly cases involving random processes and fields, e.g., \cite{shinozuka1991simulation, papakonstantinou2013probabilistic}, rendering quadrature techniques prohibitive in such settings due to the involved high-dimensionality. The required number of quadrature points in most typical cases, and thus the computational cost, increases exponentially with the dimensions, also known as the curse of dimensionality. Sparse quadrature techniques have been thus also devised \cite{smolyak1963quadrature,gerstner1998numerical}, aiming to address this limitation by utilizing sparse tensor products. A Smolyak sparse-grid quadrature scheme \cite{holtz2010sparse} is used in this work for comparison purposes with our developed approach.   

Dedicated methods for computing probabilistic moments have also been suggested in the literature, and particularly the Point Estimate Method (PEM) has a long history in engineering. The concept of PEM was first introduced by Rosenblueth \cite{rosenblueth1975point}. Since then, many PEM variations and variants have been developed, mainly of exponential and linear computationally complexity with the number of dimensions, offering generally powerful yet simple approximation methods for estimation of the first few moments of an output PDF \cite{harr1989, hong1998, julier1997new, zhao2000, julier2002scaled, julier2002spherical, tenne2003higher, mohammadi2013new, adurthi2018conjugate, papakonstantinou2022scaled}. Since their introduction almost 50 years ago, PEMs have found very broad applicability including also in energy applications, geotechnical engineering, and nonlinear filtering, among others \cite{morales2010a, delgado2014, christian1999, chang1997, evangelopoulos2014optimal, bordbari2018probabilistic, xu2021novel, xu2019adaptive, arasaratnam2009cubature, jia2013high, che2020probabilistic, jiang2019train}. PEM methods evaluate the response function at deterministic sample points, also known as sigma points, and thus their main computational cost is relevant to the number of used sigma points. In many cases, PEM-based methods can provide notably accurate mean and variance estimations with pronounced computational efficiency. However, most existing PEMs are unsuitable for estimating higher-order moments, such as skewness and kurtosis \cite{christian1999}, with non-exponential computational costs with respect to random dimensions, since, generally, numerous sigma points are required to capture higher-order moments. Hence, there is a natural trade-off between computational efficiency and accuracy in relation to the number of sigma points used. Many PEMs also often use an optimization procedure to find the suitable sigma point locations and weights \citep{rosenblueth1975point} for each problem, which can complicate solution approaches and implementation in such cases.

Motivated by the above, this work develops and presents a new general method, the Quadratic Point Estimate Method (QPEM), able to completely capture up to fifth-order input moments of the normal distribution with a quadratic expansion of the required sigma points in relation to random dimensions, and without requiring any optimization procedure to determine the sigma point characteristics. Analogous ideas have appeared, in a broad sense, in \cite{mcnamee1967construction,genz1996fully}, based, however, on a pure quadrature logic, and their rationale, algorithms, applicability, and implementations for calculating multidimensional moment integrals are completely different from our work here. The suggested QPEM can actually provide exact odd-order input moments for normal and all other symmetric distributions, due to its symmetric scheme, and can also minimize the sixth- and eighth-order input moments errors merely through appropriate scaling parameters and without any additional sampling points. 

Detailed theoretical derivations are provided in this paper to fully demonstrate the devised QPEM, and a comprehensive background section appropriately places it among other relevant methods in the literature, used also in this work for comparison purposes. The validity and outstanding performance of QPEM is also showcased on numerous mathematical and engineering examples with different characteristics, varying dimensionality, from 5 to 100 input dimensions, and linear and nonlinear mechanics, among others.   

\section{Moment integrals}
\label{sec:review}
\subsection{Quadrature techniques}
Let us consider the problem of calculating the moment integral of a real-valued function $\mathcal{M}(\mathbf{x}):\mathbb{R}^n \rightarrow \mathbb{R}$ with respect to a multivariate PDF $f(\mathbf{x})$. Quadrature techniques approximate the multidimensional moment integral with a weighted sum of function evaluations computed at specific quadrature points, $\mathbf{X}_i \in \mathbb{R}^n$, and with weights, $W_i \in \mathbb{R}$, as follows:
\begin{equation}
	\label{eq:moment_integral}
	E[\mathcal{M}(\mathbf{x})]=\int \mathcal{M}(\mathbf{x})f(\mathbf{x}) d\mathbf{x} \approx \sum_{i=1}^{N} W_i \mathcal{M}(\mathbf{X}_i)
\end{equation}
where $\mathbf{x}$ is the random vector $\{x_1, \ldots, x_n\}^T$, and $N$ is the number of quadrature points. Assuming $\mathcal{M}(\mathbf{x})$ has a valid Taylor series expansion about the mean, $\bm{\mu}_{\mathbf{x}}=\left\{\mu_1,\ldots, \mu_n\right\}^{T}$, the expectation term in Eq.~\eqref{eq:moment_integral} can be expanded as:
\begin{equation}
	\label{eq:taylor}
	E[\mathcal{M}(\mathbf{x})] = \sum_{k_1=0}^{\infty} \cdots \sum_{k_n=0}^{\infty} \dfrac{E\left[dx_1^{k_1}\cdots dx_n^{k_n}\right]}{k_1 ! \cdots k_n !}  \left. \dfrac{\partial^{k_1+\cdots+k_n} \mathcal{M}(\mathbf{x})}{\partial x_{1}^{k_1} \cdots \partial x_{n}^{k_n}} \right|_{\mathbf{x}=\bm{\mu}_{\mathbf{x}}}
\end{equation}
where $dx_j=x_j-\mu_{j}$ and $k_j$ is a positive integer or zero. An important observation in Eq.~\eqref{eq:taylor} is that evaluating the expectation of the function $\mathcal{M}(\mathbf{x})$ is related to calculating the moments of the input PDF. One can compute an increasingly accurate moment integral by increasing the number of terms in the Taylor series expansion. Substituting now the Taylor series expansion of $\mathcal{M}(\mathbf{X}_i)$ to Eq.~\eqref{eq:moment_integral} results in:
\begin{equation}
	\label{eq:quadrature_taylor}
	\sum_{i=1}^{N} W_i \mathcal{M}(\mathbf{X}_i) = \sum_{k_1=0}^{\infty} \cdots \sum_{k_n=0}^{\infty} \dfrac{\sum_{i=1}^{N} W_i \left\{dX_{i,1}^{k_1} \cdots dX_{i,n}^{k_n} \right\}}{k_1 ! \cdots k_n !}  \left. \dfrac{\partial^{k_1+\cdots+k_n} \mathcal{M}(\mathbf{x})}{\partial x_{1}^{k_1} \cdots \partial x_{n}^{k_n}} \right|_{\mathbf{x}=\bm{\mu}_{\mathbf{x}}}
\end{equation} 
where $dX_{i,j}=X_{i,j}-\mu_{j}$, and  $X_{i,j}$ is the $j$th coordinate of each $\mathbf{X}_i$. Comparing Eqs.~\eqref{eq:taylor}-\eqref{eq:quadrature_taylor} results in the following set of equations, known also as the moment constraint equations (MCEs) \cite{adurthi2018conjugate}:
\begin{equation}
	\label{eq:moment_equation}
	E\left[dx_1^{k_1}\cdots dx_n^{k_n}\right] = \sum_{i=1}^{N} W_i \left\{dX_{i,1}^{k_1} \cdots dX_{i,n}^{k_n} \right\}
\end{equation}
where $k_1+\cdots+k_n=d$ represents the order of the input moment. Alternatively, the MCEs can be also determined by stating that a set of quadrature points $\mathbf{X}_i$ with weights $W_i$ can exactly integrate polynomials of total degree $d$ or less. This is also called the moment matching method \cite{adurthi2018conjugate}. 

In practice, one may only need the minimum number of quadrature points that can accurately capture the first $d$ moments of interest. Conversely, when computing moments for non-polynomial functions, commonly used in engineering, it is generally a priori unclear how many points are needed for exact evaluation. In such cases, reducing errors in moment evaluation for non-polynomial functions may require higher-order quadrature schemes. The MCEs are typically high-order multivariate polynomial equations, the solutions of which are nontrivial and sometimes intractable. Quadrature techniques thus work well theoretically for $n$-$\mathcal{D}$ MCEs by taking the tensor product of the higher-order 1-$\mathcal{D}$ quadrature points. However, the number of quadrature points then increases exponentially with the dimensions, which is a major computational obstacle for most engineering problems.

Sparse quadrature rules have been hence developed, in contrast to quadrature techniques with full tensor products, able to achieve the same numerical accuracy with significantly fewer points. One of the most widely used sparse-grid quadrature schemes is due to Smolyak \cite{smolyak1963quadrature}, that also utilizes negative weights in the computations \cite{wasilkowski1995explicit}, with implications related to the so-called stability factor and overall accuracy, as explained subsequently.

\subsection{Point Estimate Methods}
\label{subsec:pem}
 Point Estimate Methods are dedicated probabilistic methods, specialized for such problems. Rosenblueth's Point Estimate Method (RPEM) is based on capturing selected input moments of a random variable using a two-point moment matching method \cite{rosenblueth1975point}. The RPEM computational requirements increase significantly with an increasing number of random variables, due to the use of the tensor product rule, requiring $2^n$ function evaluations in $n$-$\mathcal{D}$ problems \cite{rosenblueth1981two}. Depending on the selected input moments, an optimization procedure may also be required to determine the sigma points and associated weights, usually offering an approximate, non-unique solution \cite{tsai2005evaluation}. 

Various PEMs have been developed over the years to overcome the aforementioned drawbacks of the RPEM in $n$-$\mathcal{D}$ problems, such as Harr’s PEM (HaPEM) \cite{harr1989}, Hong’s PEM (HPEM) \cite{hong1998}, and Unscented Transformation (UT) \cite{julier1997new}, which are all offering linear complexity as dimensions grow. These PEMs generally include a decorrelation process for correlated random variables through orthogonal rotation. HaPEM uses the eigenvectors of the correlation matrix to decorrelate the input variables, and the sigma points are then selected at the intersections of the eigenvectors and $\sqrt{n}$ radius hypersphere. The total number of sigma points in HaPEM is thus $2n$, significantly reduced in comparison to RPEM and linearly increasing now with dimensions $n$. UT decorrelates the input random variables based on the square root matrix of the covariance of the $n$-$\mathcal{D}$ variables. Similar to HaPEM, symmetric sigma points on the $\sqrt{n}$ radius hypersphere are selected to retrieve the first two input moments, with an additional sigma point located now on the mean point of the random variables, for a more accurate and stable estimation, ending up with $2n+1$ total sigma points. HPEM in turn, similarly uses the square root matrix of the covariance for the decorrelation and shares many other similarities with the UT, yet it can be considered a more general UT variant, as explained subsequently. Several other notable PEMs, aiming to achieve higher accuracy or computational efficiency can be seen in \cite{zhao2000, julier2002scaled, julier2002spherical, tenne2003higher, mohammadi2013new, adurthi2018conjugate, papakonstantinou2022scaled}. Since HPEM \cite{hong1998} is likely the most general linear PEM, used also in a wide-range of applications over the years, it will be utilized for comparison purposes in this work and is thus further described accordingly. 

\begin{figure}
	\begin{subfigure}[b]{0.5\textwidth}
		\centering
		\includegraphics[trim=0in 0in 0in 0in, clip=true, width=\textwidth]{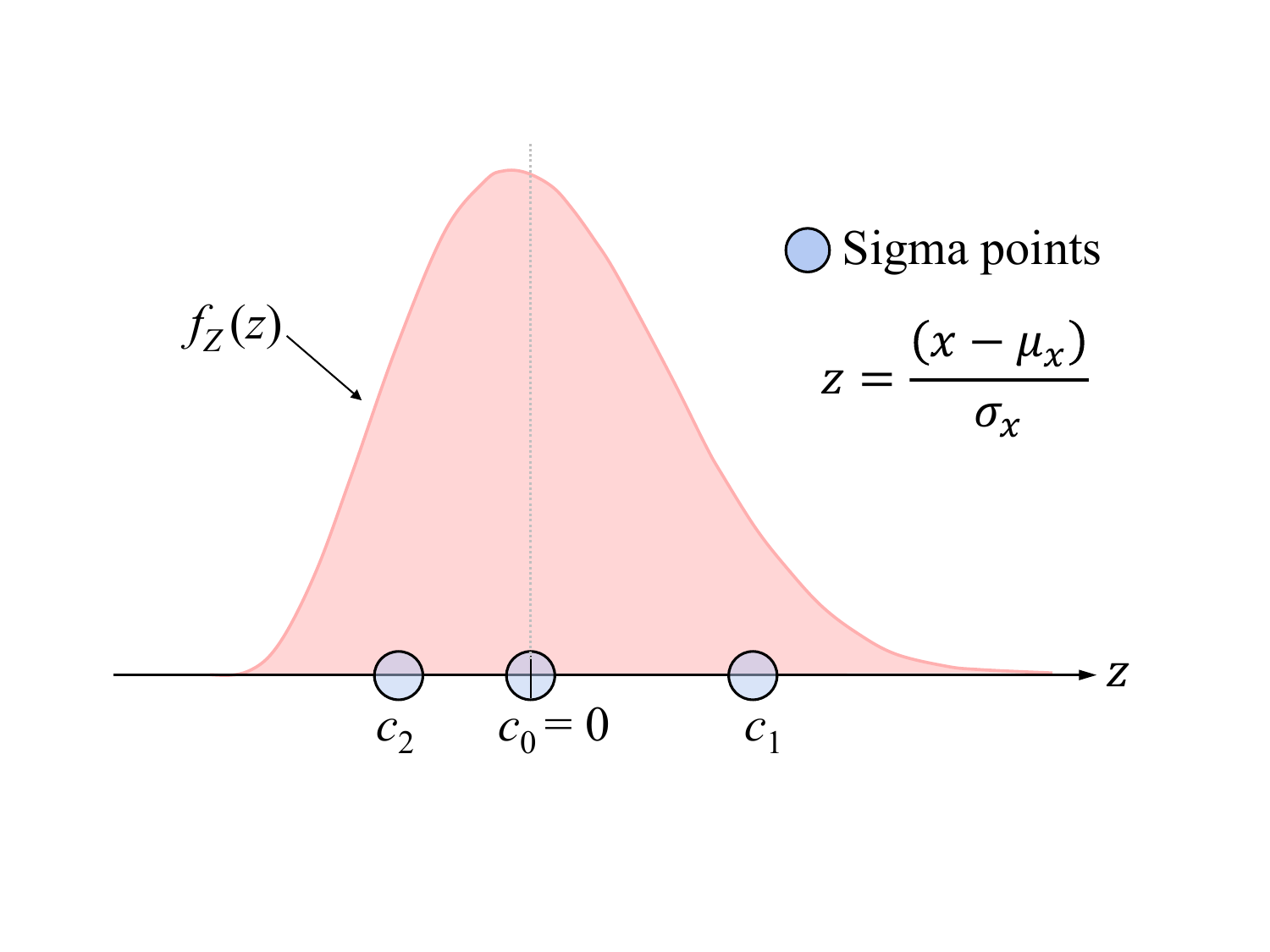}
		\caption{$1$-$\mathcal{D}$ space}
		\label{fig:hopem_a}
	\end{subfigure}
	\begin{subfigure}[b]{0.5\textwidth}
		\centering
		\includegraphics[trim=0in 0in 0in 0in, clip=true, width=\textwidth]{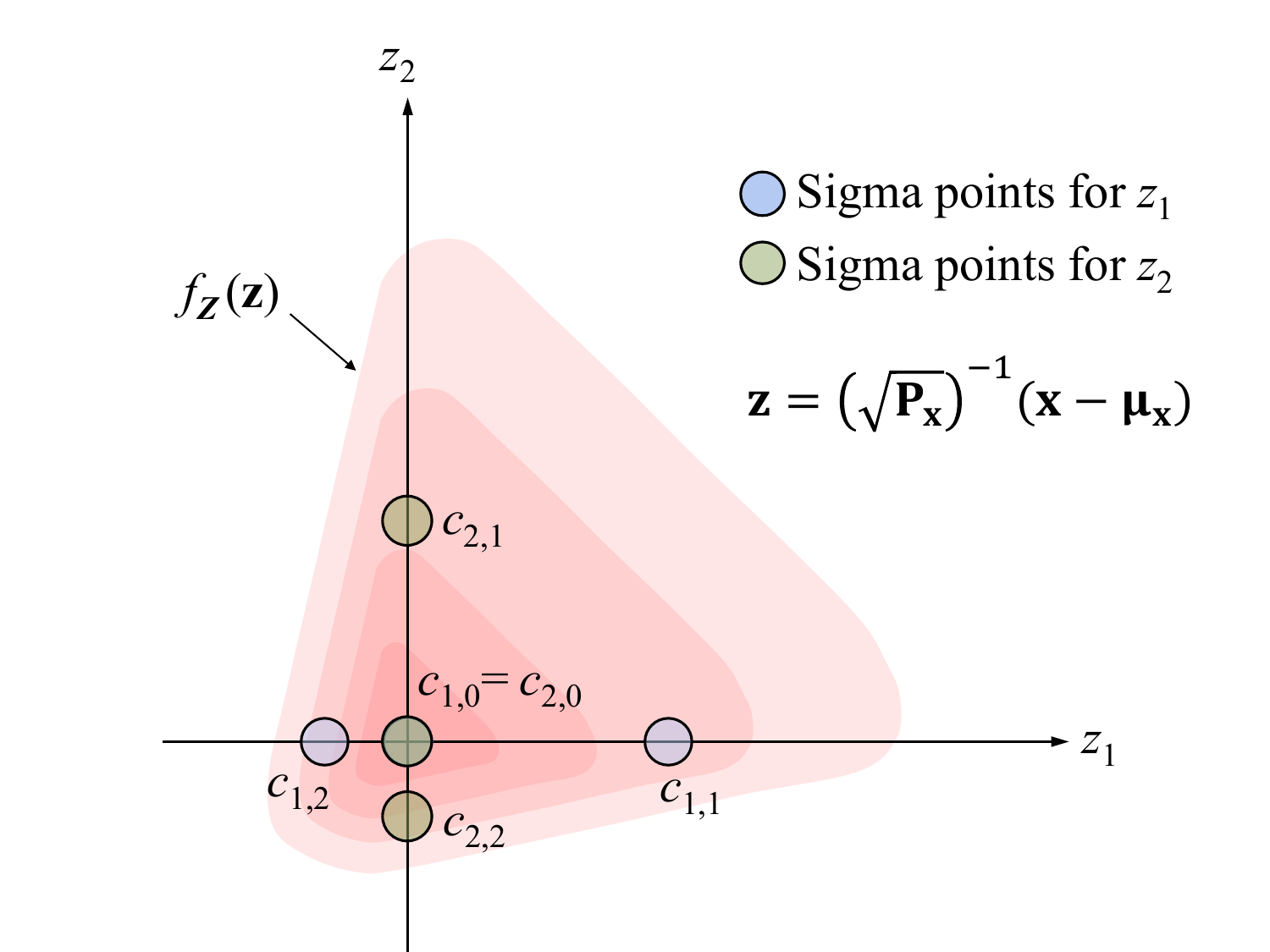}
		\caption{$2$-$\mathcal{D}$ space}
		\label{fig:hopem_b}
	\end{subfigure}
	\caption{Schematic representation of illustrative $2n+1$ PEM.}
	\label{fig:hopem}
\end{figure}

As mentioned above, HPEM was developed to overcome the computational inefficiency of the RPEM in multidimensional problems. The HPEM sigma points and associated weights are accordingly determined so that their moments match the ones of the actual PDF using the moment matching method. Without loss of generality, consider a $n$-$\mathcal{D}$ standardized random variable $\mathbf{z}$, with zero mean $\mathbf{0}$ and identity covariance $\mathbf{I}$. Any $n$-$\mathcal{D}$ random variable $\mathbf{x}$ with mean $\bm{\mu}_{\mathbf{x}}$ and covariance $\mathbf{P}_{\mathbf{x}}$ can be obtained by an affine transformation, as:
\begin{equation}
	\label{eq:affine_transformation1}
	\mathbf{z}=\sqrt{\mathbf{P}_{\mathbf{x}}}^{-1} \left(\mathbf{x}-\bm{\mu}_{\mathbf{x}}\right) %
	\qquad \leftrightarrow \qquad \mathbf{x}=\bm{\mu}_{\mathbf{x}}+\sqrt{\mathbf{P}_{\mathbf{x}}} \mathbf{z}
\end{equation}
where $\sqrt{\mathbf{P_{{x}}}}\sqrt{\mathbf{P_{{x}}}}^{T}=\mathbf{P_{x}}$. 

Each uncorrelated variable $z_i$ is then supported by three points that can capture up to the fourth-order marginal moment of $z_i$. Placing one of the three points on the central point in $z$-space, i.e., $c_{i,0}=0$, as shown in Fig.~\ref{fig:hopem}, the MCEs up to the fourth-order marginal moment for each of the $n$-$\mathcal{D}$ random variables are formulated as: 

\begin{equation}
	\label{eq:hopem_me}
	\left\{\begin{array}{lll}
		E \left[1\right] &= \sum_{l=0}^{2} w_{i,l} &= 1/n\\ [8pt]
		E \left[z_{i}\right] &= \sum_{l=0}^{2} w_{i,l} c_{i,l}&= 0\\[8pt]
		E \left[z_{i}^2\right] &= \sum_{l=0}^{2} w_{i,l} c_{i,l}^2&= 1\\[8pt]
		E \left[z_{i}^3\right] &=  \sum_{l=0}^{2} w_{i,l} c_{i,l}^3&= \gamma_{i} \\[8pt]
		E \left[z_{i}^4\right] &=  \sum_{l=0}^{2} w_{i,l} c_{i,l}^4&= \kappa_{i}
	\end{array}\right.
\end{equation}
where $c_{i,l}$ and $w_{i,l}$ are coordinate and weight, respectively, of the $l$th sigma point for the variable $z_i$. $\gamma_{i}$ and $\kappa_{i}$ are, respectively, the skewness and kurtosis of $z_i$. The first line of Eq.~\eqref{eq:hopem_me} is based on the assumption that the summation of weights for each dimension $i$, and each variable $z_i$, is the same. Solving Eq.~\eqref{eq:hopem_me} for the sigma points and weights results in:

\begin{equation}
	\label{eq:hopem}
	\left\{ \begin{array}{ll}
		{{c}_{i,0}} &=0 \\ [8pt]
		{{c}_{i,1}} &=\dfrac{{{\gamma }_{i}}}{2}+\sqrt{{{\kappa }_{i}}-\dfrac{3\gamma _{i}^{2}}{4}} \\ [8pt]
		{{c}_{i,2}} &=\dfrac{{{\gamma }_{i}}}{2}-\sqrt{{{\kappa }_{i}}-\dfrac{3\gamma _{i}^{2}}{4}} \\ 
	\end{array} \right. \quad %
	\left\{ \begin{array}{ll}
		{{w}_{i,0}} &={1}/{n}\;-\sum_{l=1}^{2}{{{w}_{i,l}}} \\ [8pt]
		{{w}_{i,1}} &=\dfrac{1}{{{c}_{i,1}}\left( {{c}_{i,1}}-{{c}_{i,2}} \right)} \\ [8pt]
		{{w}_{i,2}} &=\dfrac{-1}{{{c}_{i,2}}\left( {{c}_{i,1}}-{{c}_{i,2}} \right)} \\ 
	\end{array} \right. \quad %
	i=1, \cdots , n
\end{equation}
By noting that ${{c}_{i,0}}$ for each variable is fixed at its central point 0, $n$ out of $3n$ sigma points are then superposed at the same point, ${{\mathbf{S}}_{0}}=\mathbf{0}_{n\times1}$ as shown in Fig.~\ref{fig:hopem}, and its weight, ${{W}_{0}}$, is provided by the summation of the $c_{i,0}$ weights, as:

\begin{equation}
	{{W}_{0}}=\sum_{i=1}^{n}{{{w}_{i,0}}}=1-\sum_{i=1}^{n}{\frac{1}{{{\kappa }_{i}}-\gamma _{i}^{2}}}
\end{equation}
Hence, $2n+1$ sigma points, $\mathbf{S}$=[${{\mathbf{S}}_{0}}$, ${{\mathbf{S}}_{1}}$, $\cdots$, ${{\mathbf{S}}_{2n}}$]$_{n \times (2n+1)}$, with associated weights, $\mathbf{W}$=[${{W}_{0}}$, ${{W}_{1}}$, $\cdots$, ${{W}_{2n}}$]$_{1 \times (2n+1)}$, are selected in the uncorrelated standardized $z$-space to capture the relevant moments of the random variables $\mathbf{z}$:
\begin{equation}
	\left\{ \begin{array}{*{35}{l}}
		{{\mathbf{S}}_{0}} & ={{\mathbf{0}}_{n\times 1}} & {}  \\[8pt]
		{{\mathbf{S}}_{i}} & ={{\{0,\cdots,\ \underbrace{{{c}_{i,1}}}_{{{i}}\,\text{th component}}, \cdots, 0\}}^{T}} & {}  \\[8pt]
		{{\mathbf{S}}_{i+n}} & ={{\{0,\cdots,\ \underbrace{{{c}_{i,2}}}_{{{i}}\, \text{th component}}, \cdots, 0\}}^{T}} & {}  \\
	\end{array} \right. %
	\left\{ \begin{array}{ll}
		{{W}_{0}} & =\displaystyle \sum_{i=1}^{n}{{{w}_{i,0}}} \\[8pt]
		{{W}_{i}} & ={{w}_{i,1}}  \\[8pt]
		{{W}_{i+1}} & ={{w}_{i,2}}  \\
	\end{array} \right. \,\,\,%
	\ i=1,\cdots , n
\end{equation}
 The sigma points, $\mathbf{S}_i$, are selected in $z$-space first, and are then transformed back to the original $x$-space by an appropriate transformation $\mathbf{X}_i$=$\mathcal{T}(\mathbf{S}_i)$, such as the affine transformation in case of $x$ normal variables. The transformed points can then act as input for the computational model, $\mathbf{Y}_i=\mathbf{M} \left(\mathbf{X}_i \right)$. The moments of the response output $\mathbf{y}$ can finally be estimated by Eqs.~\eqref{eq:mean_y}-\eqref{eq:moments_y}, with Eq.~\eqref{eq:moments_y} representing the $k$th-order moment vector estimation expressed with the Kronecker product, $\otimes$, e.g., $E$[$\left( \mathbf{y}-\mathbf{\bar{y}} \right)$${{\left( \mathbf{y}-\mathbf{\bar{y}} \right)}^{T}}$] $\equiv$ $E$[$\left( \mathbf{y}-\mathbf{\bar{y}} \right)$$\otimes$$\left( \mathbf{y}-\mathbf{\bar{y}} \right)$] = $E$[${{\left( \mathbf{y}-\mathbf{\bar{y}} \right)}^{2}}$]:
\begin{equation}
	\label{eq:mean_y}
	E \left[\mathbf{y}\right]=\overline{\mathbf{y}} \approx \sum_{i=0}^{N-1} W_i \mathbf{Y}_i
\end{equation}

\begin{equation}
	\label{eq:moments_y}
	E \left[\left(\mathbf{y}-\overline{\mathbf{y}}\right)^{k}\right] \approx \sum_{i=0}^{N-1} W_i \left(\mathbf{Y}_i-\overline{\mathbf{y}}\right)^{k}
\end{equation}

The HPEM, and related PEM methods with linear complexity, are thus very computationally efficient in estimating moments with just a small set of $2n+1$ weighted sigma points in $n$-dimensional spaces. However, these sigma points lie only on the orthogonal axes, as shown in Fig.~\ref{fig:hopem}, and therefore cannot capture any of the cross moments, i.e., $E[ z_{i}^{p}z_{j}^{q} ]$=$\int{z_{i}^{p}z_{j}^{q}f\left( \mathbf{z} \right)d\mathbf{z}}$. For example, consider the fourth-order cross moment $E[ z_{i}^{2}z_{j}^{2} ]$ of the uncorrelated standard normal distribution. The HPEM always provides $E[ z_{i}^{2}z_{j}^{2} ]$=0 by construction, whereas the actual uncorrelated standard normal distribution has $E[ z_{i}^{2}z_{j}^{2} ]$=1. Hence, HPEM first admits estimation errors for the fourth-order cross moments for the normal distribution. To overcome this cross moments inaccuracy, additional sigma points are needed along different axes other than the orthogonal ones.

\subsection{Sampling-based methods}
Sampling-based numerical methods, commonly referred to as variance reduction techniques, can effectively serve as alternatives to the standard Monte Carlo method, for general sampling cases and for moment integrals evaluation as well. Quasi-Monte Carlo (QMC) is one such family of sampling methods and employs low-discrepancy sequences, such as Halton or Sobol sequences, to generate sampling points \cite{caflisch1998monte}. QMC based on Sobol sequence is used in this work, for comparison purposes. 

Latin Hypercube Sampling (LHS) is another well known and powerful variance reduction technique that incorporates desirable features of random and stratified sampling. LHS is commonly used in efforts to produce accurate results with fewer model calls compared to MC method \cite{mckay1979lhs, olsson2002latin, shields2016generalization}. LHS is also used in this paper, for comparison with our developed QPEM approach.

\section{Quadratic Point Estimate Method}
\label{sec:qpem}
In this section, our entirely new PEM termed the Quadratic Point Estimate Method (QPEM) is described in detail. QPEM is developed in order to overcome shortcomings of other methods, mentioned earlier, and to increase estimation accuracy in a computationally affordable way. QPEM can completely capture up to fifth-order moments of the uncorrelated $n$-$\mathcal{D}$ standard  normal random variables $\mathbf{z}$, it is described by fully symmetric sigma points that can thus automatically satisfy the odd-order moments of zero for normal distributions, and offers general, analytical expressions for the sigma points and associated weights, avoiding any likely undesirable optimization process by users. It can also often offer reduced error estimations of higher-order even moments without any additional computational cost.

\begin{figure}
	\begin{subfigure}[b]{0.5\textwidth}
		\centering
		\includegraphics[width=\textwidth]{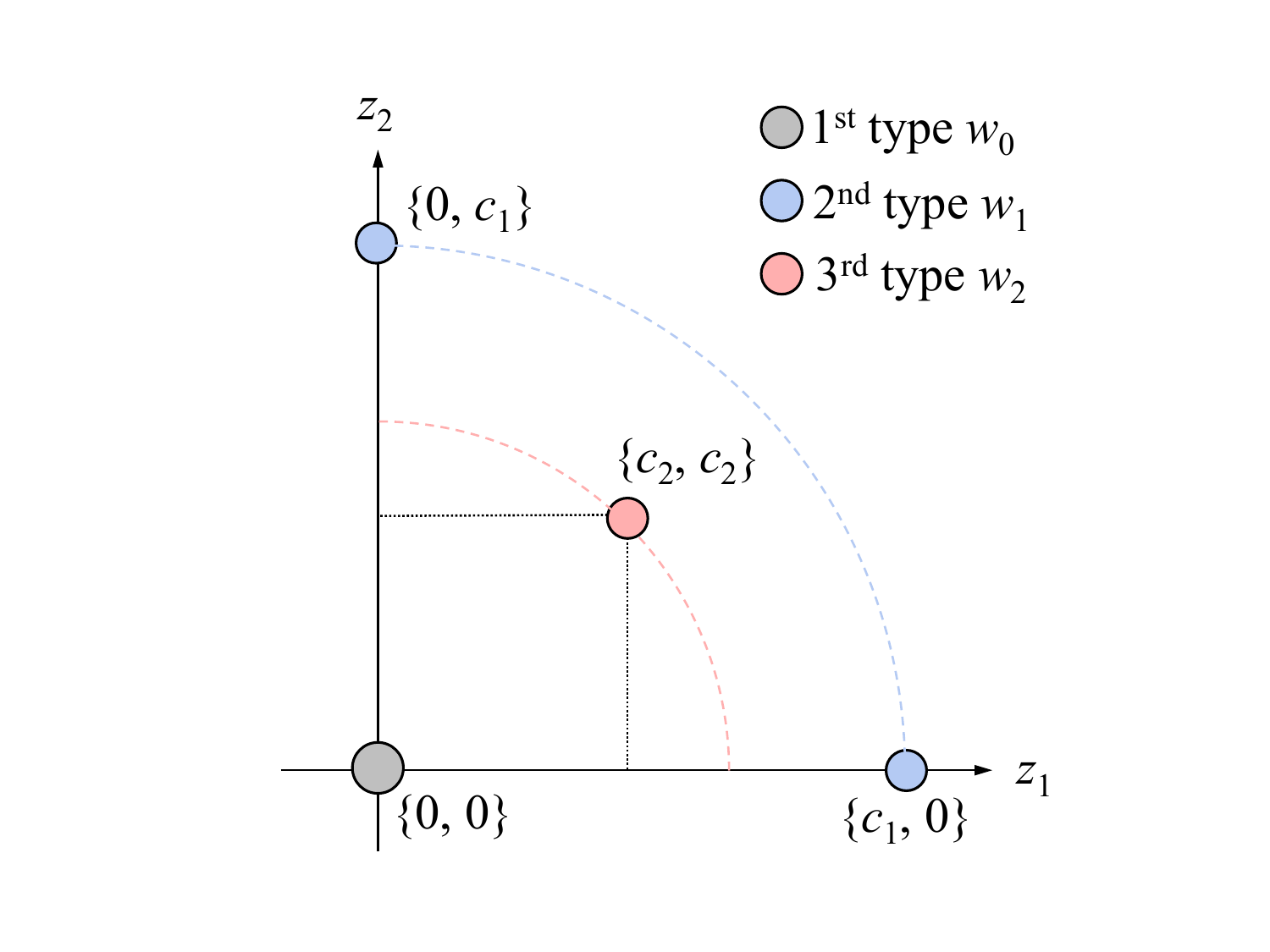}
		\caption{A quadrant in $2$-$\mathcal{D}$ space}
	\end{subfigure}
	\begin{subfigure}[b]{0.5\textwidth}
		\centering
		\includegraphics[width=\textwidth]{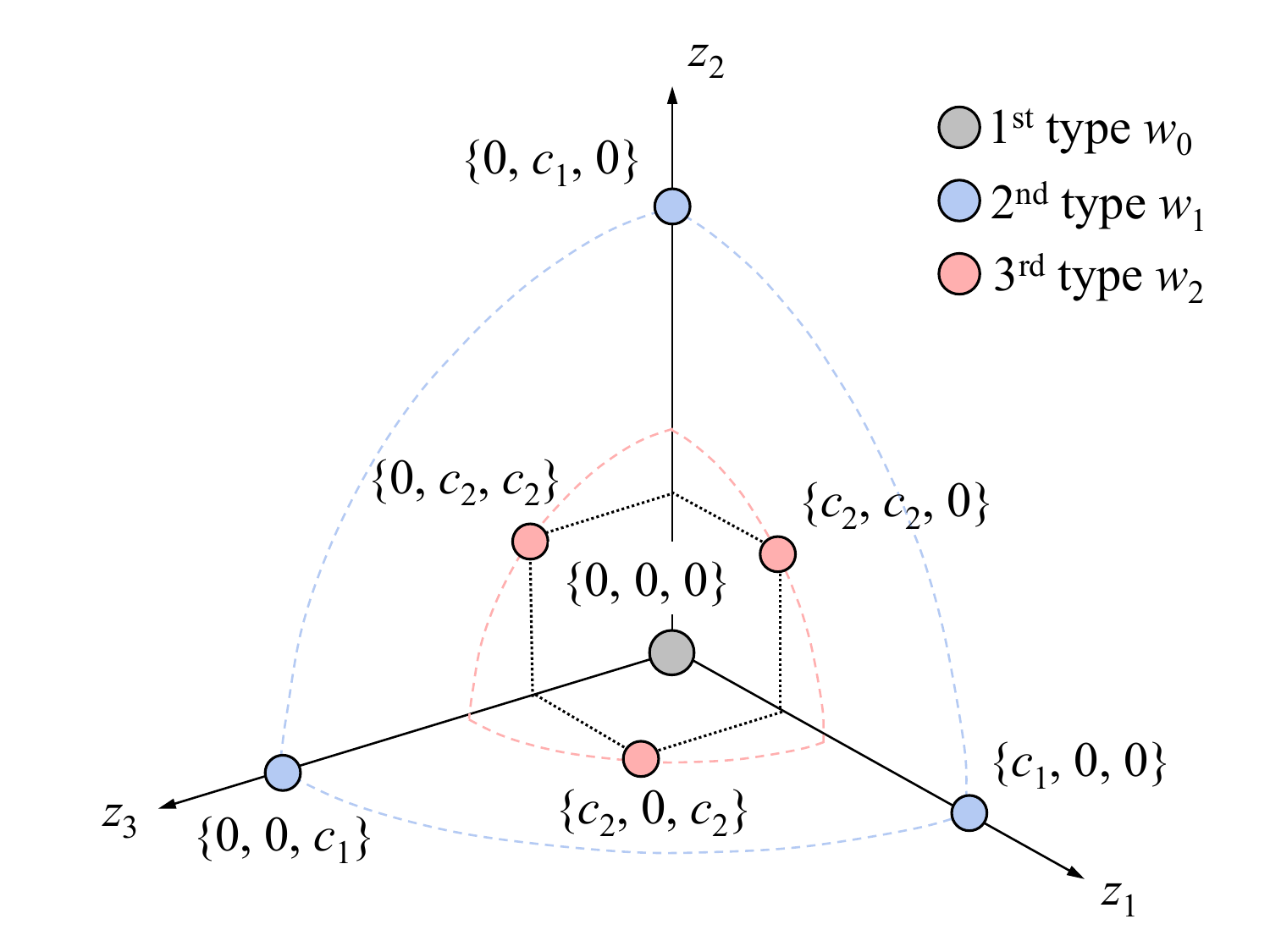}
		\caption{An octant in the $3$-$\mathcal{D}$ space}
	\end{subfigure}
	\caption{Schematic representation of the QPEM.}
	\label{fig:qpem}
\end{figure}

The QPEM fully symmetric sigma points can be seen in Fig.~\ref{fig:qpem}, together with the three types of points used: (i) The first type is a single point that lies at the central point $\mathbf{0}$, having weight $w_0$, (ii) the second type consists of $2n$ points on the orthogonal axes, with a distance of $c_1$ from the central point and with weights $w_1$, and (iii) the third type includes $2n(n-1)$ points obtained by permutations and sign changes of the coordinates ${{\left\{{{c}_{2}}, {{c}_{2}}, 0, \cdots , 0 \right\}}^{T}}$, with associated weights $w_2$. Overall, QPEM is thus described by and requires $2n^2+1$ points, $\mathbf{S}$ =[${{\mathbf{S}}_{0}}$, ${{\mathbf{S}}_{1}}$, $\cdots$, ${{\mathbf{S}}_{2n^2}}$]$_{n \times (2n^2+1)}$, with associated weights, $\mathbf{W}$ =[${{W}_{0}}$, ${{W}_{1}}$, $\cdots$, ${{W}_{2n^2}}$]$_{1 \times (2n^2+1)}$.

To match all the first four moments of the standard normal variables, the sigma points and corresponding weights of QPEM should satisfy:
\begin{equation}
	\label{eq:qpem_me}
	\left\{ \begin{array}{*{35}{lll}}
		E\left[ z_{i}^{2} \right]                 &= 2{{w}_{1}}c_{1}^{2}+4\left( n-1 \right){{w}_{2}}c_{2}^{2} &=1  \\[3pt]
		E\left[ z_{i}^{4} \right]                 &= 2{{w}_{1}}c_{1}^{4}+4\left( n-1 \right){{w}_{2}}c_{2}^{4} &=3  \\[3pt]
		E\left[ z_{i}^{2}z_{j}^{2} \right]        &= 4{{w}_{2}}c_{2}^{4}                                       &=1  \\ [3pt]
		\sum\limits_{l=0}^{2{{n}^{2}}}{{{W}_{l}}} &= {{w}_{0}}+2n{{w}_{1}}+2n\left( n-1 \right){{w}_{2}}       &=1
	\end{array} \right.
\end{equation}  
All odd-order moments of the standard normal variables are zero, which are, as mentioned, automatically captured by the fully symmetric sigma points. 

\subsection{Selection of sigma points and weights}
\label{sec:qpem_sigmapoints}
Five unknowns need to be identified so far, to fully determine the sigma points and associated weights, however, only four equations are available, as seen in Eq.~\eqref{eq:qpem_me}. As such, we assign a user-defined tuning parameter $r$ to ${{c}_{1}}$ in Eq.~\eqref{eq:qpem_me}, which indicates that the location of the second type of sigma points is predetermined at a distance $r$ from the central point, i.e.,  $c_1 = r$. This assignment reduces the number of unknowns from 5 to 4 ($c_2$, $w_0$, $w_1$, and $w_2$), matching the number of available equations in Eq.~\eqref{eq:qpem_me}. Thus, the analytical solution of the sigma points and weights can be now obtained as:
\begin{equation}
	\label{eq:qpem}
	\left\{\begin{array}{ll}
            {{c}_{0}} &=0\\ [0.3em]
		{{c}_{1}} &=r\\ [0.3em]
		{{c}_{2}} &=\left[\dfrac{r^2\left(n-1\right)}{r^2+n-4}\right]^{1/2}
	\end{array} \right. \qquad %
	\left\{ \begin{array}{ll}
            {{w}_{0}} &= 1-2nw_1-2n\left(n-1\right)w_2 \\ [0.3em]
		{{w}_{1}} &= \dfrac{4-n}{2r^4} \\ [0.3em]
		{{w}_{2}} &= \dfrac{1}{4}\left[\dfrac{r^2+n-4}{r^2\left(n-1\right)}\right]^2
	\end{array} \right. 
\end{equation}
where $r$ is a positive real number adhering to $r>\sqrt{2}$, as needed for QPEM stability. Regardless of the value of $r$, a set of sigma points obtained by Eq.~\eqref{eq:qpem} can fully capture up to fifth-order moments of the standard normal distribution in any dimensions. 

The parameter $r$ that controls the sigma points locations and weights has an important role for the estimation accuracy. Errors induced from higher than fifth-order moments could be reduced through the parameter $r$. For example, one could choose the parameter $r$ so as to satisfy the sixth-order marginal moment of the standard normal distribution, i.e., $E\left[z_i^6\right]$=15, as:
\begin{equation}
	\label{eq:qpem_lowdimension}
	\begin{array}{ll}
	    E\left[z_i^6\right] \approx \widehat{E}\left[ z_{i}^{6} \right] &=12{{w}_{1}}c_{1}^{6}+4\left(n-1\right){{w}_{2}}c_{2}^{6}\\ [0.3em]
	                     &=r^2(4-n)+(n-1)\left[\dfrac{r^2(n-1)}{r^2+n-4}\right]
	\end{array} \qquad %
\end{equation}
For spaces, with dimensions $n\geq4$, there is no real solution of $r$ that can satisfy Eq.~\eqref{eq:qpem_lowdimension}. The parameter $r$ can be then possibly determined by minimizing the squared error in the sixth-order marginal moment $e_{6}^2(r,n)$, as:

\begin{equation}
	\label{eq:qpem_highdimension1}
	{{e}_{6}^2}\left( r ,n \right)\equiv \left\{E\left[z_i^6\right]-r^2(4-n)-(n-1)\left[\dfrac{r^2(n-1)}{r^2+n-4}\right]\right\}^2
\end{equation}

\begin{equation}
	\label{eq:qpem_highdimension2}
	\dfrac{\partial }{\partial r}{{e}_{6}^2} \left( r ,n \right)=0 \,\,\,\,\,\,\,\,\,\, \text{where} \,\, r > \sqrt{2}
\end{equation}

It is interesting to note that regardless of the dimensions of the problem, $e_{6}^2(r,n)$ is always minimized when $r$ is equal to $\sqrt{3}$. However, this value of $\sqrt{3}$ is unfortunately not generally optimal in minimizing the overall error in estimating the output response function moments, since much higher-order moments can also play an important role for an accurate response moment estimation. In our study, we indeed empirically found that the $r$ value of $\sqrt{3}$ is quite small toward reducing the total estimation error.

\begin{figure}[!b]
	\begin{subfigure}[b]{0.49\textwidth}
		\centering
		\includegraphics[trim=0in 2.6in 0in 2.6in, clip=true, width=1.1\textwidth]{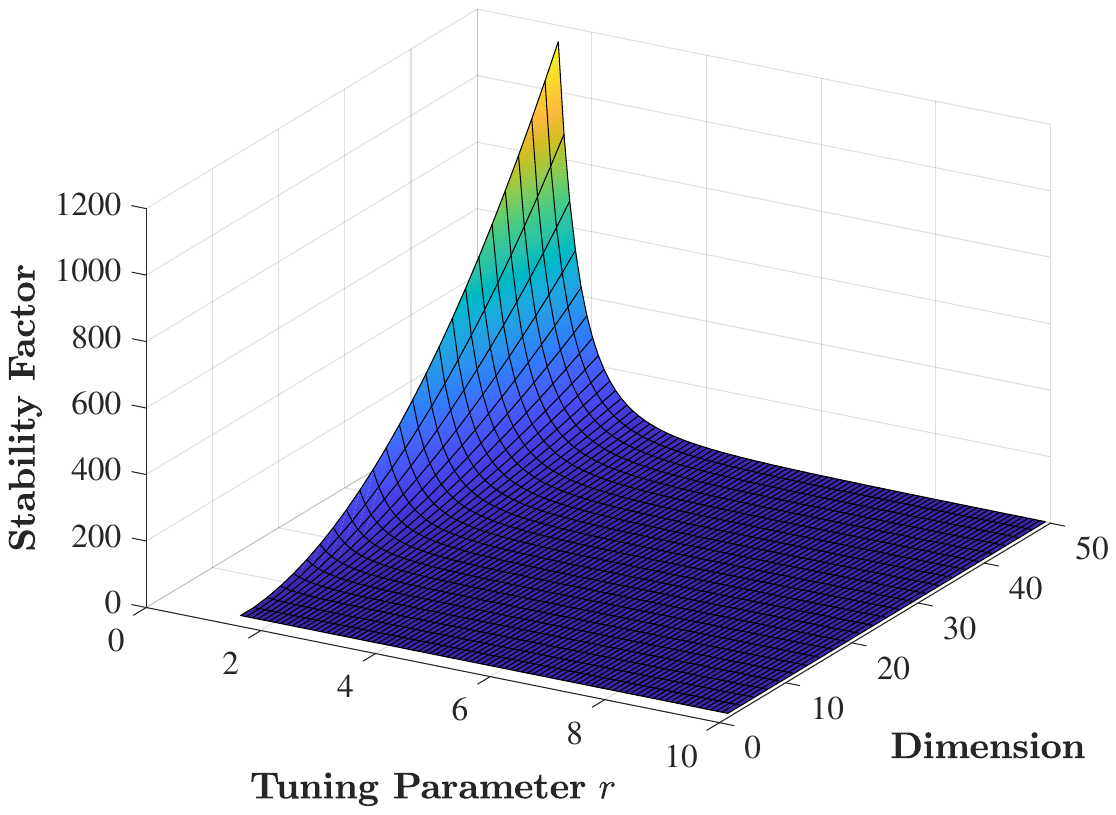}
		\caption{Stability factor depending on $r$ and dimension $n$}
	\end{subfigure}
	\begin{subfigure}[b]{0.49\textwidth}
		\centering
		\includegraphics[trim=0in 2.6in 0in 2.6in, clip=true, width=1.1\textwidth]{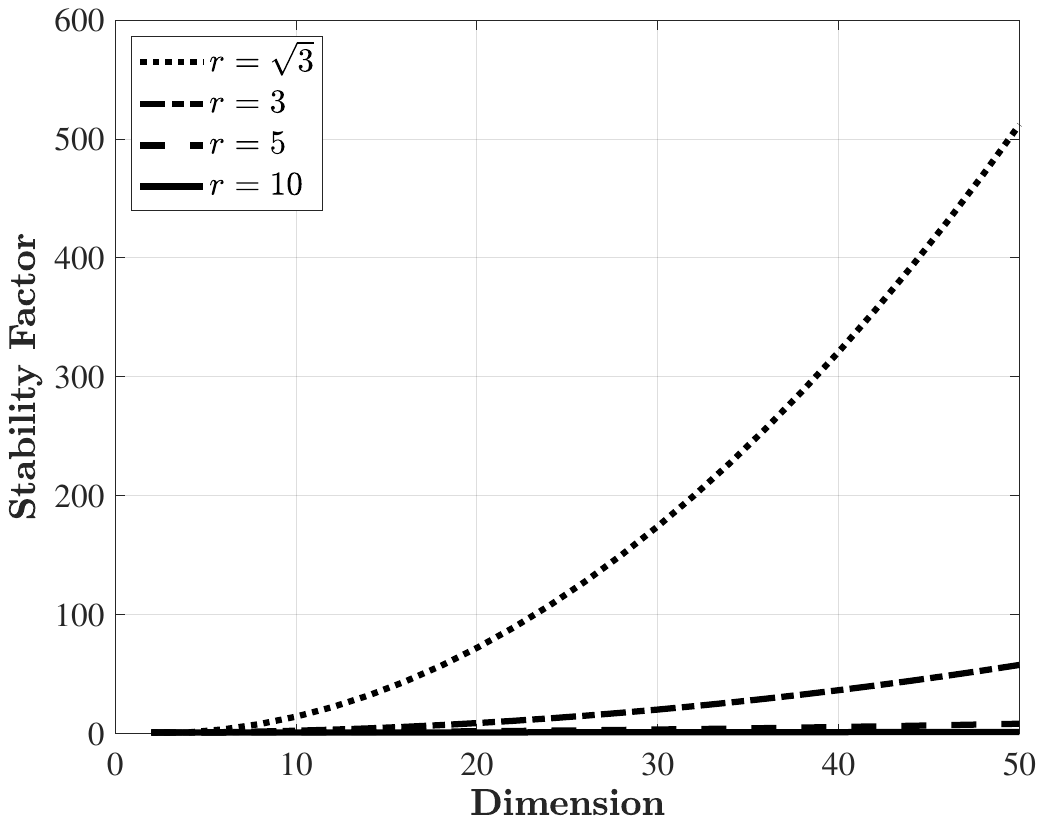}
		\caption{Comparison of stability factor for different $r$}
	\end{subfigure}
	\caption{Stability factor of the QPEM.}
	\label{fig:stability_factor}
\end{figure}

An important notion in numerical integration is the so-called stability factor, defined as the sum of the absolute values of the weights, $\sum_{i=0}^{N} |W_i|$, and has nothing to do with numerical stability. As explained in \cite{stroud1971approximate, davis2007methods,wu2006numerical}, the stability factor is generally related to the total integration error, and a value of 1, indicating only positive weights, is the most desirable for this purpose. However, so far, there is no known general method with non-exponential sampling complexity with increasing dimensions that can support sigma points with only positive weights. Fig.~\ref{fig:stability_factor} shows the stability factor of the QPEM with respect to the value of $r$ and dimensions $n$. As seen, with the value of $\sqrt{3}$ the stability factor increases considerably with increasing dimensions, in contrast to values of 3 and higher.

As seen in Fig.~\ref{fig:stability_factor}, with the stability factor getting toward the value of 1, the QPEM can nearly achieve having all positive weights when $r$ obtains larger values, e.g., $r \geq 5$. However, in such cases, the sigma points are scattered too far away from the mean, often creating numerical issues, due to the regions of the input space that are sampled, and can produce significant errors \cite{julier2002scaled}.  

Therefore, in this work, the default value of $r$ is empirically recommended as $r$=3 for the standard normal space, after extensive experimentation, which indicates that the second type of sigma points are fixed at a distance of three times the standard deviation ($3\sigma$) from the mean. As such, this value of $r$=3 works well in practice, as also shown extensively in the numerical examples section, and manages to offer a good balance between total error minimization and sigma points location in the input space. As also seen in Fig.~\ref{fig:stability_factor}(b), with $r$=3 the QPEM stability factor is very moderately increasing and remains below 100 up to 50 dimensions.  

\subsection{Incorporating higher order moments information}
\label{sec:qpem_scalingparameters}
The QPEM with $2n^2+1$ sigma points can fully capture, as explained, the first five moments of the standard normal distribution. However, the sigma points can still include partial higher-order even moments information, since all odd moments are already correctly computed to be 0, without any additional computational cost, by introducing two scaling parameters $\zeta$ and $\xi$.

Let us consider the Taylor series expansion of the nonlinear function $\mathbf{y}$=$\mathbf{M}(\mathbf{x})$ as: 
\begin{equation}
	\label{eq:function_taylor_expansion}
	\mathbf{y}=\sum\limits_{i=0}^{\infty }{\frac{1}{i!}{{\nabla }^{i}}\mathbf{M}\ }\bm{\delta}{{\mathbf{x}}^{i}}
\end{equation}
where $\bm{\delta}{\mathbf{x}}=\mathbf{x}-\bm{\mu}_{\mathbf{x}}$. The propagated mean and covariance can be then expressed by taking the expectation, as:
\begin{flalign}
    \label{eq:mean_taylor_expansion}
    \mathbf{\bar{y}} &=E\left[ \mathbf{y} \right] &&= \sum\limits_{i=0}^{\infty }{\frac{1}{i!}}{{\nabla }^{i}}\mathbf{M}\ E\left[ \bm{\delta }{{\mathbf{x}}^{i}} \right] \\
    \label{eq:covariance_taylor_expansion}    
    \mathbf{P_y} &=E\left[ {{\left( \mathbf{y}-\mathbf{\bar{y}} \right)}^{2}} \right]  &&=\sum\limits_{i=0}^{\infty }{\sum\limits_{j=0}^{\infty }{\frac{1}{i!j!}\left[ {{\nabla }^{i}}\mathbf{M}\otimes {{\nabla }^{j}}\mathbf{M} \right]\left\{ E\left[ \bm{\delta }{{\mathbf{x}}^{i}}\otimes \bm{\delta }{{\mathbf{x}}^{j}} \right]-E\left[ \bm{\delta }{{\mathbf{x}}^{i}} \right]\otimes E\left[\bm{\delta }{{\mathbf{x}}^{j}} \right] \right\}}}
\end{flalign}
with Eqs.~\eqref{eq:function_taylor_expansion}-\eqref{eq:covariance_taylor_expansion} following the notation as defined below, for ease and better readability of mathematical expressions:
\begin{align}
	&\nabla  &&=\left\{ \dfrac{\partial }{\partial {{x}_{1}}} \,\, \dfrac{\partial }{\partial {{x}_{2}}} \,\, \cdots \,\,  \dfrac{\partial }{\partial {{x}_{n}}} \right\} & \label{eq:notation1}\\
	&{{\nabla }^{k}}\mathbf{M} && =\underbrace{\nabla \otimes \nabla \otimes \cdots \otimes \nabla }_{k \,\, \text{times}}{{\left. \mathbf{M} \right|}_{\mathbf{x}=\bm{\mu}_{\mathbf{x}}}} \,\,\, =\left\{ \dfrac{\partial {{\nabla }^{k-1}}}{\partial {{x}_{1}}} \,\, \dfrac{\partial {{\nabla }^{k-1}}}{\partial {{x}_{2}}} \,\, \cdots  \,\, \dfrac{\partial {{\nabla }^{k-1}}}{\partial {{x}_{n}}} \right\} {{\left. \mathbf{M} \right|}_{\mathbf{x}=\bm{\mu}_{\mathbf{x}}}}  \\[11pt]
	&\bm{\delta \mathbf{x}} && ={{\left\{ \delta {{x}_{1}} \,\, \delta {{x}_{2}} \,\, \cdots  \,\,  \delta {{x}_{n}} \right\}}^{T}} &
\end{align}

\begin{align}
	&\bm{\delta }{{\mathbf{x}}^{k}} && =\underbrace{\bm{\delta \mathbf{x}}\otimes \bm{\delta \mathbf{x}}\otimes \cdots \otimes \bm{\delta \mathbf{x}}}_{k \,\, \text{times}} \qquad \quad ={{\left\{ \delta {{x}_{1}}\bm{\delta \mathbf{x}}^{{k-1}^{T}} \,\, \delta {{x}_{2}}\bm{\delta \mathbf{x}}^{{k-1}^{T}} \,\, \cdots \,\,  \delta {{x}_{n}}\bm{\delta \mathbf{x}}^{{k-1}^{T}} \right\}}^{T}} &\\
    &{{\left( {{\nabla }^{k}}\mathbf{g} \right)}^{\alpha }}&& =\underbrace{{{\nabla }^{k}}\mathbf{g}\otimes {{\nabla }^{k}}\mathbf{g}\otimes \cdots \otimes {{\nabla }^{k}}\mathbf{g}}_{\alpha \ \text{times}} &\\
    &E{{\left[ \bm{\delta }{{\mathbf{x}}^{k}} \right]}^{\alpha }} && =\underbrace{E\left[ \bm{\delta }{{\mathbf{x}}^{k}} \right]\otimes E\left[ \bm{\delta }{{\mathbf{x}}^{k}} \right]\otimes \cdots \otimes E\left[ \bm{\delta }{{\mathbf{x}}^{k}} \right]}_{\alpha \ \text{times}}  & \label{eq:notation2}
\end{align}
Following this notation, the full expressions for the first four moments of the output response function can be seen in the \ref{Appendix-A}, and the QPEM can satisfy all terms of the first five input moments in these expressions.

The scaling parameters $\zeta$ and $\xi$ are introduced to partially incorporate the sixth-order input moments information in the estimation of the third-order output moments, and the eighth-order input moments information in the estimation of the fourth-order output moments, respectively. Following Eq.~\eqref{eq:moments_y}, the third-order moment of the response function is obtained as:
\begin{equation}
	\label{eq:qpem_skewness}
	E \left[\left(\mathbf{y} -\bar{\mathbf{y}}\right)^{3} \right] \approx \sum_{i=0}^{2n^{2}} W_i \left(\mathbf{Y}_i - \bar{\mathbf{y}}\right)^{3}
\end{equation}
Similar to Eq.~\eqref{eq:covariance_taylor_expansion}, the third-order moments at the left-hand side of Eq.~\eqref{eq:qpem_skewness} can be expanded through Taylor series, neglecting here the odd-order input moments that are zero for normal distribution, as follows (see \ref{Appendix-A} for further details and full expressions):
\begin{equation}
	\label{eq:qpem_taylor_skewness}
	\begin{array}{ll}
		E\left[\left(\mathbf{y} -\bar{\mathbf{y}}\right)^{3} \right] &=\dfrac{1}{2}\left[ {{\left( \nabla \mathbf{M} \right)}^{2}}\otimes {{\nabla }^{2}}\mathbf{M} \right]\left\{ 3E\left[ \bm{\delta } {{\mathbf{x}}^{4}} \right]-3E{{\left[ \bm{\delta } {{\mathbf{x}}^{2}} \right]}^{2}} \right\} \\ [11pt]
		{}&+\underbrace{\frac{1}{8}{\left[{\left( {{\nabla }^{2}}\mathbf{M} \right)}^{3}\right]}\left\{%
			\begin{array}{l}
				E\left[ \bm{\delta } {{\mathbf{x}}^{2}}\otimes \bm{\delta } {{\mathbf{x}}^{2}} \otimes \bm{\delta } {{\mathbf{x}}^{2}} \right] \\
				-3E\left[ \bm{\delta } {{\mathbf{x}}^{2}} \otimes \bm{\delta } {{\mathbf{x}}^{2}} \right]\otimes E\left[ \bm{\delta } {{\mathbf{x}}^{2}} \right]+2E{{\left[ \bm{\delta } {{\mathbf{x}}^{2}} \right]}^{3}}\\ 
			\end{array}%
			\right\}}_{\equiv {{\mathbf{A}}_{6}}}+\cdots
	\end{array}
\end{equation}
where $\mathbf{A}_6$ is one of the sixth-order terms in the Taylor series expansion of the third-order output moment and the QPEM cannot capture the vector of the sixth-order input moments $E[\bm{\delta} {{\mathbf{x}}^{2}}\otimes \bm{\delta} {{\mathbf{x}}^{2}} \otimes \bm{\delta} {{\mathbf{x}}^{2}}]$. Expanding now the contribution of the central sigma point ($i$=0) at the right-hand side of Eq.~\eqref{eq:qpem_skewness} results in:
\begin{equation}
	\label{eq:qpem_taylor_central_skewness}
	\left(\mathbf{Y}_0-\bar{\mathbf{y}} \right)^{3} = -\dfrac{1}{8}\left[\left(\nabla^{2} \mathbf{M}\right)^{3}\right] \left\{E\left[\bm{\delta}\mathbf{x}^2\right]^3\right\} + \cdots
\end{equation}
Comparing Eqs.~\eqref{eq:qpem_taylor_skewness} and \eqref{eq:qpem_taylor_central_skewness}, it is observed that the term $E[\bm{\delta} \mathbf{x}^2]^3$ is contained in both. Thus, the term $E[\bm{\delta} \mathbf{x}^2]^3$ in $\mathbf{A}_6$ of Eq.~\eqref{eq:qpem_taylor_skewness} could be controlled through the central sigma point contribution and relevant scaled third-order moment expression, with the scaling parameter $\zeta$ obtained as:
\begin{equation}
	\label{eq:weighted_skewness}
	E\left[ {{\left( \mathbf{y}-\mathbf{\bar{y}} \right)}^{3}} \right]\approx \underbrace{\sum\limits_{i=0}^{2{{n}^{2}}}{{{W}_{i}}{{\left( {{\mathbf{Y}}_{i}}-\mathbf{\bar{y}} \right)}^{3}}}}_{\begin{array}{cc}
    \text{\scriptsize Accurate up to the 4th-order}\\[-5pt]
    \text{\scriptsize and all higher odd input moments}
    \end{array}}+\underbrace{\zeta {{\left( {{\mathbf{Y}}_{0}}-\mathbf{\bar{y}} \right)}^{3}}}_{\text{Added for the partial 6th-order terms}}
\end{equation}
where the summation term provides the fourth- and all higher odd-order accuracy and the additional term with the scaling parameter $\zeta$ incorporates the partial sixth-order moment effect. The error in $\mathbf{A}_6$ can now be controlled without additional sigma points, as:
\begin{equation}
    \label{eq:determine_zeta}
	\begin{array}{ll}
		\Delta {{\mathbf{A}}_{6}}&={{\mathbf{A}}_{6}}-\zeta {{\left( {{\mathbf{Y}}_{0}}-\mathbf{\bar{y}} \right)}^{3}} \\ [11pt]
		{}& =\dfrac{1}{8}{\left[{\left({{\nabla }^{2}}\mathbf{M} \right)}^{3}\right]}\left\{%
		\begin{array}{l}
			E\left[ \bm{\delta } {{\mathbf{x}}^{2}}\otimes \bm{\delta } {{\mathbf{x}}^{2}} \otimes \bm{\delta } {{\mathbf{x}}^{2}} \right] \\
			\qquad \qquad-3E\left[ \bm{\delta } {{\mathbf{x}}^{2}} \otimes \bm{\delta } {{\mathbf{x}}^{2}} \right]\otimes E\left[ \bm{\delta } {{\mathbf{x}}^{2}} \right]+(2+\zeta)E{{\left[ \bm{\delta } {{\mathbf{x}}^{2}} \right]}^{3}}\\ 
		\end{array}%
		\right\}  
	\end{array}
\end{equation}
In the case of a normal distribution, for example, the optimal value of $\zeta$ is -8. This is because $E[\bm{\delta}\mathbf{x}^2 \otimes \bm{\delta}\mathbf{x}^2 \otimes \bm{\delta}\mathbf{x}^2]$ $\approx$ $15E[\bm{\delta}\mathbf{x}^2]^3$ and $E[\bm{\delta}\mathbf{x}^2 \otimes \bm{\delta}\mathbf{x}^2]$ $\approx$ $3E[\bm{\delta}\mathbf{x}^2]^2$. The approximate sign is used here because not all elements within the vectors adhere to the aforementioned relationship in a general multivariate case. 

It should also be noted here that this scaling process and parameter $\zeta$ do not provide guarantees for a reduction in the third-order output moment error in all cases, since this scaling also induces errors to some of the sixth-order and higher-order even input moments. The logic is that all these terms cannot be controlled by QPEM anyways, so with this scaling QPEM can at least control some of the sixth-order terms. If, however, some of the higher-order terms are particularly influential, in combination with certain functions, for the output estimation, then this scaling may not be beneficial. Nonetheless, this is not something common and we have only observed this in two instances, both for skewness estimation, from all of the numerical examples in this paper. These same theoretical statements also apply accordingly for the fourth-order output moments, as described subsequently.

Along the same lines, the fourth-order output moment can be estimated by the QPEM as:
\begin{equation}
	\label{eq:qpem_kurtosis}
	E \left[\left(\mathbf{y} -\bar{\mathbf{y}}\right)^{4} \right] \approx \sum_{i=0}^{2n^{2}} W_i \left(\mathbf{Y}_i - \bar{\mathbf{y}}\right)^{4}
\end{equation}
and the following expressions are obtained by Taylor series expansions of the left-hand side of Eq.~\eqref{eq:qpem_kurtosis} and the central sigma points contribution $\left(\mathbf{Y}_0 - \bar{\mathbf{y}}\right)^{4}$  (see \ref{Appendix-A}):
\begin{equation}
	\begin{array}{ll}
		\label{eq:qpem_taylor_kurtosis}
		E\left[ {{\left( \mathbf{y}-\mathbf{\bar{y}} \right)}^{4}} \right]&=\underbrace{\cdots }_{\text{Up to 7th order}}\\
		&+\underbrace{\frac{1}{16}{\left[{\left( {{\nabla }^{2}}\mathbf{M} \right)}^{4}\right]}\left\{%
			\begin{array}{l}
				E\left[ \bm{\delta} {{\mathbf{x}}^{2}}\otimes \bm{\delta} {{\mathbf{x}}^{2}}\otimes \bm{\delta} {{\mathbf{x}}^{2}}\otimes \bm{\delta} {{\mathbf{x}}^{2}} \right]\\
				\quad  -4E\left[ \bm{\delta} {{\mathbf{x}}^{2}} \otimes \bm{\delta} {{\mathbf{x}}^{2}}\otimes \bm{\delta} {{\mathbf{x}}^{2}} \right]\otimes E\left[ \bm{\delta} {{\mathbf{x}}^{2}} \right] \\ 
				\qquad\quad \,\,\,\,\,\,\,+6E\left[ \bm{\delta} {{\mathbf{x}}^{2}} \otimes \bm{\delta} {{\mathbf{x}}^{2}} \right]\otimes E{{\left[ \bm{\delta} {{\mathbf{x}}^{2}} \right]}^{2}}\\
				\qquad\qquad\qquad\qquad\qquad \,\,\,\,\,\,\,\,\,\, -3E{{\left[ \bm{\delta}{{\mathbf{x}}^{2}} \right]}^{4}} \\ 
			\end{array}%
			\right\}}_{\equiv {{\mathbf{A}}_{8}}}+\cdots
	\end{array}
\end{equation}

\begin{equation}
	\label{eq:qpem_taylor_central_kurtosis}
	{{\left( {{\mathbf{Y}}_{0}}-\mathbf{\bar{y}} \right)}^{4}}=\frac{1}{16}{\left[{\left( {{\nabla }^{2}}\mathbf{M} \right)}^{4}\right]}\left\{E{{\left[ \bm{\delta} {{\mathbf{x}}^{2}} \right]}^{4}}\right\}+\cdots
\end{equation}
$\mathbf{A}_8$ is one of the eighth-order terms of the Taylor series expansion of the fourth-order output moments and the QPEM cannot capture the vectors of the eighth- and sixth-order input moments, $E[ \bm{\delta} {{\mathbf{x}}^{2}}\otimes \bm{\delta} {{\mathbf{x}}^{2}}\otimes \bm{\delta} {{\mathbf{x}}^{2}}\otimes \bm{\delta} {{\mathbf{x}}^{2}}]$ and $E[ \bm{\delta} {{\mathbf{x}}^{2}} \otimes \bm{\delta} {{\mathbf{x}}^{2}}\otimes \bm{\delta} {{\mathbf{x}}^{2}} ]$. The last term, $E[\bm{\delta}\mathbf{x}^2]^4$, of $\mathbf{A}_8$ is, however, included in both Eqs.~\eqref{eq:qpem_taylor_kurtosis} and \eqref{eq:qpem_taylor_central_kurtosis}, and the term could be controlled through the scaling parameter $\xi$, to incorporate the partial eighth-order input moment information for the estimation of the fourth-order output moments, similarly as before:
\begin{equation}
	E\left[ {{\left( \mathbf{y}-\mathbf{\bar{y}} \right)}^{4}} \right]\approx \underbrace{\sum\limits_{i=0}^{2{{n}^{2}}}{{{W}_{i}}{{\left( {{\mathbf{Y}}_{i}}-\mathbf{\bar{y}} \right)}^{4}}}}_{\begin{array}{cc}
    \text{\scriptsize Accurate up to the 4th-order}\\[-5pt]
    \text{\scriptsize and all higher odd input moments}
    \end{array}}+\underbrace{\xi {{\left( {{\mathbf{Y}}_{0}}-\mathbf{\bar{y}} \right)}^{4}}}_{\text{Added for the partial 8th-order terms}}
\end{equation}

\begin{equation}
    \label{eq:determine_xi}
	\begin{array}{ll}
		\Delta {{\mathbf{A}}_{8}} &={{\mathbf{A}}_{8}}-\xi {{\left( {{\mathbf{Y}}_{0}}-\mathbf{\bar{y}} \right)}^{4}} \\ 
		{}&=\dfrac{1}{16}{\left[{\left( {{\nabla }^{2}}\mathbf{M} \right)}^{4}\right]}\left\{%
		\begin{array}{l}
			E\left[ \bm{\delta} {{\mathbf{x}}^{2}}\otimes \bm{\delta} {{\mathbf{x}}^{2}}\otimes \bm{\delta} {{\mathbf{x}}^{2}}\otimes \bm{\delta} {{\mathbf{x}}^{2}} \right]\\
			\qquad \qquad-4E\left[ \bm{\delta} {{\mathbf{x}}^{2}} \otimes \bm{\delta} {{\mathbf{x}}^{2}}\otimes \bm{\delta} {{\mathbf{x}}^{2}} \right]\otimes E\left[ \bm{\delta} {{\mathbf{x}}^{2}} \right] \\ 
			\qquad \qquad \qquad  \qquad+6E\left[ \bm{\delta} {{\mathbf{x}}^{2}} \otimes \bm{\delta} {{\mathbf{x}}^{2}} \right]\otimes E{{\left[ \bm{\delta} {{\mathbf{x}}^{2}} \right]}^{2}}\\
			\qquad \qquad \qquad \qquad \qquad \qquad \qquad -(3+\xi)E{{\left[ \bm{\delta}{{\mathbf{x}}^{2}} \right]}^{4}} \\ 
		\end{array}%
		\right\} 
	\end{array}
\end{equation}
For the special case of a normal input distribution, the optimal value of $\xi$ is 60 because $E[\bm{\delta}\mathbf{x}^2 \otimes \bm{\delta}\mathbf{x}^2 \otimes \bm{\delta}\mathbf{x}^2 \otimes \bm{\delta}\mathbf{x}^2]$ $\approx$ $105E[\bm{\delta}\mathbf{x}^2]^4$, $E[\bm{\delta}\mathbf{x}^2 \otimes \bm{\delta}\mathbf{x}^2 \otimes \bm{\delta}\mathbf{x}^2]$ $\approx$ $15E[\bm{\delta}\mathbf{x}^2]^3$ and $E[\bm{\delta}\mathbf{x}^2 \otimes \bm{\delta}\mathbf{x}^2]$ $\approx$ $3E[\bm{\delta}\mathbf{x}^2]^2$. 
\par
All parameters, the sigma points locations and associated weights for the QPEM can now be summarized in Table~\ref{table:qpem_points_weights}. The method can be straightforwardly applied now and the first four output moments are obtained as:
\begin{align}
    \label{eq:qpem_mean_equation}
    E\left[ \mathbf{y} \right]\equiv \mathbf{\bar{y}} & =\displaystyle \sum\limits_{i=0}^{2{{n}^{2}}}{W_{i}^{\left( 1 \right)}{{\mathbf{Y}}_{i}}} \\ \label{eq:qpem_moments_equation}
    E\left[ {{\left( \mathbf{y}-\mathbf{\bar{y}} \right)}^{k}} \right] & =\displaystyle \sum\limits_{i=0}^{2{{n}^{2}}}{W_{i}^{\left( k \right)}\left(\mathbf{Y}_i-\overline{\mathbf{y}}\right)^{k}} \,\,\,\,\,\,\,\,\,\, \text{for}\,\,\, k=2,3,4 
\end{align}

\begin{table}[!t]
	\caption{Sigma points locations and weights for the QPEM}\label{table:qpem_points_weights}
	\begin{tabular*}{\textwidth}[t]{@{}lll@{}}
		\toprule
		\begin{tabular}[t]{l}Weights \end{tabular}&  & \begin{tabular}[t]{|l} 
			$W_{0}^{(1)}=W_{0}^{(2)}=w_0$   \quad $W_{0}^{(3)}=w_0+\zeta$  \quad $W_{0}^{(4)}=w_0+\xi$\\
			$W_i^{(k)\,\,a} = \begin{cases}
				w_1, & \text{for $i=1,\cdots, 2n$} \\ 
				w_2, & \text{for $i=2n+1,\cdots, 2n^2$} 
			\end{cases}$ \end{tabular} \\
            && \begin{tabular}[t]{|l} where, $\left\{ \begin{array}{ll}
            {{w}_{0}} &= 1-2nw_1-2n\left(n-1\right)w_2 \\
		{{w}_{1}} &= \frac{4-n}{2r^4} \\
		{{w}_{2}} &= \frac{1}{4}\left[\frac{r^2+n-4}{r^2\left(n-1\right)}\right]^2
	\end{array} \right.$  \end{tabular}\\
		\midrule
		\multirow{2}{*}{\begin{tabular}[t]{l} Sigma\\ Point \end{tabular}} & 1st type & \begin{tabular}[t]{|l} 
			$\mathbf{S}_0=\mathbf{0}_{n \times 1}$
		\end{tabular}\\
		\begin{tabular}[t]{l}	\end{tabular}							& 2nd type & \begin{tabular}[t]{|lll}
			$\mathbf{S}_1=\{c_1, 0,\cdots, 0\}^T$& $\cdots$ & $\mathbf{S}_n=\{0,\cdots, 0, c_1\}^T$\\
			$\mathbf{S}_{n+1}=\{-c_1, 0,\cdots, 0\}^T$& $\cdots$ & $\mathbf{S}_{2n}=\{0,\cdots, 0, -c_1\}^T$\end{tabular} \\
		\begin{tabular}[t]{l}	\end{tabular}							& 3rd type$^{b}$ &  \begin{tabular}[t]{|lll}
			$\mathbf{S}_{2n+1}=\{c_2, c_2, 0, \cdots, 0\}^T$ & $\cdots$ & $\mathbf{S}_{2n^2-3}=\{0,\cdots, 0, c_2, c_2\}^T$\\
			$\mathbf{S}_{2n+2}=\{-c_2, c_2, 0,\cdots, 0\}^T$& $\cdots$ & $\mathbf{S}_{2n^2-2}=\{0,\cdots, 0, -c_2, c_2\}^T$\\
			$\mathbf{S}_{2n+3}=\{c_2, -c_2, 0,\cdots, 0\}^T$& $\cdots$ & $\mathbf{S}_{2n^2-1}=\{0,\cdots, 0, c_2, -c_2\}^T$\\
			$\mathbf{S}_{2n+4}=\{-c_2, -c_2, 0,\cdots, 0\}^T$& $\cdots$ & $\mathbf{S}_{2n^2}=\{0,\cdots, 0, -c_2, -c_2\}^T$\end{tabular} \\
            &&\begin{tabular}[t]{|l} where, $\left\{\begin{array}{ll}
            {{c}_{0}} &=0\\
		{{c}_{1}} &=r\\
		{{c}_{2}} &=\left[\frac{r^2\left(n-1\right)}{r^2+n-4}\right]^{1/2}
	\end{array} \right.$ \end{tabular}\\
            \midrule
            \multirow{2}{*}{\begin{tabular}[t]{l} Parameters \end{tabular}} & \qquad $r$ & \begin{tabular}[t]{|l} Default value of 3 in this work \end{tabular}\\
		\begin{tabular}[t]{l}	\end{tabular}	&  \qquad $\zeta$ & \begin{tabular}[t]{|l} Determined by Eq.~\eqref{eq:determine_zeta} or -8 for normal input distribution \end{tabular}\\
		\begin{tabular}[t]{l}	\end{tabular}	&  \qquad $\xi$ &  \begin{tabular}[t]{|l} Determined by Eq.~\eqref{eq:determine_xi} or 60 for normal input distribution \end{tabular}\\
		\bottomrule
		\multicolumn{3}{l}{\begin{tabular}{l} $^a$ Associated weights for $k^{th}$ moments of the response function\\$^b$ Obtained by permutations and by sign alterations based on $\mathbf{S}_{2n+1}$ \end{tabular}}\\ 
	\end{tabular*}
\end{table}
The QPEM with nonzero $\zeta$ and $\xi$ is called \textit{scaled} QPEM in this work. As an important note here, for further clarity, although Taylor series expansion has been used to prove the logic behind QPEM, the scaling parameters and their values, the actual implementation of the method just simply relies on Table~\ref{table:qpem_points_weights} and Eqs.~\eqref{eq:qpem_mean_equation}-\eqref{eq:qpem_moments_equation}, and can also work for non-differentiable, non-smooth, and discontinuous functions, as other PEMs, e.g., \cite{papakonstantinou2022scaled,rosenblueth1981two, julier2004unscented}.

\section{Numerical examples}
In this section, several numerical examples are studied in order to evaluate the QPEM performance in estimating the first four output moments, i.e., mean, standard deviation, skewness coefficient, and kurtosis coefficient. Overall, three types of problems are examined, without loss of generality; (i) a polynomial function with varying dimensions from 5 to 100-$\mathcal{D}$, (ii) two structural analysis examples, and (iii) three Finite Element Analysis (FEA) examples, involving spatial stochastic fields represented by the Karhunen-Loève (KL) expansion.

When an exact solution to the problem is not available, the MC method with $10^6$ samples has been used as a reference solution. The QPEM is compared here against variance reduction techniques, and in particular, Latin Hypercube Sampling (LHS) and Sobol's Quasi-Monte Carlo (QMC), a sparse grid quadrature technique, i.e., Smolyak's Gauss-Hermite quadrature with three quadrature points in a single dimension (SGH3), and HPEM, as one of the most well-known and popular PEMs.

\subsection{Polynomial function with varying dimensions}
The first benchmark example is a second-order polynomial function with varying dimensionality \cite{shields2016generalization}:
\begin{equation}
	\label{eq:schwefel_eq}
	y=\sum_{i=1}^{n} \left(\sum_{j=1}^{i} x_j \right)^2
\end{equation}
Here, $\mathbf{x}=\left\{x_1, x_2, \cdots, x_n \right\}^T$ represents an $n$-$\mathcal{D}$ random variable, where $x_i \sim N(5, 1)$. Dimensions ranging from 5 to 100 have been considered.

\begin{table}[!tb]
	\centering
	\caption{Sample size used for each method according to space dimensions}\label{table:sample_size}
	\begin{tabular*}{\textwidth}{p{0.12\textwidth}p{0.06\textwidth}p{0.06\textwidth}p{0.06\textwidth}p{0.06\textwidth}p{0.06\textwidth}p{0.06\textwidth}p{0.06\textwidth}p{0.06\textwidth}p{0.06\textwidth}p{0.06\textwidth}}
		\toprule
		Dimension& 5  & 10  & 15  & 20  & 30   & 40   & 50   & 60   & 70   & 100   \\\midrule
		LHS       & 51 & 201 & 451 & 801 & 1,801 & 3,201 & 5,001 & 7,201 & 9,801 & 20,001 \\
		QMC      & 51 & 201 & 451 & 801 & 1,801 & 3,201 & 5,001 & 7,201 & 9,801 & 20,001 \\
		SGH3      & 61 & 221 & 481 & 841 & 1,861 & 3,281 & 5,101 & 7,321 & 9,941 & 20,201 \\
		HPEM        & 11 & 21  & 31  & 41  & 61   & 81   & 101  & 121  & 141  & 201   \\
		QPEM      & 51 & 201 & 451 & 801 & 1,801 & 3,201 & 5,001 & 7,201 & 9,801 & 20,001 \\ \bottomrule
	\end{tabular*}
\end{table}

The QPEM has been applied with $r$=3, as the default value reported in Section~\ref{sec:qpem_sigmapoints}, and $\zeta$=-8 and $\xi$=60 for its scaled version, as explained in Section~\ref{sec:qpem_scalingparameters}. For all methods the number of samples is kept similar to the one needed by QPEM, except HPEM which has a lower computational cost. The sample sizes for all methods according to the dimensions are summarized in Table~\ref{table:sample_size}. The relative errors of each method with respect to the exact solutions in this case are provided in Fig.~\ref{fig:ex1}, with the axes of the relative errors in logarithmic scale. The zero error values for mean and standard deviation imply that the methods could compute the exact solution in this case. Some observations from the results in Fig.~\ref{fig:ex1} follow:

\begin{figure}[!tb]
	\centering
	\begin{subfigure}[tb]{0.45\textwidth}
		\includegraphics[width=\textwidth]{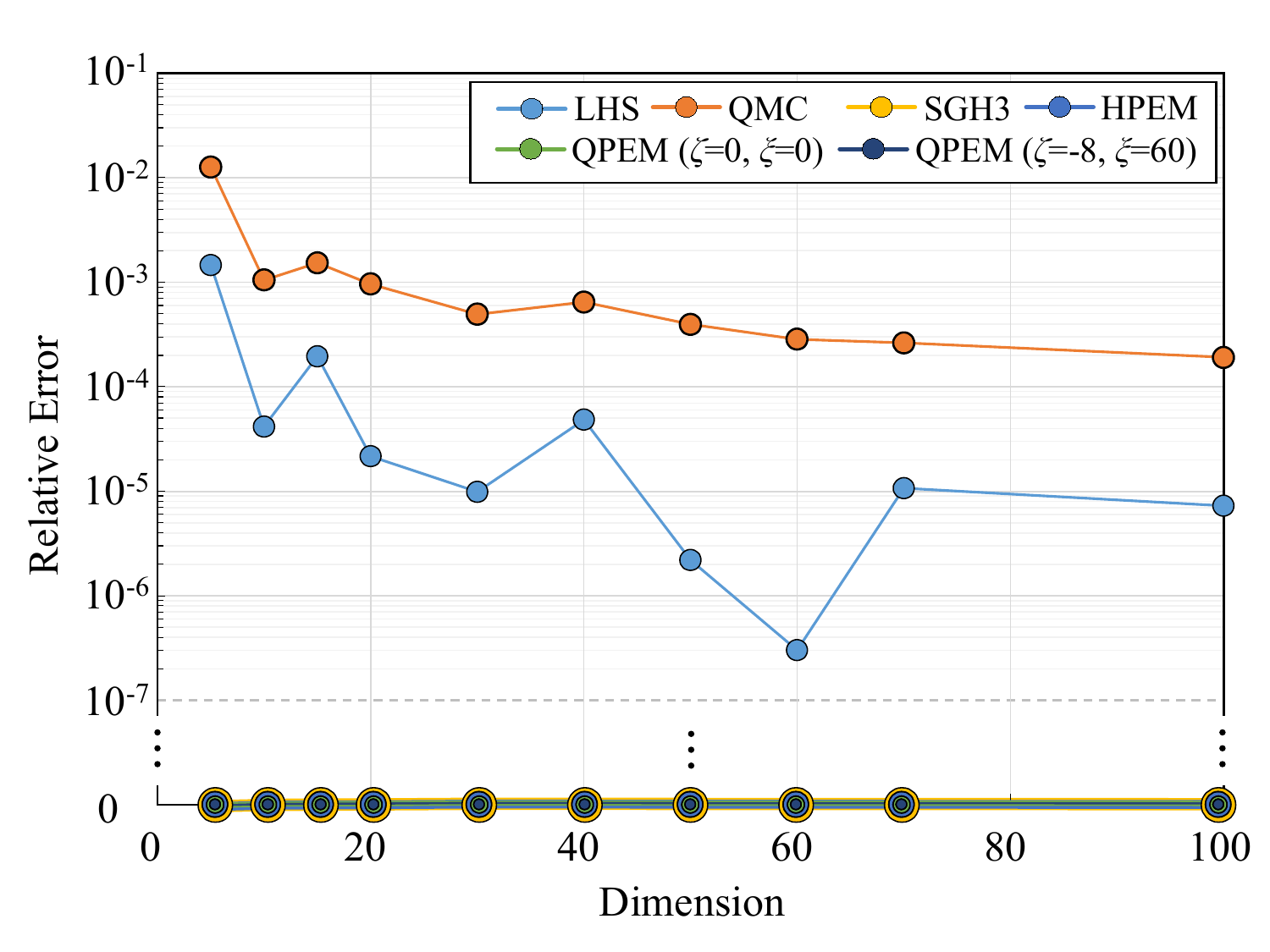}
		\caption{Mean}
	\end{subfigure}
	\begin{subfigure}[tb]{0.45\textwidth}
		\includegraphics[width=\textwidth]{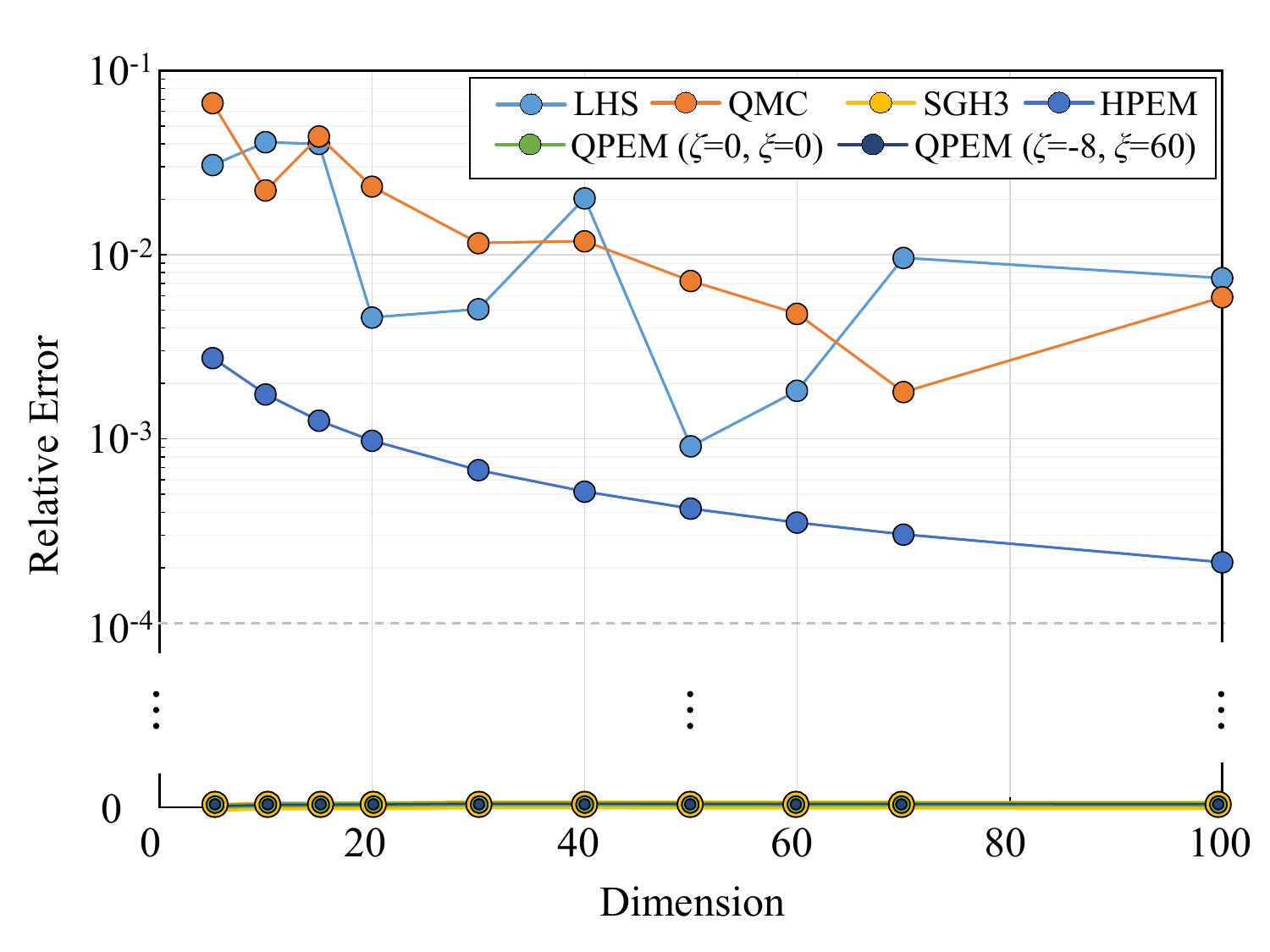}
		\caption{Standard deviation}
	\end{subfigure}
	\begin{subfigure}[tb]{0.45\textwidth}
		\includegraphics[width=\textwidth]{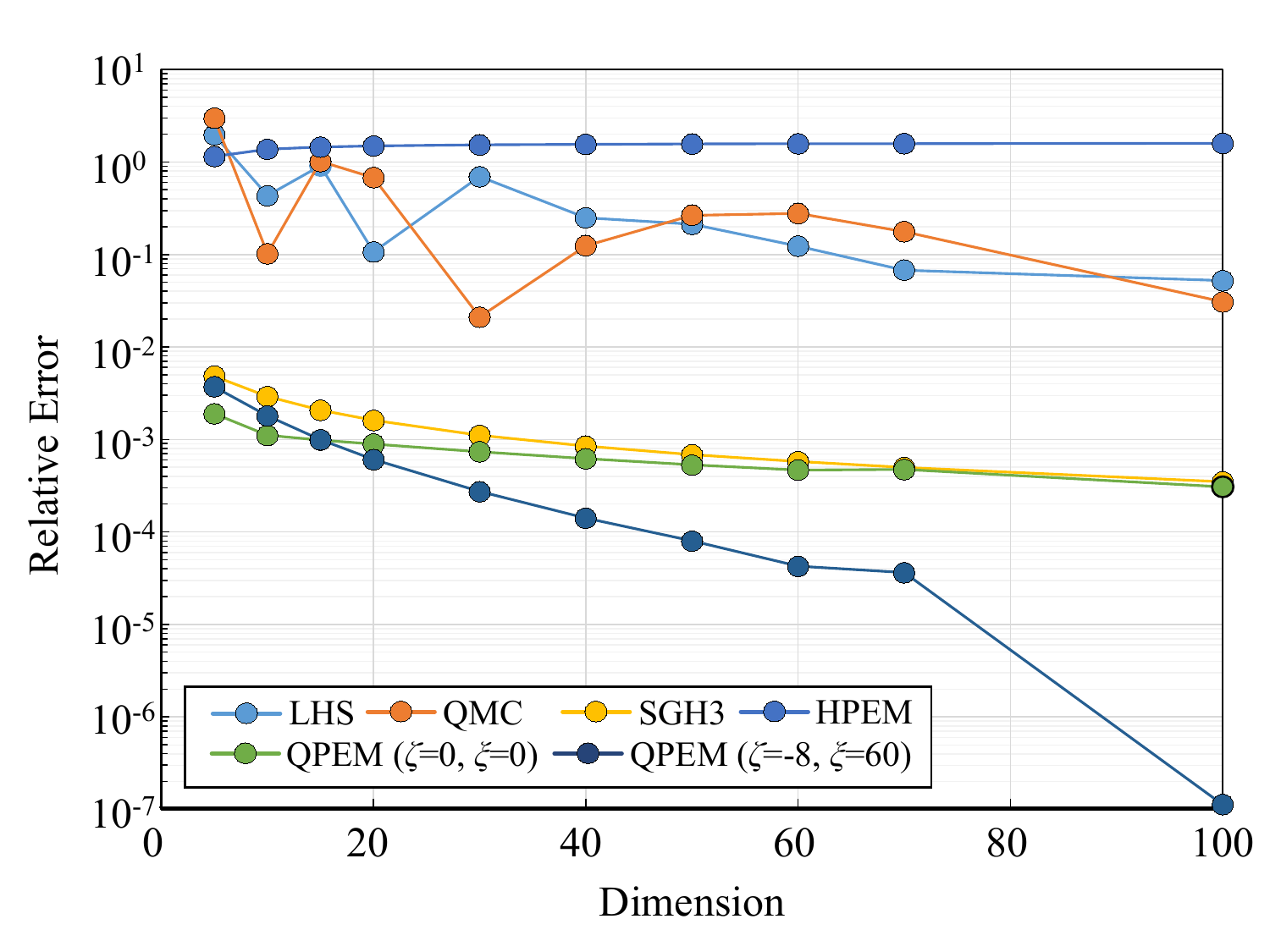}
		\caption{Skewness coefficient}
	\end{subfigure}
	\begin{subfigure}[tb]{0.45\textwidth}
		\includegraphics[width=\textwidth]{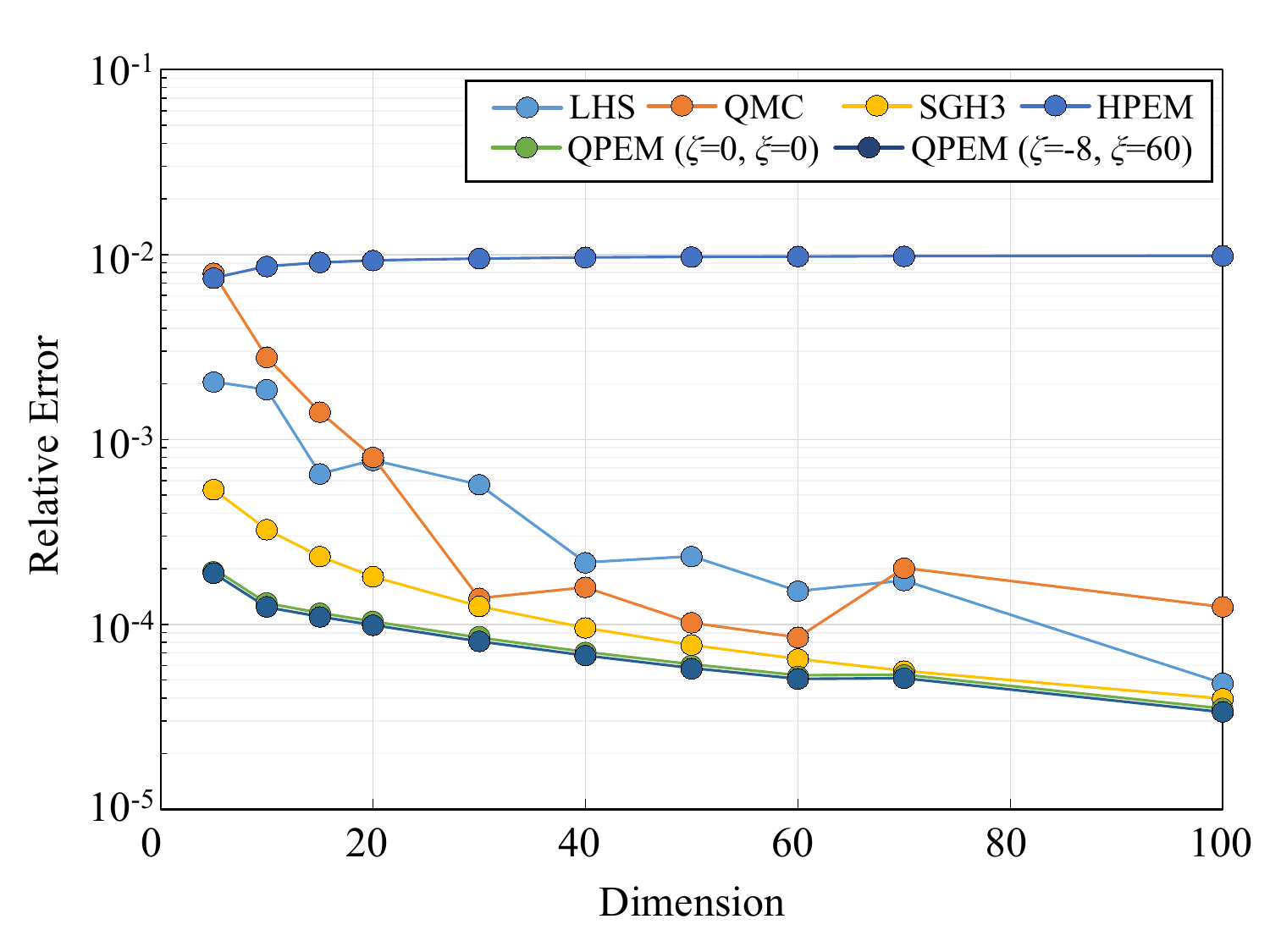}
		\caption{Kurtosis coefficient}
	\end{subfigure}
	\caption{Moments estimation for polynomial functions with varying dimensions.}
	\label{fig:ex1}
\end{figure}

\begin{itemize}
	\item[(a)] Since the response function comprises second-order polynomials, SGH3 and QPEM give exact results for the first two moments, and HPEM for the first moment. HPEM demonstrates great performance in estimating mean and standard deviation, achieving this with the fewest sigma points. However, its accuracy in estimating skewness and kurtosis significantly deteriorates, as expected based on its formulation and targeted accuracy. The errors of SGH3 for skewness and kurtosis are slightly larger than those of the QPEM. 
	\item[(b)] The scaling parameters $\zeta$ and $\xi$ of QPEM assist, in general, in reducing the relative error for the skewness and kurtosis estimations. For example, the scaling parameters improve the skewness results significantly in higher dimensions. As explained, however, these parameters cannot guarantee that the moment evaluations are always improved. For example, for the 5- and 10-$\mathcal{D}$ cases, skewness estimations become less accurate with the scaling parameter $\zeta$. In relation, even though all kurtosis estimations are improved by the parameter $\xi$, the effects are not substantial.	
	\item[(c)] Overall, the scaled QPEM presents the best performance in estimating the moments of the response function in all dimensions. As expected, LHS and QMC errors are more significant, in relation to QPEM, in the lower-dimensional examples here, because of the insufficient sample size of $2n^2+1$, used also for these methods. Thus, QPEM can be particularly preferable and competitive in problems with low to medium dimensions, e.g., up to 20-40 or so, where its computational efficiency is pronounced. However, as also shown in other examples, the performance of QPEM is extremely good and competitive in problems with higher dimensions as well.
\end{itemize}

\subsection{Structural example 1 – Roof truss}

\begin{figure}[!b]
	\centering
	\includegraphics[trim=0in 1.5in 0in 0.8in, clip=true, width=4in]{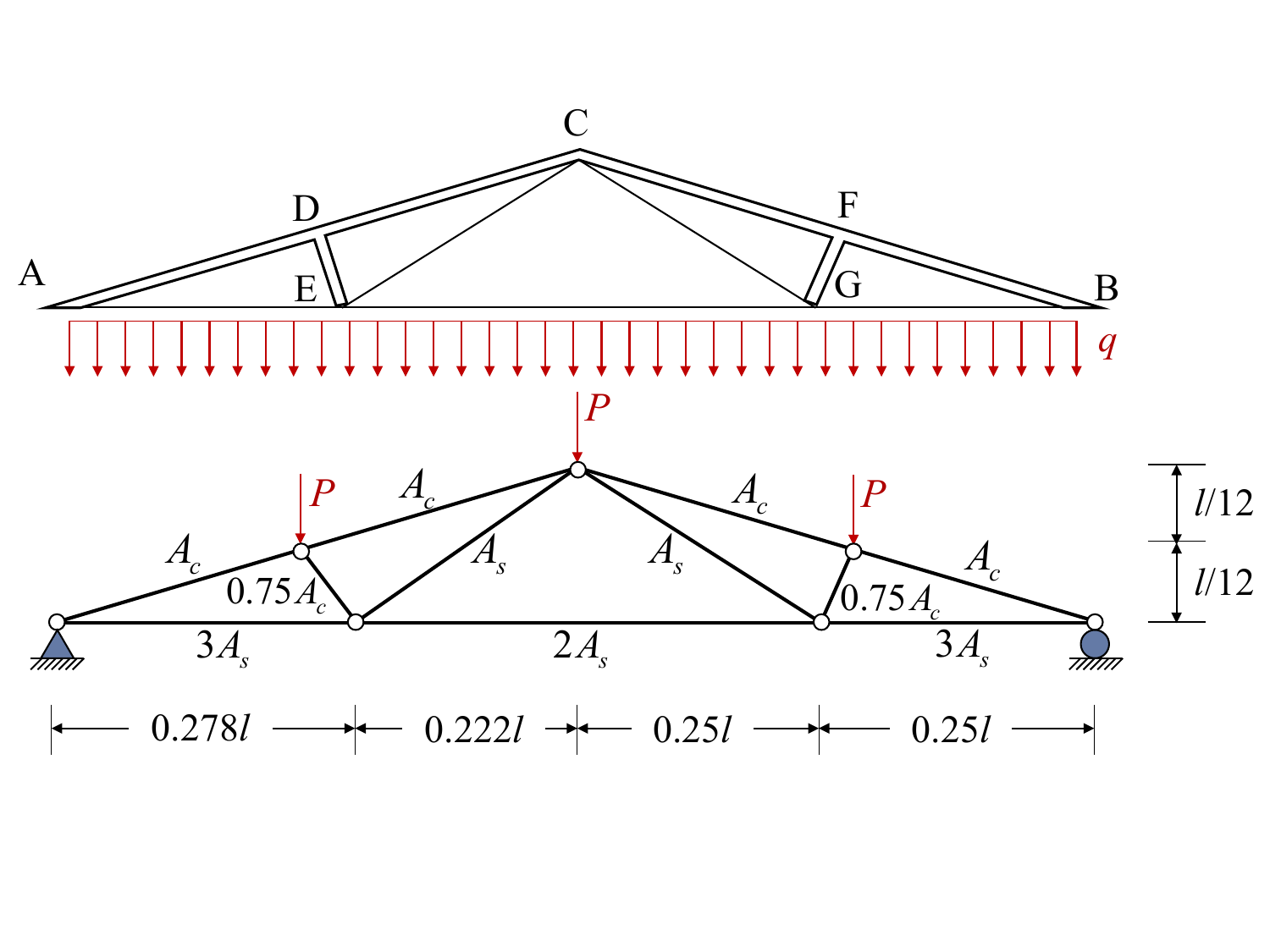}
	\caption{Schematic drawing of the roof structure.}
	\label{fig:ex2_rooftruss}
\end{figure}

\begin{table}[!b]
	\caption{Random variables for the roof structure example ($n=6$)}\label{table:ex2} 
	\begin{tabular*}{\textwidth}[t]{p{0.15\textwidth}p{0.35\textwidth}p{0.15\textwidth}p{0.15\textwidth}p{0.1\textwidth}}
		\toprule
		Variable (unit)     & Description            &Distribution & Mean         & COV$^{\text{a}}$ \\\midrule
		$q \,(N/m)$              & Uniformly distributed load      & Normal                & 200,000                & 0.07             \\
		$l \,(m)$                & Roof span                       & Normal                & 12                    & 0.01             \\
		$A_s \,(m^2)$            & Cross sectional area (steel)    & Normal                & $9.82 \times 10^{-4}$ & 0.06             \\
		$A_c \,(m^2)$            & Cross sectional area (concrete) & Normal                & $400 \times 10^{-4}$  & 0.12             \\
		$E_s \,(N/m^2)$          & Elastic modulus (steel)         & Normal                & $2 \times 10^{11}$    & 0.06             \\
		$E_c \,(N/m^2)$          & Elastic modulus (concrete)      & Normal                & $3 \times 10^{10}$    & 0.06             \\\bottomrule
		\multicolumn{5}{l}{$^{\text{a}}$ COV = coefficient of variation}  
	\end{tabular*}
\end{table}

\begin{table}[!t]
	\caption{Moment estimations for the roof truss structure ($n=6$)}\label{table:ex2_results} 
	\begin{tabular*}{\textwidth}[t]{p{0.30\textwidth} p{0.08\textwidth} p{0.12\textwidth} p{0.12\textwidth} p{0.12\textwidth} p{0.12\textwidth}}
		\toprule
		Method   & No. of \newline Points & Mean\newline[MPa]& STD$^{\text{b}}$\newline[MPa] & Skewness & Kurtosis\\\midrule
		MC \newline (95\% Upper CI\newline/95\% Lower CI$^{\text{c}}$)     & $10^6$   & 23.6689 \newline (23.6749 \newline/23.6648) &2.6027 \newline (2.6045\newline/2.5968)		& 0.3550 \newline (0.3582\newline/0.3470)&3.2633\newline (3.2890\newline/3.2543)\\  \hdashline[2pt/2pt]
		LHS                                     & 73      & 23.6752       & 2.6153       & 0.1682        & 2.7188 \\
		QMC                                    & 73       & 23.7111       & 2.6119        & 0.6210        & 5.5348 \\
		SGH3                                    & 85        & 23.6687       & 2.5977        & 0.3102        & 2.7435 \\
		HPEM               & 13        & 23.6703       & 2.5847        & 0.0082        & 1.7741 \\
		Unscaled QPEM ($r=\sqrt{3}$)& 73        & 23.6688       & 2.5995        & 0.3286        & 2.9724 \\
		Scaled QPEM$^{\text{a}}$ ($r=\sqrt{3}$)& 73        & 23.6688       & 2.5995        & 0.3350       & 2.9768 \\
		Unscaled QPEM ($r=3$)& 73        & 23.6688       & 2.5995        & 0.3368        & 3.0869 \\
		Scaled QPEM$^{\text{a}}$ ($r=3$)& 73        & 23.6688       & 2.5995        & 0.3432        & 3.0913 \\\bottomrule
		\multicolumn{6}{l}{$^{\text{a}}$ $\zeta=-8$ and $\xi=60$}\\  
            \multicolumn{6}{l}{$^{\text{b}}$ STD = standard deviation}\\
            \multicolumn{6}{l}{$^{\text{c}}$ CI = Bootstrap Confidence Interval}
	\end{tabular*}
\end{table}

\begin{figure}[!t]
	\centering
	\includegraphics[width=5.0in]{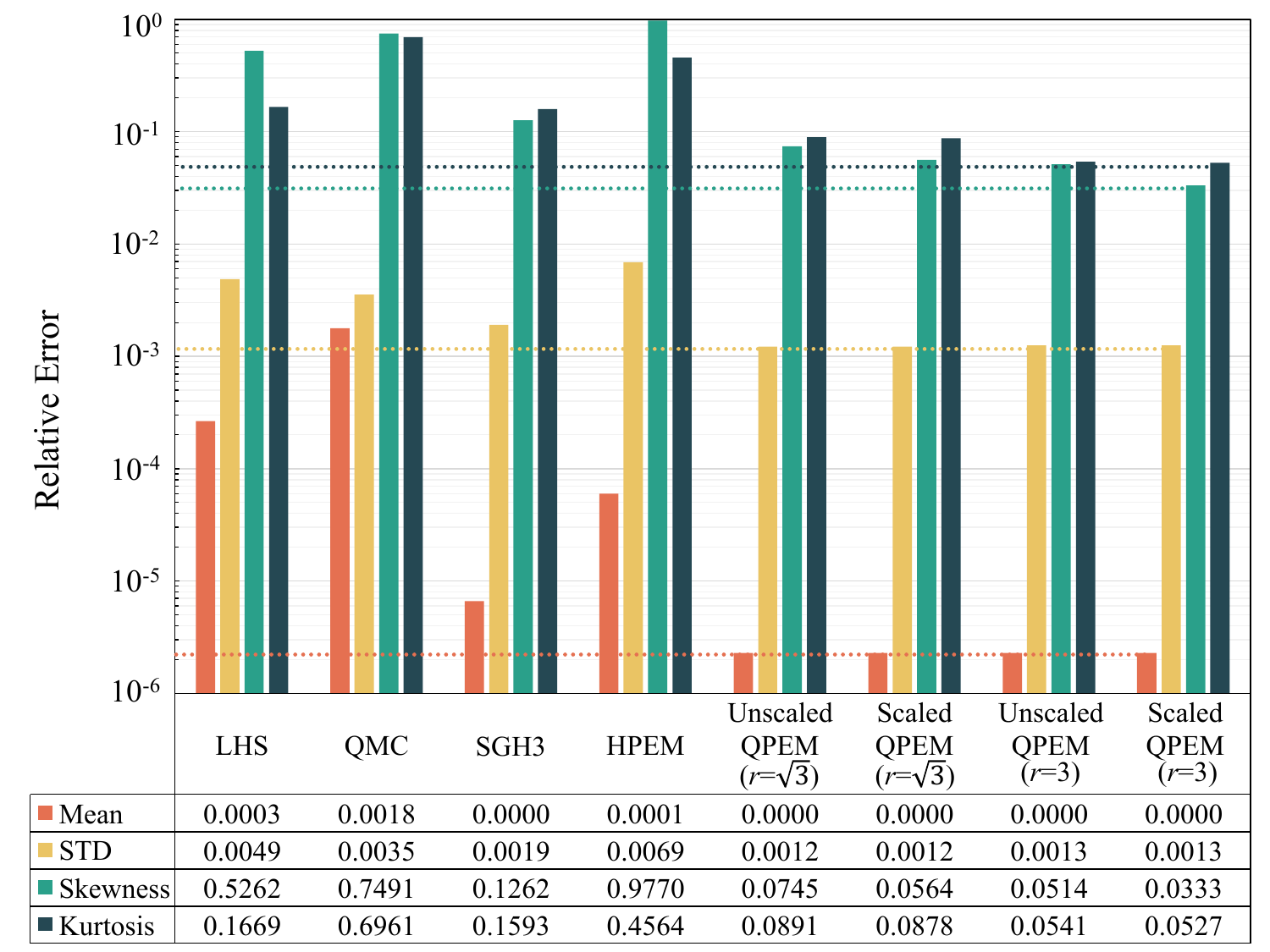}
	\caption{Relative errors for the roof truss structure example.}
	\label{fig:ex2_results}
\end{figure}

A roof truss structure adapted from ~\cite{song2009subset, zhou2013new} is investigated in this example, shown in Fig.~\ref{fig:ex2_rooftruss}, where the compression and tension members consist of concrete and steel, respectively. The uniformly distributed load $q$ is applied as nodal load $P$=$ql/4$ at nodes C, D, and F. The vertical deflection $y$ of the peak (node C) is obtained as:
\begin{equation}
	\label{eq:ex2_rooftruss}
	y=\dfrac{ql^2}{2}\left(\dfrac{3.81}{A_c E_c}+\dfrac{1.13}{A_s E_s} \right)
\end{equation}
with descriptions and distributions of the variables: summarized in Table~\ref{table:ex2} for this 6-$\mathcal{D}$ problem ($n$=6). Pairwise correlation coefficients are given as $\rho_{lA_{s}}$=$\rho_{lA_{c}}$=0.3 and $\rho_{A_{c}A_{s}}$=$\rho_{E_{c}E_{s}}$=0.5, while the remaining variables are independent ($\rho$=0). As before, the first four moments have been estimated, of the vertical deflection at the peak this time. The QPEM has been examined with two values for the parameter $r$=$\sqrt{3}$ and 3, for both scaled and unscaled cases, to showcase their influence on results. Results are reported in Table~\ref{table:ex2_results} and Fig.~\ref{fig:ex2_results}. The relative errors in Fig.~\ref{fig:ex2_results} are now with respect to MC reference solutions with $10^6$ samples, and the 95\% confidence intervals in Table~\ref{table:ex2_results} are computed based on $10^4$ bootstrap samples with the Bias-Corrected Accelerated (BCA) bootstrap method \cite{diciccio1996bootstrap}. BCA involves resampling with replacement, from the samples obtained in the initial MC simulation and recomputing the moments for each resampled dataset. BCA does not strictly rely on the assumption of normality and is versatile in providing the distribution of various estimators, including ones involving higher-order moments \cite{diciccio1996bootstrap}. 

Commenting on results, generally, all methods approximate the first two moments well. For skewness and kurtosis, however, the LHS and QMC have substantial relative errors, ranging from 16.69\% to 74.91\%, since their 73 sample points, as many as the QPEM one, are not sufficient to provide accurate enough results. The unscaled QPEM with $r$=$\sqrt{3}$ has much lower relative errors for skewness (7.45\%) and kurtosis (8.91\%), and applying the scaling parameters $\zeta$ and $\xi$ improves the errors (5.64\% for skewness and 8.78\% for kurtosis) in this case. The unscaled QPEM with $r$=3 has even lower errors for skewness and kurtosis (5.14\% and 5.41\%, respectively), and the scaling parameters $\zeta$ and $\xi$ also work to reduce the errors greatly in this case with $r$=3, to 3.33\% and 5.27\%, respectively.

Overall, the QPEM is the most efficient method in this low-dimensional example, the scaled version offers improvements in this case, and the default $r$=3 value performs better than the $r$=$\sqrt{3}$ one. 

\subsection{Structural example 2 – Six-story structure under horizontal static loads}

\begin{figure}[!b]
	\centering
	\includegraphics[trim=0.7in 1in 0.5in 0.7in, clip=true, width=4.5in]{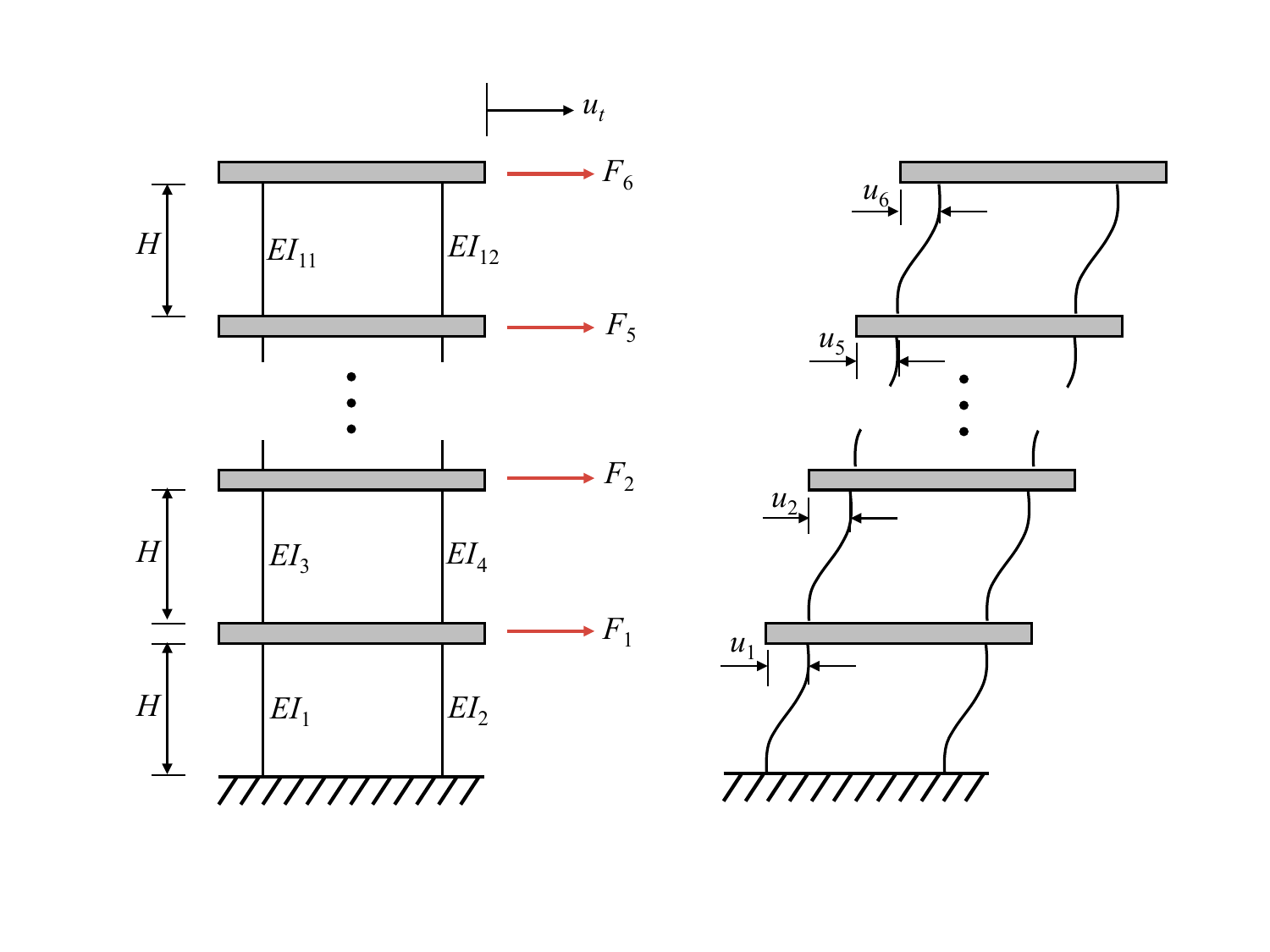}
	\caption{Six-story structure subjected to horizontal static loads.}
	\label{fig:ex3_six_story_structure}
\end{figure}

A structural example of a six-story structure is studied in this case, shown in Fig.~\ref{fig:ex3_six_story_structure}, as adopted and modified from \cite{bucher2009computational}. The structure is subjected to six horizontal static loads $F_i$ ($i$=$1,\cdots,6$), the floor slabs are assumed to be rigid, and the columns with a height $H$=$4m$ have bending stiffnesses $EI_k$ ($k$=$1,\cdots,12$). The horizontal loads and the bending stiffnesses are considered as random variables resulting in $n$=18 dimensions for this problem. The loads are normally distributed with $F_i\sim N(20,\,6)~[kN]$, and the stiffnesses are also normally distributed with $EI_k \sim N(10^4,\,10^3)~[N \cdot m^2]$. The random variables are assumed to be pairwise linearly correlated, as $\rho_{F_iF_j}$=$0.5$ and $\rho_{EI_iEI_j}$=$0.1$, where $i \neq j$, and there is no correlation between loads $F$ and stiffnesses $EI$, i.e., $\rho_{F_iEI_j}$=$0$. Structural analysis is performed based on linear elastic behavior assumption, neglecting the effects of gravity. The total horizontal displacement at the top, $u_t$, is considered the response function in this example, as follows:
\begin{equation}
	\label{eq:total_horizontal_deflection_six_story_structure}
	u_t= \sum_{i=1}^{6} u_i
\end{equation}
where:
\begin{equation}
	\label{eq:horizontal_deflection}
	\begin{split}
		u_6 &= \dfrac{F_6 H^3}{12(EI_{11}+EI_{12})}\\
		u_5 &= \dfrac{(F_5+F_6)H^3}{12(EI_{9}+EI_{10})}\\
		u_4 &= \dfrac{(F_4+F_5+F_6)H^3}{12(EI_{7}+EI_{8})}
	\end{split}
	\qquad \qquad
	\begin{split}
		u_3 &= \dfrac{(F_3+F_4+F_5+F_6) H^3}{12(EI_{5}+EI_{6})}\\
		u_2 &= \dfrac{(F_2+F_3+F_4+F_5+F_6) H^3}{12(EI_{3}+EI_{4})}\\
		u_1 &= \dfrac{(F_1+F_2+F_3+F_4+F_5+F_6)H^3}{12(EI_{1}+EI_{2})}
	\end{split}
\end{equation}

\begin{figure}[!t]
	\centering
	\includegraphics[width=5.0in]{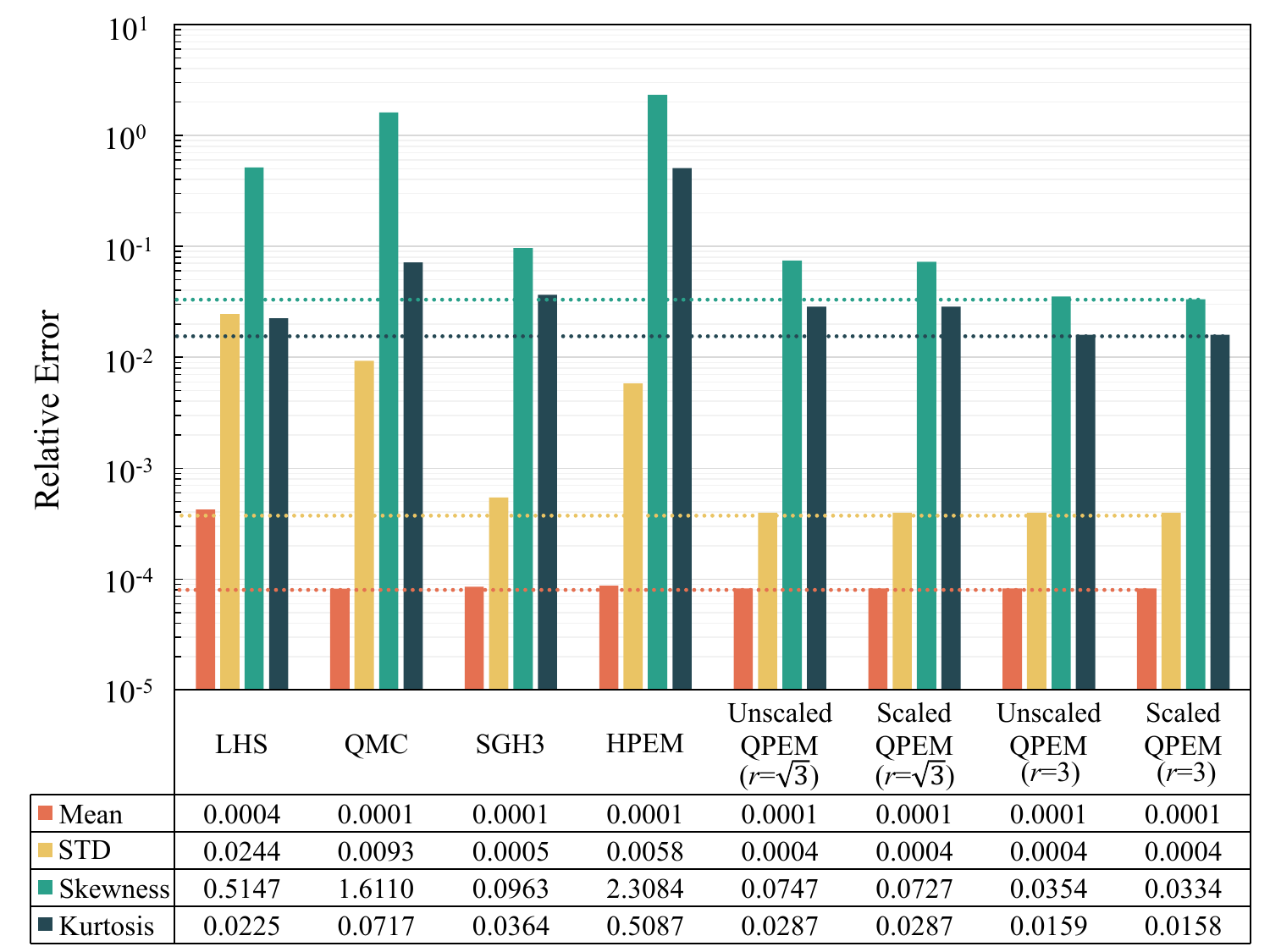}
	\caption{Relative errors for the six-story structure example.}
	\label{fig:relative_error_six_story_structure}
\end{figure}

The first four moments of the top horizontal displacement are estimated by all methods used for comparison, and QPEM is again implemented with $r$=$\sqrt{3}$ and 3, for both scaled and unscaled versions. Sample sizes and results can be seen in Table~\ref{table:moments_six_story_structure} and Fig.~\ref{fig:relative_error_six_story_structure}, and, same as in the previous example, the relative errors in Fig.~\ref{fig:relative_error_six_story_structure} are again with respect to MC reference solutions with $10^6$ samples, and the 95\% confidence intervals in Table~\ref{table:moments_six_story_structure} are computed based on $10^4$ bootstrap samples with the BCA method. The QPEM shows again the best performance in estimating all four moments for this 18-$\mathcal{D}$ problem, with the LHS and QMC sampling methods not exhibiting that accurate results, especially for the skewness. Among the QPEM variants, the scaled one with $r$=3 provides again the best results.

\begin{table}[!b]
	\caption{Moment estimations for the total horizontal displacement ($n=18$)}\label{table:moments_six_story_structure} 
	\begin{tabular*}{\textwidth}[t]{p{0.30\textwidth} p{0.08\textwidth} p{0.12\textwidth} p{0.12\textwidth} p{0.12\textwidth} p{0.12\textwidth}}
		\toprule
		Method   										&No. of \newline Points	 & Mean\newline[mm]										& STD$^{\text{b}}$ \newline[mm] & Skewness & Kurtosis\\\midrule
		MC \newline (95\% Upper CI\newline/95\% Lower CI$^{\text{c}}$)   & $10^6$    	& 112.6359 \newline (112.6881\newline/112.5843) & 26.7536	\newline (26.7913\newline/26.7178)		&  0.0509\newline (0.0544\newline/0.0447) & 3.0285\newline (3.0294\newline/3.0092)\\\hdashline[2pt/2pt]
		LHS                                     & 649       & 112.5881       & 26.0997       & 0.0771       & 2.9603 \\
		QMC                                    & 649       & 112.6451       & 27.0027        & 0.1329        & 2.8113 \\
		SGH3                                    & 685      & 112.6263      & 26.7390        & 0.0460        & 2.9182 \\
		HPEM               & 37        & 112.6457       & 26.5981        & -0.0666       & 4.5691 \\
		Unscaled QPEM ($r=\sqrt{3}$) & 649       & 112.6266       & 26.7430       & 0.0471        & 2.9417 \\
		Scaled QPEM$^{\text{a}}$ ($r=\sqrt{3}$)& 649       & 112.6266       & 26.7430        & 0.0472        & 2.9417 \\
		Unscaled QPEM ($r=3$) & 649       & 112.6266       & 26.7430       & 0.0491        & 2.9804 \\
		Scaled QPEM$^{\text{a}}$ ($r=3$)& 649       & 112.6266       & 26.7430        & 0.0492        & 2.9805 \\\bottomrule
		\multicolumn{6}{l}{$^{\text{a}}$ $\zeta=-8$ and $\xi=60$}\\
            \multicolumn{6}{l}{$^{\text{b}}$ STD = standard deviation}\\
            \multicolumn{6}{l}{$^{\text{c}}$ CI = Bootstrap Confidence Interval}  
	\end{tabular*}
\end{table}

\begin{figure}[!b]
	\centering
	\includegraphics[trim=1in 2in 1in 2in, clip=true, width=3.15in]{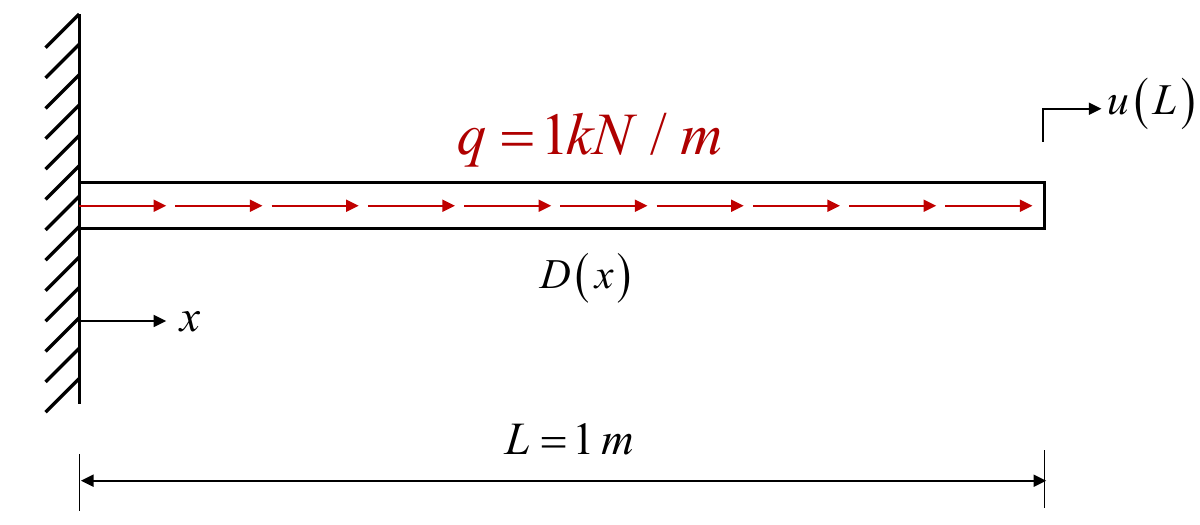}
	\caption{Elastic bar with random field axial rigidity.}
	\label{fig:ex4_1d_bar}
\end{figure}

\subsection{Linear finite element analysis 1 – Elastic bar with random axial rigidity}

In this example, a finite element analysis (FEA) of a linear elastic bar with random axial rigidity is studied, as adapted from \cite{papaioannou2019pls}, and shown in Fig.~\ref{fig:ex4_1d_bar}. The displacement of the bar $u\left(x\right)$ is given by:

\begin{equation}
	\label{eq:ex4_1d_bar}
	-\dfrac{d}{dx}\left(D(x)\dfrac{d}{dx}u\right)=q(x) \,\,\,\,x\in[0,L]
\end{equation}
The axial rigidity $D(x)$ is modeled as a homogeneous normal random field with mean $\mu_D$=$100kN$ and standard deviation $\sigma_D$=$10kN$. The autocorrelation function of the random field $D$ is given as $\rho_D (\Delta x)$= $\exp\left(-\Delta x^2 / l_{x}^2\right)$ with the correlation parameter $l_x$ equal to $0.20m$. The random field $D(x)$ is represented by the KL expansion with 20 terms, i.e., 20 independent normal variables are involved in this problem. A deterministic load $q$=$1\,kN/m$ is applied on the elastic bar, that has a length of $L$=$1m$. Eq.~\eqref{eq:ex4_1d_bar} is particularly solved for the displacement at the tip, $u(L)$, by FEA with 100 linear bar elements. 

\begin{table}[!b]
	\caption{Moment estimations for the displacement at the tip of the elastic bar ($n=20$)}\label{table:ex4_results} 
	\begin{tabular*}{\textwidth}[t]{p{0.30\textwidth} p{0.08\textwidth} p{0.12\textwidth} p{0.12\textwidth} p{0.12\textwidth} p{0.12\textwidth}}
		\toprule
		Method   &No. of\newline Points & Mean\newline[mm]& STD$^{\text{b}}$\newline[mm] & Skewness & Kurtosis\\\midrule
		MC\newline (95\% Upper CI\newline/95\% Lower CI$^{\text{c}}$)   & $10^6$    &5.0512 \newline (5.0523\newline/5.0510)  &
		0.3263 \newline (0.3264\newline/0.3254) &0.4651\newline (0.4744\newline/0.4616)  & 3.4844\newline (3.5291\newline/3.4755)  \\ \hdashline
		LHS                                     & 801       & 5.0517       & 0.3306        & 0.3563        & 3.0398\\
		QMC                                    & 801       & 5.0479       & 0.3250        & 0.3768        & 2.9702 \\
		SGH3                                    & 871       & 5.0515       & 0.3253        & 0.3510        & 2.6187 \\
		HPEM            						& 71         & 5.0504       & 0.3160        & -0.0305        & 1.8970 \\
		Unscaled QPEM ($r=\sqrt{3}$) 					& 801       & 5.0515       & 0.3258        & 0.3950        & 3.0292 \\
		Scaled QPEM$^{\text{a}}$ ($r=\sqrt{3}$)& 801       & 5.0515       & 0.3258        & 0.4266        & 3.0666 \\
		Unscaled QPEM ($r=3$) 					& 801       & 5.0515       & 0.3265        & 0.4484       & 3.4350 \\
		Scaled QPEM$^{\text{a}}$ ($r=3$)& 801       & 5.0515       & 0.3265       & 0.4799       & 3.5007 \\\bottomrule
		\multicolumn{6}{l}{$^{\text{a}}$ $\zeta=-8$ and $\xi=60$}\\
            \multicolumn{6}{l}{$^{\text{b}}$ STD = standard deviation}\\
            \multicolumn{6}{l}{$^{\text{c}}$ CI = Bootstrap Confidence Interval}     
	\end{tabular*}
\end{table}

\begin{figure}[!b]
	\centering
	\includegraphics[width=5.0in]{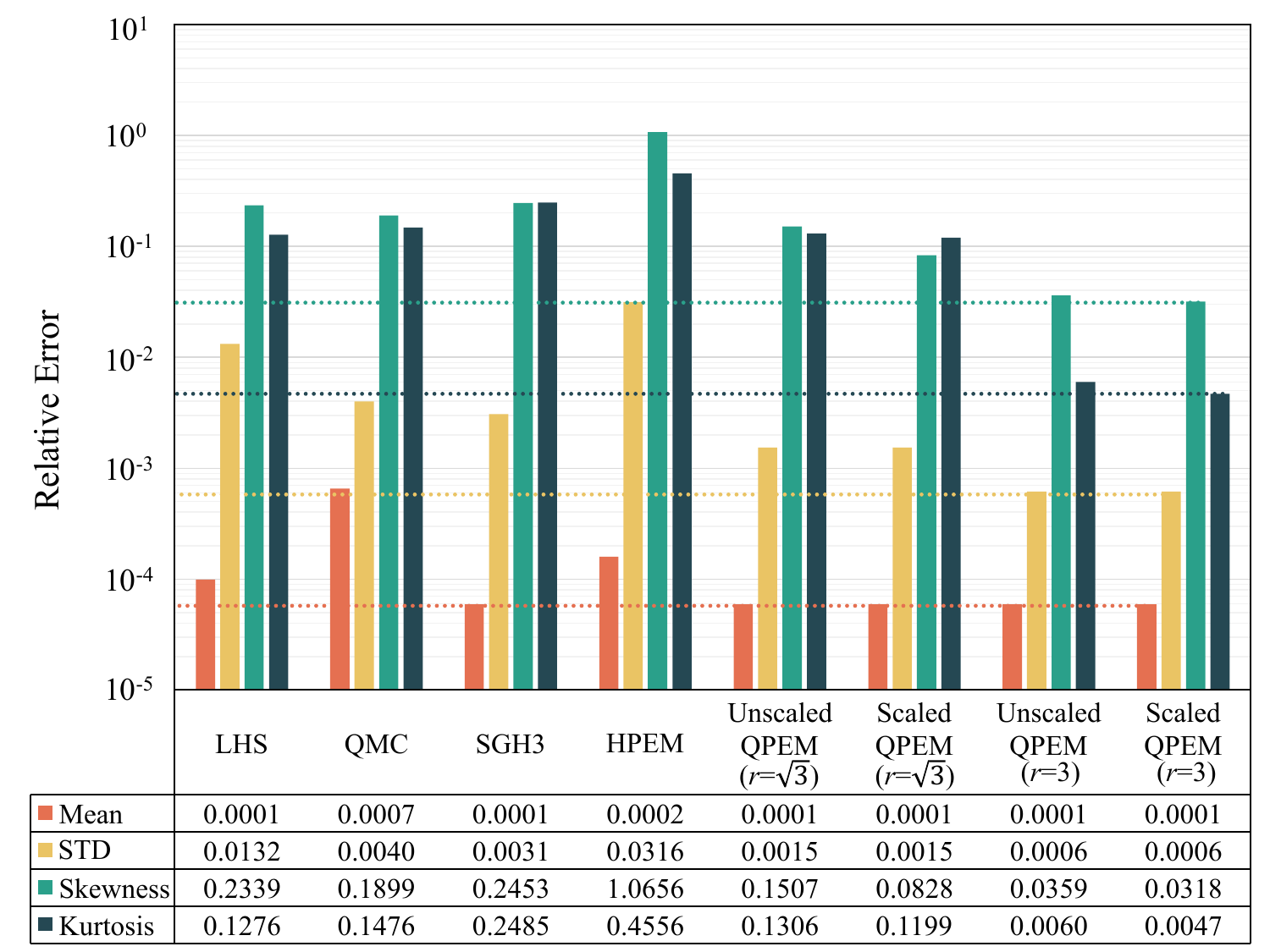}
	\caption{Relative errors for the elastic bar example.}
	\label{fig:ex4_results}
\end{figure}

Table~\ref{table:ex4_results} and Figure~\ref{fig:ex4_results} show the computed results by each method. MC samples and confidence intervals are exactly used and computed as in previous examples, and all methods other than the QPEM exhibit significant errors in the skewness and kurtosis estimations. The unscaled QPEM with $r$=$\sqrt{3}$ estimates a somewhat more accurate skewness and kurtosis, in comparison to other methods, while the scaled QPEM version significantly improves the third and fourth-order moments estimations. The parameter $r$=3 again leads to important improvements, compared to $r$=$\sqrt{3}$, being able to better minimize high-order terms errors, and its scaled version provides the best results.

Overall, this example demonstrates the excellent performance and applicability of the QPEM in problems with random fields, and the superiority of the method in relation to variance reduction sampling and sparse grid quadrature techniques with same/analogous computational cost.

\subsection{Linear finite element analysis 2 – Plate with a hole with random elastic modulus} 

\begin{figure}[!t]
	\centering
	\includegraphics[trim=1.5in 0.8in 1.5in 0.8in, clip=true, width=3.5in]{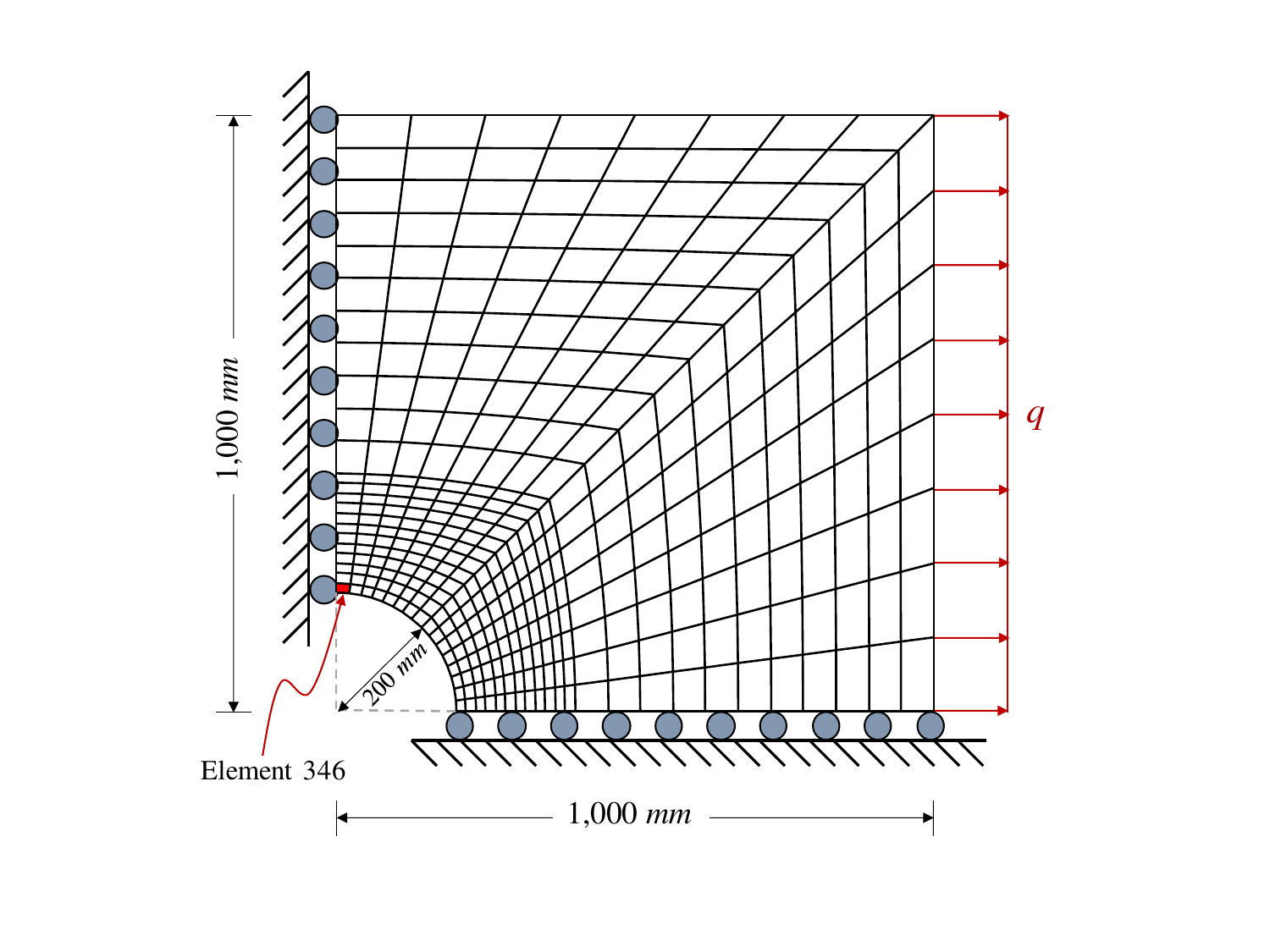}
	\caption{Quarter-plate with a circular hole and random field elastic modulus.}
	\label{fig:ex_plate_hole}
\end{figure}

A low-carbon steel plate is considered in this example, shown in Fig.~\ref{fig:ex_plate_hole} \cite{papaioannou2019pls}. The geometry of the plate is $2,000\,mm \times 2,000\,mm$, and a circular hole with a radius of $200\,mm$ is located at the center. Due to symmetry in geometry only a quarter of the plate is analyzed. The Poisson ratio is assumed to be $\nu$=0.3 and the plate thickness is $t$=$10\,mm$. A deterministic uniform load of $q$=$100\,N/mm^2$ is applied on the right edge of the plate. The elastic modulus  $E(x,y)$ of the plate is modeled as a two-dimensional homogeneous normal random field with mean $\mu_E$=$2\times10^5\,N/mm^2$ and standard deviation $\sigma_E$=$2 \times 10^4\,N/mm^2$. The autocorrelation function of the random field $E$ is given as:
\begin{equation}
	\label{eq:autocorrelation}
	{{\rho }_{E}}\left( \Delta x,\Delta y \right)=\exp \left[ -{{\left( \frac{\Delta x}{{{l}_{x}}} \right)}^{2}}-{{\left( \frac{\Delta y}{{{l}_{y}}} \right)}^{2}} \right]
\end{equation}
where $l_x$ and $l_y$ are the correlation parameters equal to $300\,mm$ in this example. Further increasing the analysis dimensions in this case, the KL expansion with 30 terms is used to describe the spatial stochastic fields, i.e., $n=30$. The stress, strain, and displacement fields of the plate should satisfy the Navier-Cauchy equilibrium equation \cite{johnson2012numerical}, simplified under the plane stress assumptions as:
\begin{equation}
	\label{eq:navier_cauchy}
	G(x,y)\bm{\nabla}^2 \mathbf{u}+\dfrac{E(x,y)}{2(1-\nu)}\bm{\nabla} (\bm{\nabla} \cdot \mathbf{u})+\mathbf{f}=\mathbf{0}
\end{equation}
where $G(x, y)$ is the shear modulus, i.e., $E(x,y)/[2(1+\nu)]$, $\mathbf{u}$ is the displacement vector, and $\mathbf{f}$ is the body force vector acting on the plate. Eq.~\eqref{eq:navier_cauchy} is solved through FEA with a mesh of 368 four-node quadrilateral elements, shown in Fig.~\ref{fig:ex_plate_hole}. The maximum first principal plane stress in the plate is considered as the response output in this example, most often occurring at element 346 (see red marker in Fig.~\ref{fig:ex_plate_hole}).

\begin{table}[!t]
	\caption{Moment estimations for the maximum first principal stress of the plate with a hole ($n=30$)}
	\label{table:ex_results} 
	\begin{tabular*}{\textwidth}[t]{p{0.30\textwidth} p{0.08\textwidth} p{0.12\textwidth} p{0.12\textwidth} p{0.12\textwidth} p{0.12\textwidth}}
		\toprule
		Method   &No. of \newline points& Mean\newline[MPa]& STD$^{\text{b}}$\newline[MPa] & Skewness & Kurtosis\\\midrule
		MC\newline (95\% Upper CI\newline/95\% Lower CI$^{\text{c}}$)    & $10^6$    & 316.7064 \newline (316.7417\newline/316.6719)  & 17.6459 \newline (17.6708\newline/17.6205)  &0.0544 \newline (0.0597\newline/0.0492)  & 3.1224 \newline (3.1142\newline/3.0914) \\ \hdashline[2pt/2pt]
		LHS                                     & 1,801       & 316.6956    & 17.0706        & 0.0420       & 3.2442 \\
		QMC                                    & 1,801       & 316.4213       & 17.7969        & 0.0463        & 3.2485 \\
		SGH3                                    & 1,861       & 316.7053       & 17.6299        & 0.0492        & 2.8977 \\
		HPEM              						& 61         & 316.6972       & 17.4672        & 0.0158        & 0.7383 \\
		Unscaled QPEM ($r=\sqrt{3}$)  					& 1,801       & 316.7054      & 17.6326       & 0.0507        & 2.9866 \\
		Scaled QPEM$^{\text{a}}$ ($r=\sqrt{3}$)	& 1,801       & 316.7054       & 17.6326        & 0.0508        & 2.9867 \\
		Unscaled QPEM ($r=3$)  					& 1,801       & 316.7058      & 17.6352       & 0.0532        & 3.1467 \\
		Scaled QPEM$^{\text{a}}$ ($r=3$)	& 1,801       & 316.7058       & 17.6352        & 0.0533        & 3.1467 \\\bottomrule
		\multicolumn{6}{l}{$^{\text{a}}$ $\zeta=-8$ and $\xi=60$}\\
            \multicolumn{6}{l}{$^{\text{b}}$ STD = standard deviation}\\
            \multicolumn{6}{l}{$^{\text{c}}$ CI = Bootstrap Confidence Interval} 
	\end{tabular*}
\end{table}

\begin{figure}[!t]
	\centering
	\includegraphics[width=5in]{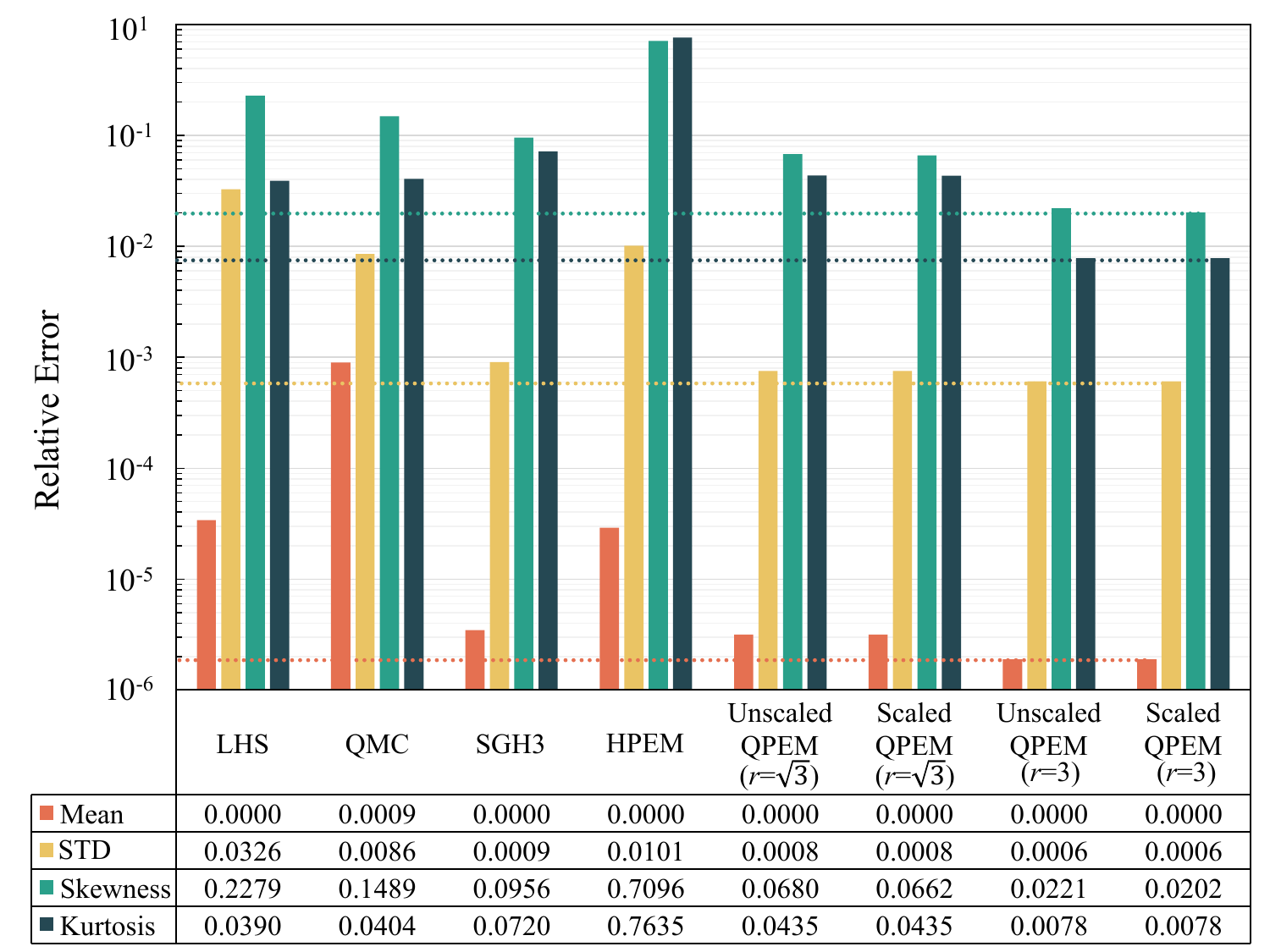}
	\caption{Relative errors for the plate with a hole example.}
	\label{fig:ex_results}
\end{figure}

Further analysis details also follow here what has already been done and reported in previous examples.  As shown in Table~\ref{table:ex_results} and Fig.~\ref{fig:ex_results}, the QPEM with $r$=3 demonstrates the best performance again over all other methods. In general, all methods estimate the mean, standard deviation, and kurtosis accurately, except for the HPEM, which, due to its design and targeted accuracy, cannot adequately capture output moments higher than the variance. For skewness, both LHS and QMC exhibit significant relative errors of 22.79\% and 14.89\%, respectively, while the unscaled QPEM with $r$=3 provides a quite accurate skewness estimation, with an error of only 2.21\%, and the relevant scaled version only minimally improves the skewness estimator in this case, with a reported 2.02\% relative error.  

\begin{figure}[!t]
	\begin{subfigure}{0.49\textwidth}
		\includegraphics[trim=0.4in 0.4in 0.4in 0.4in, clip=true, width=3in]{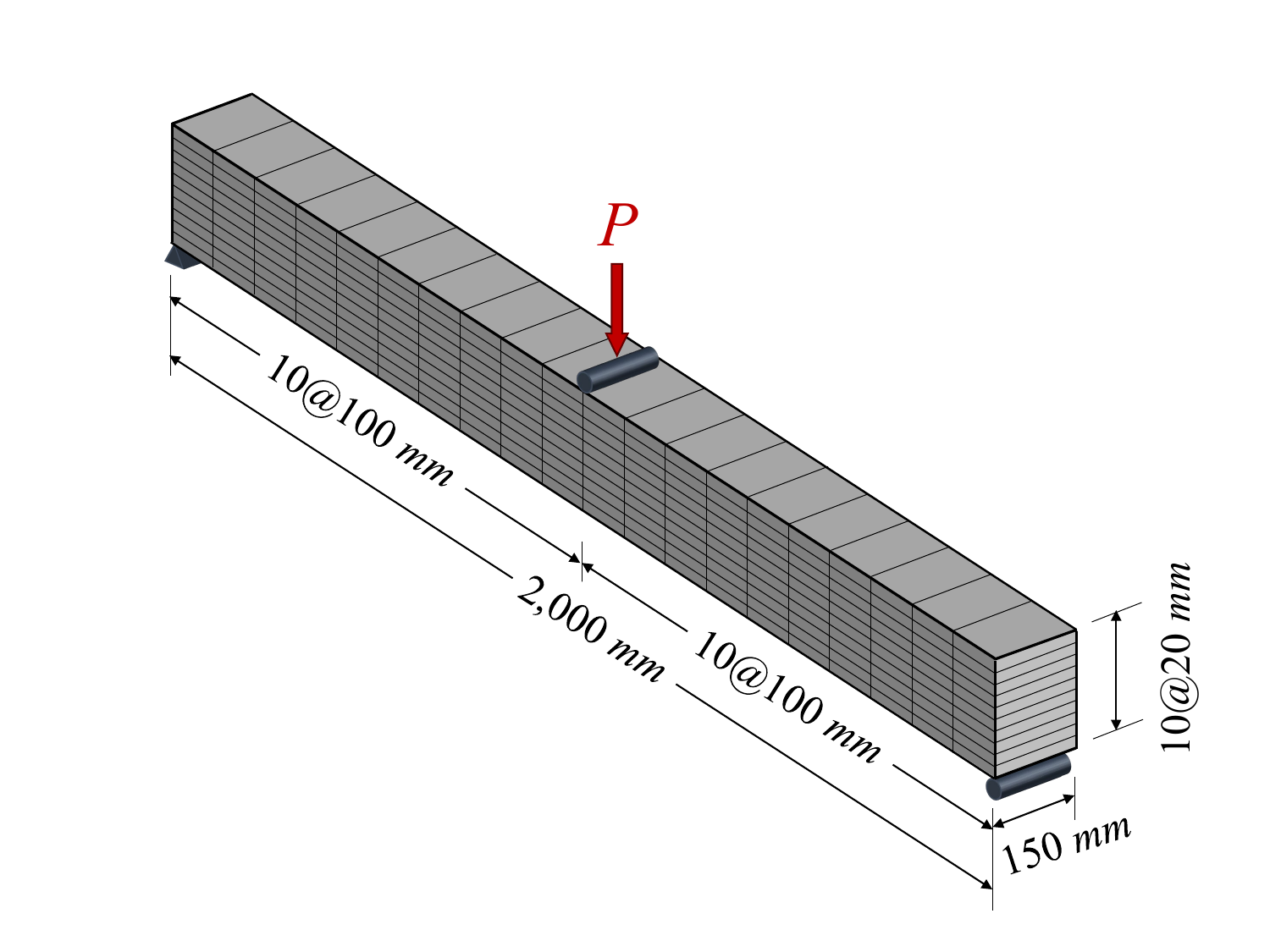}	
		\caption{Geometry and mesh}
	\end{subfigure}
	\begin{subfigure}{0.49\textwidth}
		\includegraphics[trim=0in 0.3in 0in 0.4in, clip=true, width=3in]{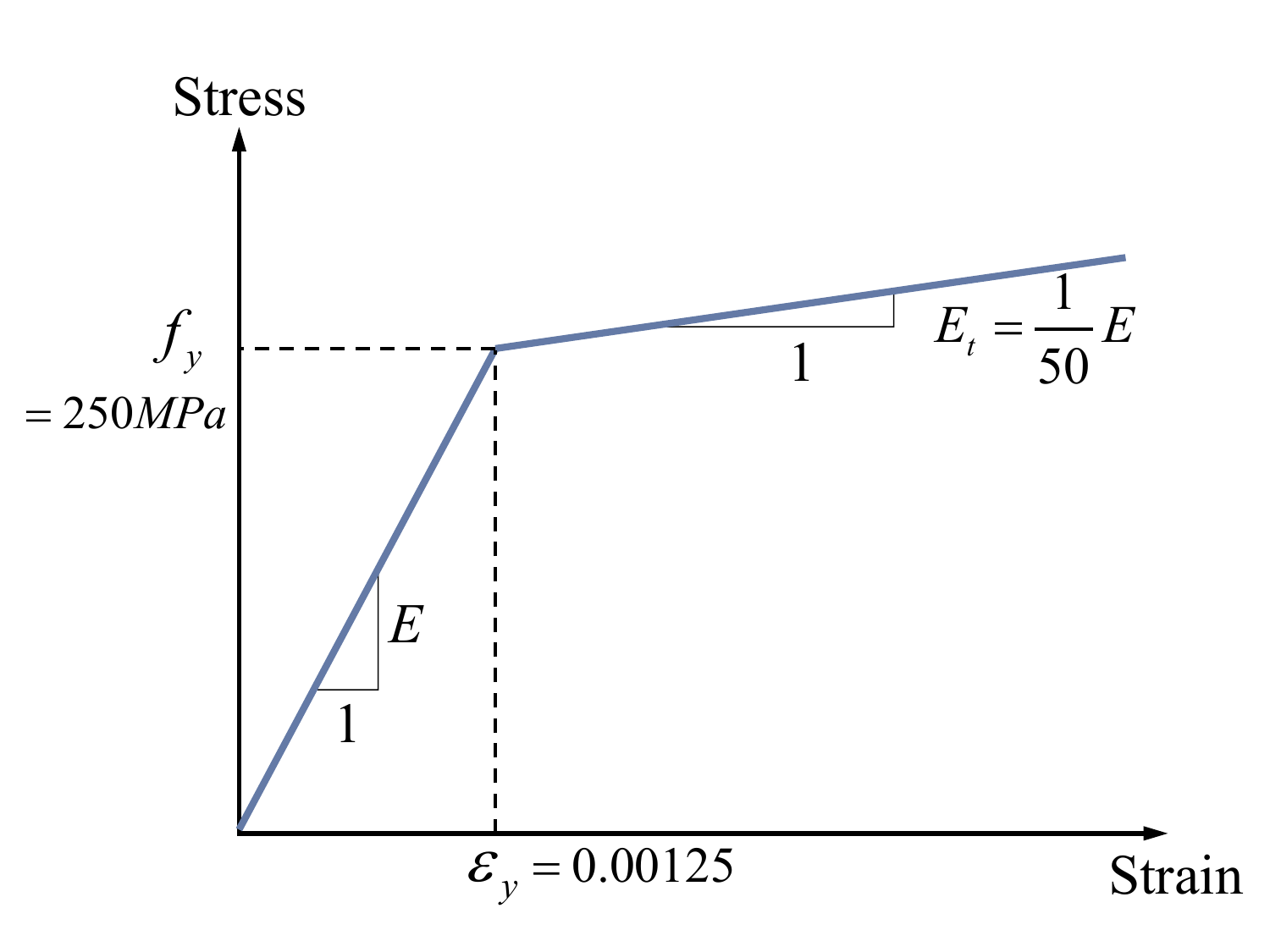}
		\caption{Elastoplastic bi-linear model}
	\end{subfigure}
	\caption{Simply supported beam with material nonlinearity and random field elastic modulus.}
	\label{fig:ex6_beam}
\end{figure}

\begin{table}[!b]
	\caption{Moment estimations for maximum beam displacement ($n=40$)}\label{table:ex6_results} 
	\begin{tabular*}{\textwidth}[t]{p{0.30\textwidth} p{0.08\textwidth} p{0.12\textwidth} p{0.12\textwidth} p{0.12\textwidth} p{0.12\textwidth}}
		\toprule
		Method   &No. of \newline points& Mean\newline[$mm$]& STD$^{\text{b}}$\newline[$mm$] & Skewness & Kurtosis\\\midrule
		MC\newline (95\% Upper CI\newline/95\% Lower CI$^{\text{c}}$)   & $10^6$    & 28.1685 \newline (28.1728\newline/28.1662)& 1.6707\newline (1.6721\newline/1.6671)&0.3767\newline (0.3833\newline/0.3724)& 3.2737\newline (3.2927\newline/3.2591)\\ \hdashline[2pt/2pt]
		LHS                                     & 3,201       & 28.1657      & 1.6647        & 0.3218        & 3.0620 \\
		QMC                                    & 3,201       & 28.1782       & 1.6616        & 0.3846        & 3.0615 \\
		SGH3                                    & 3,281       & 28.2014       & 1.6675        & 0.2620       & 2.8273\\
		HPEM              					& 81         & 28.1587       & 1.6481        & 0.0660        & 1.8801 \\
		Unscaled QPEM ($r=\sqrt{3}$)  & 3,201       & 28.1820       & 1.6683        & 0.3240        & 3.0516 \\
		Scaled QPEM$^{\text{a}}$ ($r=\sqrt{3}$) & 3,201        & 28.1820       & 1.6683        & 0.3326        & 3.0581 \\
		Unscaled QPEM ($r=3$)  & 3,201       & 28.1680       & 1.6697        & 0.3685        & 3.1575 \\
		Scaled QPEM$^{\text{a}}$ ($r=3$) & 3,201        & 28.1680       & 1.6697        & 0.3751        & 3.1621 \\\bottomrule
		\multicolumn{6}{l}{$^{\text{a}}$ $\zeta=-8$ and $\xi=60$}\\
            \multicolumn{6}{l}{$^{\text{b}}$ STD = standard deviation}\\
            \multicolumn{6}{l}{$^{\text{c}}$ CI = Bootstrap Confidence Interval}  
	\end{tabular*}
\end{table}

\begin{figure}[!t]
	\centering
	\includegraphics[width=5.0in]{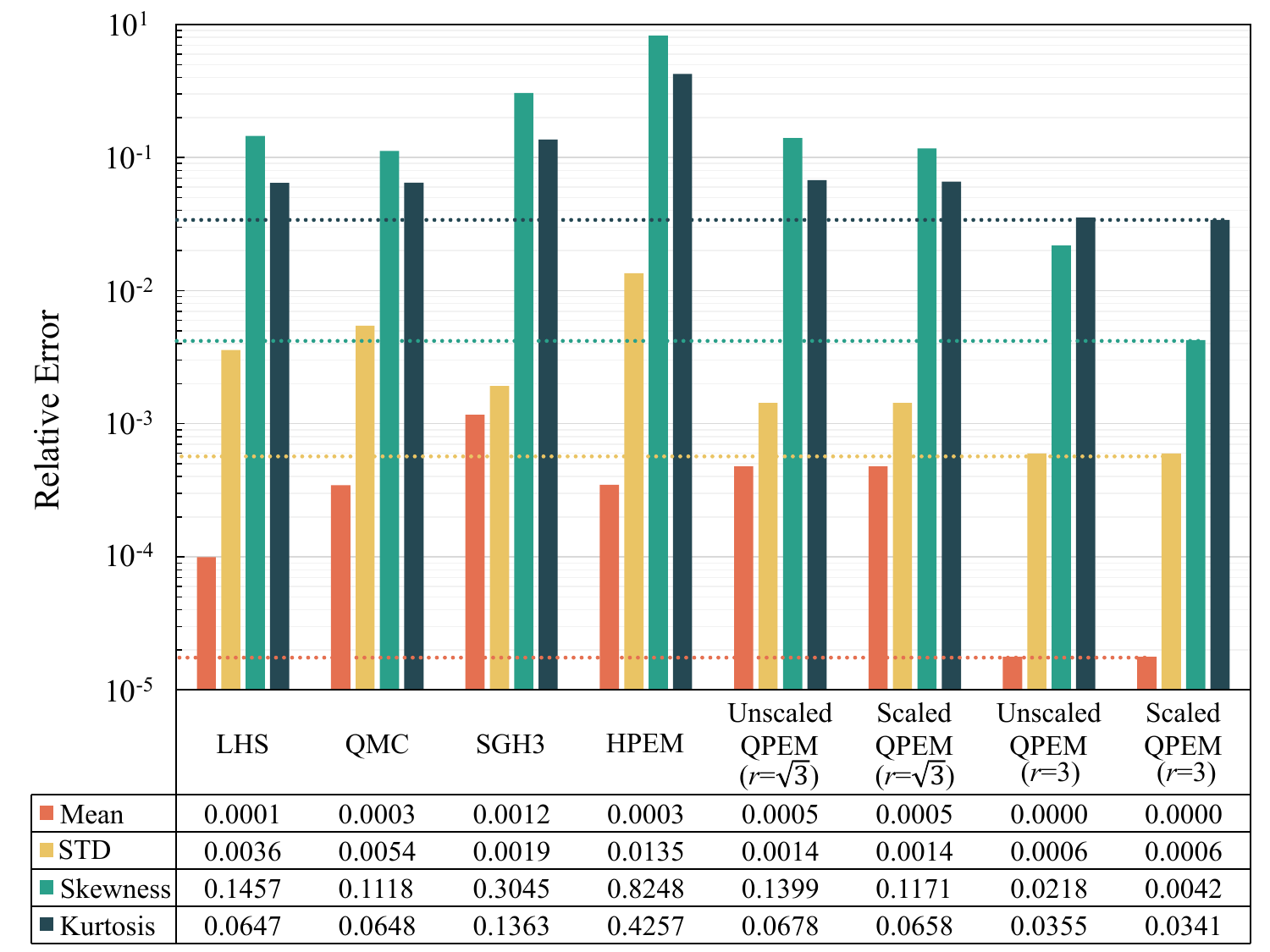}
	\caption{Relative errors for a simply supported beam example.}
	\label{fig:ex6_results}
\end{figure}

\subsection{Nonlinear finite element analysis – Simply supported beam with material nonlinearity} 
A simply supported beam involving material nonlinearity is analyzed by a FEA in this case. The beam has a length of $L$=$2,000\,mm$ and a rectangular cross-section with width of $b$=$150\,mm$ and height of $h$=$200\,mm$. The beam is modeled by 200 eight-noded hexahedron solid elements and its mesh for the nonlinear FEA is shown in Fig.~\ref{fig:ex6_beam}(a). The Poisson ratio is assumed to be $\nu=0.3$. The material nonlinearity follows an elastoplastic bi-linear model with yield strain of $\varepsilon_y=0.00125$ and yield stress of $f_y=250 \, N/mm^2$, as shown in Fig.~\ref{fig:ex6_beam}(b), and the von Mises yield criterion is considered for the yield condition. The elastic modulus $E(x,y)$ is modeled here as a two-dimensional homogeneous normal random field, in the length-height plane, through KL expansion with 40 terms in this example ($n$=40), with mean $\mu_E$=$2\times10^5\, N/mm^2$ and standard deviation $\sigma_E$=$2\times 10^4 \, N/mm^2$. The autocorrelation function of the random field is again described by Eq.~\eqref{eq:autocorrelation}, where $l_x$=$2,000/3\,mm$ and $l_y$=$200/3\,mm$ are used. The beam is subjected to a deterministic load $P$ that increases from 0 to 1,100 $kN$. The maximum beam deflection when the maximum load of $P$=$1,100kN$ is applied is the output quantity of interest in this example.

Table~\ref{table:ex6_results} and Fig.~\ref{fig:ex6_results} show the computed estimation results by each method. MC reference solutions with $10^6$ samples and confidence intervals are exactly used and computed as in previous examples. The scaled QPEM with $r$=3 has again the best performance in this case, with the scaling parameter $\zeta$ significantly improving skewness estimation accuracy, from a 2.18\% to a 0.42\% reported relative error. Overall, the QPEM shows once more an excellent performance and superior results, in this demanding example involving material nonlinearity and a stochastic field supported by 40 random dimensions.  

\section{Conclusions}
The Quadratic Point Estimate Method (QPEM) is introduced in this work, requiring $2n^2 + 1$ sampling (sigma) points in $n$-dimensional problems for estimating the first four output moments of quantities of interest. QPEM is a new, substantial and unique addition to the family of PEMs and Unscented Transform methods that are usually described by linear or exponential increase of sigma points with increasing dimensions $n$. As such, QPEM provides a distinctive balance between estimation accuracy and computational efficiency, as showcased both theoretically and through numerous numerical examples in this work. The QPEM can exactly capture all the input moments up to fifth-order for a standard normal multivariate distribution, can accurately describe all zero odd moments up to any order, and can also incorporate scaling parameters in its formulation that can partially control sixth- and eighth-order input term errors without additional computational cost. Based on derivations and numerical experimentation, all parameters of the QPEM are provided in analytical forms in this work and the method can be readily and straightforwardly applied in all stochastic problems. Due to the required number of sampling points, QPEM can be particularly competitive in low to middle-dimensional problems, e.g., $n \leq 40$, without a need for any dimensionality-reduction, but the method also works adeptly and accurately for high-dimensional problems as well. As such, users can decide on its use based on the known computational cost and expected accuracy of methods in different dimensions. Several numerical examples are provided in this work, to show the implementation and superior performance of QPEM in comparison to variance reduction sampling techniques (LHS and QMC), sparse quadrature (SGH3), and a linear PEM (HPEM). Example features vary from 5 to 100 dimensions, closed-form output expressions and finite element analyses, linear elasticity and nonlinear mechanics, and examples with one- and two-dimensional random fields, among others. In all cases, QPEM offered exceptional performance and better solutions in relation to all other methods. Ongoing work is currently extending the use and applicability of QPEM in various diverse settings, including problems related to complex non-Gaussian inputs \cite{ko2023copula} and complete output distributions estimations. 
    
\appendix
\section{Taylor series expansion for moments}
\label{Appendix-A}
This work theoretically analyzes the output moments estimation accuracy mainly through their Taylor series expansions and a standard normal input. Therefore, the Taylor series expansion of $\mathbf{y}=\mathbf{M \left( x  \right)}$  is utilized, as:
\begin{align}
	\label{eq:A.1}
	\mathbf{y} & \approx \mathbf{M} \left( {\mathbf{\bar{x}}} \right) + \nabla \mathbf{M \left( x- \bar{x} \right)} +\frac{1}{2} \nabla^2	\mathbf{M \left(x- \bar{x}\right)}^2+\frac{1}{3!} \nabla^3	\mathbf{M \left(x- \bar{x}\right)}^3+\frac{1}{4!} \nabla^4	\mathbf{M \left(x- \bar{x}\right)}^4+ \dots \\
	\label{eq:A.2}
	&=\mathbf{M} \left( {\mathbf{\bar{x}}} \right) + \nabla \mathbf{M} \, \bm{\delta x} + \frac{1}{2}\nabla^2 \mathbf{M} \, \bm{\delta x}^2 + \frac{1}{3!}\nabla^3 \mathbf{M} \, \bm{\delta x}^3 + \frac{1}{4!}\nabla^4 \mathbf{M} \, \bm{\delta x}^4 + \dots 
\end{align}
where $\mathbf{\bar{x}}=E[\mathbf{x}]=\bm{\mu}_{\mathbf{x}}$. Taking expectation of both sides of Eq.~\eqref{eq:A.2} provides the Taylor expansion of mean $E[\mathbf{y}]$, with the odd-order moments crossed-out as zero in the case of a normal distribution input:
\begin{equation}
    \label{eq:A.3}
    E\left[ \mathbf{y} \right]=\mathbf{\bar{y}}=\mathbf{M}\left( {\mathbf{\bar{x}}} \right)+\cancel{\nabla \mathbf{M}\ E\left[ \mathbf{\bm{\delta} x} \right]}+\frac{1}{2}{{\nabla }^{2}}\mathbf{M}\ E\left[ \bm{\delta }{{\mathbf{x}}^{2}} \right]+\cancel{\frac{1}{3!}{{\nabla }^{3}}\mathbf{M}\ E\left[ \mathbf{\bm{\delta} }{{\mathbf{x}}^{3}} \right]}+\frac{1}{4!}{{\nabla }^{4}}\mathbf{M}\ E\left[ \mathbf{\bm{\delta} }{{\mathbf{x}}^{4}} \right]+\cdots 
\end{equation}

Expressing the expansion of the squared deviation with the Kronecker product, as in Eqs.~\eqref{eq:notation1}-\eqref{eq:notation2}, results in:
\begin{align}
	\label{eq:A.4}
	\left( \mathbf{y}-\mathbf{\bar{y}} \right)\otimes \left( \mathbf{y}-\mathbf{\bar{y}} \right) 
	& \equiv {{\left( \mathbf{y}-\mathbf{\bar{y}} \right)}^{2}} \nonumber \\ 
        & = \left[\left(\nabla \mathbf{M}\right)^2\right] \left\{\bm{\delta}\mathbf{x}^2\right\} \nonumber \\ 
        {} & \,\,\,\, +\frac{2}{2}\left[ \nabla \mathbf{M}\otimes {{\nabla }^{2}}\mathbf{M} \right]\left\{ \bm{\delta}{{\mathbf{x}} \otimes \bm{\delta}{\mathbf{x}}^{2}}+E\left[\bm{\delta}{\mathbf{x}}\right] \otimes E\left[\bm{\delta}{\mathbf{x}}^{2}\right] \right. \nonumber\\ 
	{} &\qquad \qquad \qquad \qquad \qquad \qquad \qquad %
	\left.-\bm{\delta x}\otimes E\left[ \bm{\delta }{{\mathbf{x}}^{2}} \right] - E\left[\bm{\delta x}\right] \otimes  \bm{\delta }{{\mathbf{x}}^{2}} \right\} \nonumber \\ 
	{} & \,\,\,\, +\frac{1}{2^2}{\left[{\left( {{\nabla }^{2}}\mathbf{M} \right)}^{2}\right]}\left\{ \bm{\delta}{\mathbf{x}}^{2} \otimes \bm{\delta}{\mathbf{x}}^{2}-\bm{\delta}{{\mathbf{x}}^{2}}\otimes E\left[ \bm{\delta }{{\mathbf{x}}^{2}} \right] \right. \nonumber\\ 
	{} &\qquad \qquad \qquad \qquad \qquad \qquad \qquad %
	\left.-E\left[\bm{\delta }{{\mathbf{x}}^{2}} \right]\otimes \bm{\delta }{{\mathbf{x}}^{2}}+E{{\left[ \bm{\delta }{{\mathbf{x}}^{2}} \right]}^{2}} \right\} \nonumber \\ 
	\, & \,\,\,\, +\frac{2}{3!}\left[ \nabla \mathbf{M}\otimes {{\nabla }^{3}}\mathbf{M} \right]\left\{\bm{\delta}{\mathbf{x}} \otimes \bm{\delta}{\mathbf{x}}^{3}\right\}+\cdots 
\end{align}
Taking expectation then results in a vectorized covariance of $\mathbf{y}$, as:
\begin{align}
	\label{eq:A.5}
	E\left[ {{\left( \mathbf{y}-\mathbf{\bar{y}} \right)}^{2}} \right]&=vec\left( E\left[ \left( \mathbf{y}-\mathbf{\bar{y}} \right){{\left( \mathbf{y}-\mathbf{\bar{y}} \right)}^{T}} \right] \right) \nonumber \\ 
	{}& ={\left[{\left( \nabla \mathbf{M} \right)}^{2}\right]}\left\{E\left[ \bm{\delta }{{\mathbf{x}}^{2}} \right]\right\}+\left[ \nabla \mathbf{M}\otimes {{\nabla }^{2}}\mathbf{M} \right] \left\{  \cancel{E\left[ {\bm{\delta}{\mathbf{x}} \otimes \bm{\delta}{\mathbf{x}}^{2}} \right]}-\cancel{E\left[ \bm{\delta }\mathbf{x} \right]\otimes E\left[ \bm{\delta }{{\mathbf{x}}^{2}} \right]} \right\} \nonumber \\ 
	{}& \,\,\,\, +\frac{1}{2^2}{\left[{\left( {{\nabla }^{2}}\mathbf{M} \right)}^{2}\right]}\left\{ E\left[\bm{\delta}{\mathbf{x}}^{2} \otimes \bm{\delta}{\mathbf{x}}^{2} \right]-E{{\left[ \bm{\delta }{{\mathbf{x}}^{2}} \right]}^{2}} \right\} \nonumber\\
	{}&\,\,\,\,+\frac{1}{3}\left[ \nabla \mathbf{M}\otimes {{\nabla }^{3}}\mathbf{M} \right] \left\{E\left[ \bm{\delta }{{\mathbf{x}}} \otimes \bm{\delta }{{\mathbf{x}}^{3}}  \right]\right\}+\cdots   
\end{align}
where $vec \left( \cdot \right)$ is an operator of the vectorization (or flattening) of a matrix.

Similarly, the vectorized third and fourth-order central moments of $\mathbf{y}$ can be expressed as:

\begin{align}
	\label{eq:A.6}
	{{\left( \mathbf{y}-\mathbf{\bar{y}} \right)}^{3}} &=\left[ {{\left( \nabla \mathbf{M} \right)}^{3}}\right]\left\{\bm{\delta }{{\mathbf{x}}^3}-3\bm{\delta }{{\mathbf{x}}^2}\otimes E\left[\bm{\delta }{{\mathbf{x}}} \right]+3\bm{\delta }{{\mathbf{x}}}\otimes E\left[\bm{\delta }{{\mathbf{x}}^2} \right]-E\left[\bm{\delta }{{\mathbf{x}}^3} \right]\right\}  \nonumber \\
        {}&\,\,\,\,  +\frac{1}{2}\left[ {{\left( \nabla \mathbf{M} \right)}^{2}}\otimes {{\nabla }^{2}}\mathbf{M} \right]\left\{3\bm{\delta }{{\mathbf{x}}^{2}} \otimes \bm{\delta }{{\mathbf{x}}^{2}}+3\bm{\delta }{{\mathbf{x}}^{2}} \otimes E\left[\bm{\delta }{{\mathbf{x}}}\right]^2-6\bm{\delta}\mathbf{x}\otimes \bm{\delta}\mathbf{x}^2 \otimes E\left[\bm{\delta}\mathbf{x}\right] \right.  \nonumber \\
        {} & \qquad \qquad \qquad %
	\left.-3\bm{\delta }{{\mathbf{x}}^{2}}\otimes E\left[\bm{\delta }{{\mathbf{x}}^{2}} \right]-3 E\left[\bm{\delta }{{\mathbf{x}}^{2}}\right] \otimes E\left[\bm{\delta }{{\mathbf{x}}}\right]^2+6\bm{\delta}\mathbf{x}\otimes E\left[\bm{\delta}\mathbf{x}^2\right] \otimes E\left[\bm{\delta}\mathbf{x}\right] \right\} \nonumber \\
        {}&\,\,\,\,  +\frac{1}{2^2}\left[ {{ \nabla \mathbf{M}}}\otimes \left({\nabla }^{2}\mathbf{M}\right)^2 \right]\left\{3\bm{\delta }{{\mathbf{x}}} \otimes \bm{\delta }{{\mathbf{x}}^{2}}\otimes \bm{\delta }{{\mathbf{x}}^{2}}+3\bm{\delta }{{\mathbf{x}}} \otimes E\left[\bm{\delta }{{\mathbf{x}}^{2}}\right]^2-6\bm{\delta}\mathbf{x}\otimes \bm{\delta}\mathbf{x}^2 \otimes E\left[\bm{\delta}\mathbf{x}^2\right] \right.  \nonumber \\
        {} & \qquad \qquad \qquad %
	\left.-3 E\left[\bm{\delta }{{\mathbf{x}}}\right] \otimes \bm{\delta }{{\mathbf{x}}^{2}}\otimes \bm{\delta }{{\mathbf{x}}^{2}}-3 E\left[\bm{\delta }{{\mathbf{x}}}\right] \otimes E\left[\bm{\delta }{{\mathbf{x}}^{2}}\right]^2+6 E\left[\bm{\delta}\mathbf{x}\right]\otimes \bm{\delta}\mathbf{x}^2 \otimes E\left[\bm{\delta}\mathbf{x}^2\right] \right\} \nonumber \\
        {}&\,\,\,\,  +\frac{1}{3!}\left[ {{\left( \nabla \mathbf{M} \right)^2}}\otimes {\nabla }^{3}\mathbf{M}\right]\left\{3\bm{\delta }{{\mathbf{x}}}^2 \otimes \bm{\delta }{{\mathbf{x}}^{3}}+3 E\left[\bm{\delta }{{\mathbf{x}}}\right]^2 \otimes \bm{\delta }{{\mathbf{x}}^3}-6\bm{\delta}\mathbf{x}\otimes \bm{\delta}\mathbf{x}^3\otimes E\left[\bm{\delta}\mathbf{x}\right]  \right.  \nonumber \\
        {} & \qquad \qquad \qquad %
	\left.-3\bm{\delta }{{\mathbf{x}}}^2 \otimes E\left[\bm{\delta }{{\mathbf{x}}^{3}}\right]-3 E\left[\bm{\delta }{{\mathbf{x}}}\right]^2 \otimes E\left[\bm{\delta }{{\mathbf{x}}^3}\right] +6 \bm{\delta}\mathbf{x}\otimes E\left[\bm{\delta}\mathbf{x}^3\right]\otimes E\left[\bm{\delta}\mathbf{x}\right]  \right\} \nonumber \\
	{} &  \,\,\,\,+\frac{1}{{{2}^{3}}}{\left[{\left( {{\nabla }^{2}}\mathbf{M} \right)}^{3}\right]}\left\{ \bm{\delta }{{\mathbf{x}}^{2}} \otimes \bm{\delta }{{\mathbf{x}}^{2}}\otimes \bm{\delta }{{\mathbf{x}}^{2}}-3\bm{\delta }{{\mathbf{x}}^{2}}\otimes \bm{\delta }{{\mathbf{x}}^{2}}\otimes E\left[ \bm{\delta }{{\mathbf{x}}^{2}} \right] \right. \nonumber \\
	{} & \qquad \qquad \qquad \qquad \qquad \qquad \qquad %
	\left.+3\bm{\delta }{{\mathbf{x}}^{2}}\otimes E{{\left[ \bm{\delta }{{\mathbf{x}}^{2}} \right]}^{2}}-E{{\left[ \bm{\delta }{{\mathbf{x}}^{2}} \right]}^{3}} \right\}+\cdots
\end{align}

\begin{align} 
	E\left[ {{\left( \mathbf{y}-\mathbf{\bar{y}} \right)}^{3}} \right] &=\left[ {{\left( \nabla \mathbf{M} \right)}^{3}}\right]\left\{\cancel{ E\left[\bm{\delta }{{\mathbf{x}}^3}\right]}-3\cancel{E\left[\bm{\delta }{{\mathbf{x}}^2}\right]\otimes E\left[\bm{\delta }{{\mathbf{x}}} \right]}+3\cancel{E\left[\bm{\delta }{{\mathbf{x}}}\right]\otimes E\left[\bm{\delta }{{\mathbf{x}}^2} \right]}-\cancel{E\left[\bm{\delta }{{\mathbf{x}}^3} \right]}\right\}  \nonumber \\
        {}&\,\,\,\,  +\frac{1}{2}\left[ {{\left( \nabla \mathbf{M} \right)}^{2}}\otimes {{\nabla }^{2}}\mathbf{M} \right] \left\{3 E\left[\bm{\delta }{{\mathbf{x}}^{2}} \otimes \bm{\delta }{{\mathbf{x}}^{2}}\right]+3\cancel{E\left[\bm{\delta }{{\mathbf{x}}^{2}}\right] \otimes E\left[\bm{\delta }{{\mathbf{x}}}\right]^2} \right.  \nonumber \\
        {} & \qquad \qquad \qquad %
	-6\cancel{E\left[\bm{\delta}\mathbf{x}\otimes \bm{\delta}\mathbf{x}^2\right] \otimes E\left[\bm{\delta}\mathbf{x}\right]}-3E\left[\bm{\delta }{{\mathbf{x}}^{2}}\right]\otimes E\left[\bm{\delta }{{\mathbf{x}}^{2}} \right]\nonumber \\
        {}& \qquad \qquad \qquad %
        \left. -3 \cancel{E\left[\bm{\delta }{{\mathbf{x}}^{2}}\right] \otimes E\left[\bm{\delta }{{\mathbf{x}}}\right]^2} +6 \cancel{E\left[\bm{\delta}\mathbf{x}\right]\otimes E\left[\bm{\delta}\mathbf{x}^2\right] \otimes E\left[\bm{\delta}\mathbf{x}\right]} \right\} \nonumber \\
        {}&\,\,\,\,  +\frac{1}{2^2}\left[ {{ \nabla \mathbf{M}}}\otimes \left({\nabla }^{2}\mathbf{M}\right)^2 \right]\left\{3\cancel{E\left[\bm{\delta }{{\mathbf{x}}} \otimes \bm{\delta }{{\mathbf{x}}^{2}}\otimes \bm{\delta }{{\mathbf{x}}^{2}}\right]}+3\cancel{E\left[\bm{\delta }{{\mathbf{x}}}\right] \otimes E\left[\bm{\delta }{{\mathbf{x}}^{2}}\right]^2} \right.  \nonumber \\
        {} & \qquad \qquad \qquad %
	-6\cancel{E\left[\bm{\delta}\mathbf{x}\otimes \bm{\delta}\mathbf{x}^2\right] \otimes E\left[\bm{\delta}\mathbf{x}^2\right]} -3 \cancel{E\left[\bm{\delta }{{\mathbf{x}}}\right] \otimes E\left[\bm{\delta }{{\mathbf{x}}^{2}}\otimes \bm{\delta }{{\mathbf{x}}^{2}}\right]} \nonumber \\
        {} & \qquad \qquad \qquad %
        \left.-3 \cancel{E\left[\bm{\delta }{{\mathbf{x}}}\right] \otimes E\left[\bm{\delta }{{\mathbf{x}}^{2}}\right]^2}+6 \cancel{E\left[\bm{\delta}\mathbf{x}\right]\otimes E\left[\bm{\delta}\mathbf{x}^2\right] \otimes E\left[\bm{\delta}\mathbf{x}^2\right]} \right\} \nonumber
\end{align}
\begin{align}
    \label{eq:A.7}
        {}&\,\,\,\,  +\frac{1}{3!}\left[ {{\left( \nabla \mathbf{M} \right)^2}}\otimes {\nabla }^{3}\mathbf{M}\right]\left\{3\cancel{E\left[\bm{\delta }{{\mathbf{x}}}^2 \otimes \bm{\delta }{{\mathbf{x}}^{3}}\right]}+3 \cancel{E\left[\bm{\delta }{{\mathbf{x}}}\right]^2 \otimes E\left[\bm{\delta }{{\mathbf{x}}^3}\right]} \right.  \nonumber \\
        {} & \qquad \qquad \qquad %
	-6 \cancel{E\left[\bm{\delta}\mathbf{x}\otimes \bm{\delta}\mathbf{x}^3\right]\otimes E\left[\bm{\delta}\mathbf{x}\right]} -3\cancel{E\left[\bm{\delta }{{\mathbf{x}}}^2\right] \otimes E\left[\bm{\delta }{{\mathbf{x}}^{3}}\right]} \nonumber \\
        {}& \qquad \qquad \qquad %
        \left.-3 \cancel{E\left[\bm{\delta }{{\mathbf{x}}}\right]^2 \otimes E\left[\bm{\delta }{{\mathbf{x}}^3}\right]} +6 \cancel{E\left[\bm{\delta}\mathbf{x}\right]\otimes E\left[\bm{\delta}\mathbf{x}^3\right]\otimes E\left[\bm{\delta}\mathbf{x}\right]} \right\} \nonumber \\ 
	{} &  \,\,\,\,+\frac{1}{{{2}^{3}}}{\left[{\left( {{\nabla }^{2}}\mathbf{M} \right)}^{3}\right]}\left\{E\left[\bm{\delta }{{\mathbf{x}}^{2}} \otimes \bm{\delta }{{\mathbf{x}}^{2}}\otimes \bm{\delta }{{\mathbf{x}}^{2}}\right]-3E\left[\bm{\delta }{{\mathbf{x}}^{2}}\otimes \bm{\delta }{{\mathbf{x}}^{2}}\right]\otimes E\left[ \bm{\delta }{{\mathbf{x}}^{2}} \right] +2E{{\left[ \bm{\delta }{{\mathbf{x}}^{2}} \right]}^{3}} \right\} \nonumber \\
	{} & \,\,\,\,+\cdots
\end{align}

\begin{align} 
	{{\left( \mathbf{y}-\mathbf{\bar{y}} \right)}^{4}}&={\left[{\left( \nabla \mathbf{M} \right)}^{4}\right]} \left\{\bm{\delta }{{\mathbf{x}}}^4-4\bm{\delta }{{\mathbf{x}}}^3 \otimes E\left[\bm{\delta }{{\mathbf{x}}}\right]+6\bm{\delta}\mathbf{x}^2 \otimes E\left[\bm{\delta}\mathbf{x}\right]^2-4\bm{\delta }{{\mathbf{x}}}\otimes E\left[\bm{\delta }{{\mathbf{x}}}\right]^3 +E\left[\bm{\delta}\mathbf{x}\right]^4\right\} \nonumber \\
	\,& \,\,\,\,+\frac{1}{2}\left[ {{\left( \nabla \mathbf{M} \right)}^{3}}\otimes {{\nabla }^{2}}\mathbf{M} \right]\left\{4\bm{\delta }{{\mathbf{x}}}^3\otimes \bm{\delta }{{\mathbf{x}}^{2}} -4\bm{\delta }{{\mathbf{x}}}^3\otimes E \left[\bm{\delta }{{\mathbf{x}}^{2}}\right]-12\bm{\delta}\mathbf{x}^2\otimes \bm{\delta}\mathbf{x}^2 \otimes E\left[\bm{\delta}\mathbf{x}\right] \right.  \nonumber\\ 
        \,& \qquad \qquad + 12\bm{\delta}\mathbf{x}^2\otimes E \left[\bm{\delta}\mathbf{x}^2\right] \otimes E\left[\bm{\delta}\mathbf{x}\right]+12\bm{\delta}\mathbf{x}\otimes \bm{\delta}\mathbf{x}^2 \otimes E\left[\bm{\delta}\mathbf{x}\right]^2 \nonumber\\
        \,& \qquad \qquad  \left. -12\bm{\delta}\mathbf{x}\otimes E\left[\bm{\delta}\mathbf{x}\right]^2\otimes E\left[\bm{\delta}\mathbf{x}^2\right] -4 \bm{\delta}\mathbf{x}^2 \otimes E\left[\bm{\delta}\mathbf{x}\right]^3 +4E \left[\bm{\delta}\mathbf{x}\right]^3 \otimes E\left[\bm{\delta}\mathbf{x}^2\right] \right\} \nonumber \\
        \,& \,\,\,\,+\frac{1}{3!}\left[ {{\left( \nabla \mathbf{M} \right)}^{3}}\otimes {{\nabla }^{3}}\mathbf{M} \right]\left\{4\bm{\delta }{{\mathbf{x}}}^3\otimes \bm{\delta }{{\mathbf{x}}^{3}} -4\bm{\delta }{{\mathbf{x}}}^3\otimes E \left[\bm{\delta }{{\mathbf{x}}^{3}}\right]-12\bm{\delta}\mathbf{x}^2\otimes \bm{\delta}\mathbf{x}^3 \otimes E\left[\bm{\delta}\mathbf{x}\right] \right.  \nonumber\\ 
        \,& \qquad \qquad + 12\bm{\delta}\mathbf{x}^2\otimes E \left[\bm{\delta}\mathbf{x}^3\right] \otimes E\left[\bm{\delta}\mathbf{x}\right]+12\bm{\delta}\mathbf{x}\otimes \bm{\delta}\mathbf{x}^3 \otimes E\left[\bm{\delta}\mathbf{x}\right]^2 \nonumber\\
        \,& \qquad \qquad  \left. -12\bm{\delta}\mathbf{x}\otimes E\left[\bm{\delta}\mathbf{x}\right]^2\otimes E\left[\bm{\delta}\mathbf{x}^3\right] -4 \bm{\delta}\mathbf{x}^3 \otimes E\left[\bm{\delta}\mathbf{x}\right]^3 +4E \left[\bm{\delta}\mathbf{x}\right]^3 \otimes E\left[\bm{\delta}\mathbf{x}^3\right] \right\} \nonumber \\
        \,& \,\,\,\,+\frac{1}{2^2}\left[ {{\left( \nabla \mathbf{M} \right)}^{2}}\otimes {{\left( \nabla^2 \mathbf{M} \right)}^{2}} \right]\left\{6\bm{\delta }{{\mathbf{x}}}^2 \otimes \bm{\delta }{{\mathbf{x}}^{2}}\otimes \bm{\delta }{{\mathbf{x}}^{2}} -12 \bm{\delta}\mathbf{x}\otimes \bm{\delta}\mathbf{x}^2 \otimes \bm{\delta}\mathbf{x}^2 \otimes E\left[\bm{\delta}\mathbf{x}\right] \right.  \nonumber\\ 
        \,& \qquad -12 \bm{\delta}\mathbf{x}^2\otimes \bm{\delta}\mathbf{x}^2 \otimes E\left[\bm{\delta}\mathbf{x}^2\right] -12 \bm{\delta}\mathbf{x}\otimes E\left[\bm{\delta}\mathbf{x}\right] \otimes E\left[\bm{\delta}\mathbf{x}^2\right]^2  \nonumber\\
        \,& \qquad \qquad  -12 \bm{\delta}\mathbf{x}^2\otimes E\left[\bm{\delta}\mathbf{x}\right]^2 \otimes E\left[\bm{\delta}\mathbf{x}^2\right] +24 \bm{\delta}\mathbf{x}\otimes \bm{\delta}\mathbf{x}^2 \otimes E\left[\bm{\delta}\mathbf{x}\right] \otimes E\left[\bm{\delta}\mathbf{x}^2\right]  \nonumber \\
        \,& \qquad \qquad \left. +6 \bm{\delta}\mathbf{x}^2\otimes \bm{\delta}\mathbf{x}^2 \otimes E\left[\bm{\delta}\mathbf{x}\right]^2 +6 \bm{\delta}\mathbf{x}^2\otimes E\left[\bm{\delta}\mathbf{x}^2\right]^2  +6 E\left[\bm{\delta}\mathbf{x}\right]^2 \otimes E\left[\bm{\delta}\mathbf{x}^2\right]^2 \right\} \nonumber\\
        \,& \,\,\,\,+\frac{1}{2^3}\left[ {\nabla \mathbf{M}}\otimes \left({{\nabla }^{2}}\mathbf{M}\right)^3 \right]\left\{4\bm{\delta }{{\mathbf{x}}}\otimes \bm{\delta }{{\mathbf{x}}^{2}}\otimes \bm{\delta }{{\mathbf{x}}^{2}}\otimes \bm{\delta }{{\mathbf{x}}^{2}} -4\bm{\delta }{{\mathbf{x}}}\otimes E \left[\bm{\delta }{{\mathbf{x}}^{2}}\right]^3 \right.  \nonumber\\ 
        \,& \qquad \qquad -12\bm{\delta}\mathbf{x}\otimes \bm{\delta}\mathbf{x}^2\otimes \bm{\delta}\mathbf{x}^2 \otimes E\left[\bm{\delta}\mathbf{x}^2\right]+ 12\bm{\delta}\mathbf{x}^2\otimes\bm{\delta}\mathbf{x}^2\otimes E \left[\bm{\delta}\mathbf{x}\right] \otimes E\left[\bm{\delta}\mathbf{x}^2\right]\nonumber\\
        \,& \qquad \qquad -12\bm{\delta}\mathbf{x}^2\otimes E\left[\bm{\delta}\mathbf{x}\right]\otimes E\left[\bm{\delta}\mathbf{x}^2\right]^2 +12\bm{\delta}\mathbf{x}\otimes \bm{\delta}\mathbf{x}^2 \otimes E\left[\bm{\delta}\mathbf{x}^2\right]^2 \nonumber \\
        \,& \qquad \qquad  \left. -4 \bm{\delta}\mathbf{x}^2 \otimes\bm{\delta}\mathbf{x}^2 \otimes\bm{\delta}\mathbf{x}^2 \otimes E\left[\bm{\delta}\mathbf{x}\right] +4E \left[\bm{\delta}\mathbf{x}\right] \otimes E\left[\bm{\delta}\mathbf{x}^2\right]^3 \right\} \nonumber\\
        \,& \,\,\,\,+\frac{1}{4!}\left[ \left( \nabla \mathbf{M} \right)^3\otimes {{\nabla }^{4}}\mathbf{M}\right]\left\{4 \bm{\delta}\mathbf{x}^3 \otimes \bm{\delta}\mathbf{x}^4 -4 \bm{\delta}\mathbf{x}^3 \otimes E\left [\bm{\delta}\mathbf{x}^4\right] -4 \bm{\delta}\mathbf{x}^4 \otimes  E\left[\bm{\delta}\mathbf{x}\right]^3 \right.  \nonumber\\
        \,& \qquad \qquad -12 \bm{\delta}\mathbf{x}^2 \otimes \bm{\delta}\mathbf{x}^4 \otimes E \left[\bm{\delta}\mathbf{x}\right] +12\bm{\delta}\mathbf{x} \otimes \bm{\delta}\mathbf{x}^4 \otimes E\left[\bm{\delta}\mathbf{x}^2\right] +12 \bm{\delta}\mathbf{x}^2 \otimes E\left[\bm{\delta}\mathbf{x}\right] \otimes E\left[\bm{\delta}\mathbf{x}^4\right] \nonumber \\
        \,& \qquad \qquad \left. -12 \bm{\delta}\mathbf{x}\otimes E\left[\bm{\delta}\mathbf{x}\right]^2 \otimes E\left[\bm{\delta}\mathbf{x}^4\right] +4 E\left[\bm{\delta}\mathbf{x}\right]^3\otimes E\left[\bm{\delta}\mathbf{x}^4\right] \right \} \nonumber\\
        \,& \,\,\,\,+\frac{1}{2!3!}\left[ {{\left( \nabla \mathbf{M} \right)}^{2}}\otimes {{\nabla }^{2}}\mathbf{M} \otimes  {{\nabla }^{3}}\mathbf{M} \right]\left\{12\bm{\delta }{{\mathbf{x}}}^2 \otimes \bm{\delta }{{\mathbf{x}}}^2 \otimes \bm{\delta }{{\mathbf{x}}}^3 -12 \bm{\delta }{{\mathbf{x}}}^2 \otimes \bm{\delta }{{\mathbf{x}}}^3 \otimes E\left[\bm{\delta }{{\mathbf{x}}}^2\right]  \right.  \nonumber\\ 
    \,& \qquad \qquad -12 \bm{\delta }{{\mathbf{x}}}^2 \otimes \bm{\delta }{{\mathbf{x}}}^2 \otimes E\left[\bm{\delta }{{\mathbf{x}}}^3\right] -24 {\bm{\delta}\mathbf{x}} \otimes {\bm{\delta}\mathbf{x}}^2 \otimes {\bm{\delta}\mathbf{x}}^3 \otimes E\left[{\bm{\delta}\mathbf{x}}\right] \nonumber \\
    \,& \qquad \qquad +24 {\bm{\delta}\mathbf{x}}\otimes {\bm{\delta}\mathbf{x}}^3 \otimes E\left[{\bm{\delta}\mathbf{x}}\right] \otimes E\left[{\bm{\delta}\mathbf{x}}^2\right]+24{\bm{\delta}\mathbf{x}}\otimes {\bm{\delta}\mathbf{x}}^2 \otimes E\left[{\bm{\delta}\mathbf{x}}\right] \otimes E\left[{\bm{\delta}\mathbf{x}}^3\right] \nonumber \\
    \,& \qquad \qquad +12 {\bm{\delta}\mathbf{x}}^2 \otimes E\left[{\bm{\delta}\mathbf{x}}^2\right] \otimes E\left[{\bm{\delta}\mathbf{x}}^3 \right] -12 {\bm{\delta}\mathbf{x}}^2 \otimes E \left[{\bm{\delta}\mathbf{x}}\right]^2 E\left[{\bm{\delta}\mathbf{x}}^3\right] \nonumber \\
    \,& \qquad \qquad \left.-24 {\bm{\delta}\mathbf{x}} \otimes E\left[{\bm{\delta}\mathbf{x}}\right] \otimes E\left[{\bm{\delta}\mathbf{x}}^2 \right] \otimes E\left[{\bm{\delta}\mathbf{x}}^3 \right] +12 E\left[{\bm{\delta}\mathbf{x}} \right]^2 \otimes E\left[{\bm{\delta}\mathbf{x}}^2 \right] \otimes E\left[{\bm{\delta}\mathbf{x}}^3 \right] \right\} \nonumber
\end{align}
\begin{align}     
    \label{eq:A.8}
    \,& \,\,\,\,+\frac{1}{2^4}\left[ {{\left( \nabla^2 \mathbf{M} \right)}^{4}}\right] \left\{ {\bm{\delta}\mathbf{x}}^2\otimes {\bm{\delta}\mathbf{x}}^2 \otimes {\bm{\delta}\mathbf{x}}^2 \otimes {\bm{\delta}\mathbf{x}}^2 -4 {\bm{\delta}\mathbf{x}}^2 \otimes E \left[{\bm{\delta}\mathbf{x}}^2\right]^3  \right. \nonumber \\
    \,& \qquad  \qquad \left. +6 {\bm{\delta}\mathbf{x}}^2 \otimes {\bm{\delta}\mathbf{x}}^2 \otimes E\left[{\bm{\delta}\mathbf{x}}^2\right]^2-4 {\bm{\delta}\mathbf{x}}^2 \otimes {\bm{\delta}\mathbf{x}}^2 \otimes{\bm{\delta}\mathbf{x}}^2 \otimes E \left[{\bm{\delta}\mathbf{x}}^2\right]  + E\left[{\bm{\delta}\mathbf{x}}^2\right]^4 \right\} +\cdots
\end{align}

\begin{align} 
	E\left[{{\left( \mathbf{y}-\mathbf{\bar{y}} \right)}^{4}}\right]&={\left[{\left( \nabla \mathbf{M} \right)}^{4}\right]} \left\{E\left[\bm{\delta }{{\mathbf{x}}}^4\right] -4 \cancel{E\left[\bm{\delta }{{\mathbf{x}}}^3 \right]\otimes E\left[\bm{\delta }{{\mathbf{x}}}\right]}\right. \nonumber \\
        \,& \qquad \qquad \left.+6\cancel{E\left[\bm{\delta}\mathbf{x}^2\right]\otimes E\left[\bm{\delta}\mathbf{x}\right]^2} -4\cancel{E\left[\bm{\delta }{{\mathbf{x}}}\right] \otimes E\left[\bm{\delta }{{\mathbf{x}}}\right]^3} +\cancel{E\left[\bm{\delta}\mathbf{x}\right]^4}\right\} \nonumber \\
	\,& \,\,\,\,+\frac{1}{2}\left[ {{\left( \nabla \mathbf{M} \right)}^{3}}\otimes {{\nabla }^{2}}\mathbf{M} \right]\left\{4 \cancel{E\left[\bm{\delta }{{\mathbf{x}}}^3\otimes \bm{\delta }{{\mathbf{x}}^{2}}\right]} -4\cancel{E\left[\bm{\delta }{{\mathbf{x}}}^3 \right] \otimes E \left[\bm{\delta }{{\mathbf{x}}^{2}}\right]}\right.  \nonumber\\ 
        \,& \qquad \qquad  -12 \cancel{E\left[\bm{\delta}\mathbf{x}^2\otimes \bm{\delta}\mathbf{x}^2\right] \otimes E\left[\bm{\delta}\mathbf{x}\right]} + 12\cancel{E\left[\bm{\delta}\mathbf{x}\right] \otimes E\left[\bm{\delta}\mathbf{x}^2\right]^2} \nonumber\\
        \,& \qquad \qquad  \left. +12\cancel{E\left[\bm{\delta}\mathbf{x}\otimes \bm{\delta}\mathbf{x}^2\right] \otimes E\left[\bm{\delta}\mathbf{x}\right]^2} -12\cancel{E\left[\bm{\delta}\mathbf{x}\right]^3\otimes E\left[\bm{\delta}\mathbf{x}^2\right]} \right\} \nonumber \\
        \,& \,\,\,\,+\frac{1}{3!}\left[ {{\left( \nabla \mathbf{M} \right)}^{3}}\otimes {{\nabla }^{3}}\mathbf{M} \right]\left\{4E\left[\bm{\delta }{{\mathbf{x}}}^3\otimes \bm{\delta }{{\mathbf{x}}^{3}}\right] -4\cancel{E\left[\bm{\delta }{{\mathbf{x}}}^3\right]\otimes E \left[\bm{\delta }{{\mathbf{x}}^{3}}\right]} \right. \nonumber \\
        \,& \qquad \qquad -12 \cancel{E\left[\bm{\delta}\mathbf{x}^2\otimes \bm{\delta}\mathbf{x}^3 \right]\otimes E\left[\bm{\delta}\mathbf{x}\right]} + 12 \cancel{E\left[ \bm{\delta}\mathbf{x}^2\right]\otimes E \left[\bm{\delta}\mathbf{x}^3\right] \otimes E\left[\bm{\delta}\mathbf{x}\right]} \nonumber\\ 
        \,& \qquad \qquad \left.+12 \cancel{\left[\bm{\delta}\mathbf{x}\otimes \bm{\delta}\mathbf{x}^3 \right]\otimes E\left[\bm{\delta}\mathbf{x}\right]^2}-12 \cancel{E\left[\bm{\delta}\mathbf{x}\right]^3\otimes E\left[\bm{\delta}\mathbf{x}^3\right]} \right\} \nonumber\\
        \,& \,\,\,\,+\frac{1}{2^2}\left[ {{\left( \nabla \mathbf{M} \right)}^{2}}\otimes {{\left( \nabla^2 \mathbf{M} \right)}^{2}} \right]\left\{6E\left[\bm{\delta }{{\mathbf{x}}}^2 \otimes \bm{\delta }{{\mathbf{x}}^{2}}\otimes \bm{\delta }{{\mathbf{x}}^{2}}\right] \right.  \nonumber\\ 
        \,& \qquad \qquad -12 \cancel{E\left[\bm{\delta}\mathbf{x}\otimes \bm{\delta}\mathbf{x}^2 \otimes \bm{\delta}\mathbf{x}^2\right] \otimes E\left[\bm{\delta}\mathbf{x}\right]} -12 E\left[\bm{\delta}\mathbf{x}^2\otimes \bm{\delta}\mathbf{x}^2\right] \otimes E\left[\bm{\delta}\mathbf{x}^2\right]  \nonumber\\
        \,& \qquad \qquad  -24 \cancel{E\left[\bm{\delta}\mathbf{x}\right]^2 \otimes E\left[\bm{\delta}\mathbf{x}^2\right]^2} +24 \cancel{E\left[\bm{\delta}\mathbf{x}\otimes \bm{\delta}\mathbf{x}^2\right] \otimes E\left[\bm{\delta}\mathbf{x}\right] \otimes E\left[\bm{\delta}\mathbf{x}^2\right]}  \nonumber \\
        \,& \qquad \qquad \left. +6 \cancel{E\left[\bm{\delta}\mathbf{x}^2\otimes \bm{\delta}\mathbf{x}^2 \right]\otimes E\left[\bm{\delta}\mathbf{x}\right]^2} +6 \cancel{E\left[\bm{\delta}\mathbf{x}\right]^2 \otimes E\left[\bm{\delta}\mathbf{x}^2\right]^2} +6 E\left[\bm{\delta}\mathbf{x}^2\right]^3 \right\} \nonumber\\
        \,& \,\,\,\,+\frac{1}{2^3}\left[ {\nabla \mathbf{M}}\otimes \left({{\nabla }^{2}}\mathbf{M}\right)^3 \right]\left\{4 \cancel{E\left[\bm{\delta }{{\mathbf{x}}}\otimes \bm{\delta }{{\mathbf{x}}^{2}}\otimes \bm{\delta }{{\mathbf{x}}^{2}}\otimes \bm{\delta }{{\mathbf{x}}^{2}}\right]} -4\cancel{E\left[\bm{\delta }{{\mathbf{x}}}\right]\otimes E \left[\bm{\delta }{{\mathbf{x}}^{2}}\right]^3} \right.  \nonumber\\ 
        \,& \qquad \qquad -12\cancel{E\left[\bm{\delta}\mathbf{x}\otimes \bm{\delta}\mathbf{x}^2\otimes \bm{\delta}\mathbf{x}^2 \right] \otimes E\left[\bm{\delta}\mathbf{x}^2\right]}+ 12\cancel{E\left[\bm{\delta}\mathbf{x}^2\otimes\bm{\delta}\mathbf{x}^2\right]\otimes E \left[\bm{\delta}\mathbf{x}\right] \otimes E\left[\bm{\delta}\mathbf{x}^2\right]}\nonumber\\
        \,& \qquad \qquad -8\cancel{E\left[\bm{\delta}\mathbf{x}\right]\otimes E\left[\bm{\delta}\mathbf{x}^2\right]^3} +12\cancel{\left[\bm{\delta}\mathbf{x}\otimes \bm{\delta}\mathbf{x}^2\right]\otimes E\left[\bm{\delta}\mathbf{x}^2\right]^2} \nonumber \\
        \,& \qquad \qquad  \left. -4 \cancel{E\left[\bm{\delta}\mathbf{x}^2 \otimes\bm{\delta}\mathbf{x}^2 \otimes\bm{\delta}\mathbf{x}^2\right] \otimes E\left[\bm{\delta}\mathbf{x}\right]} \right\} \nonumber \\
        \,& \,\,\,\,+\frac{1}{4!}\left[ \left( \nabla \mathbf{M} \right)^3\otimes {{\nabla }^{4}}\mathbf{M}\right]\left\{4 \cancel{E\left[\bm{\delta}\mathbf{x}^3 \otimes \bm{\delta}\mathbf{x}^4\right]} -4 \cancel{E\left[\bm{\delta}\mathbf{x}^3\right]\otimes E\left [\bm{\delta}\mathbf{x}^4\right]} \right.  \nonumber\\ 
        \,& \qquad \qquad -12 \cancel{E\left[\bm{\delta}\mathbf{x}^2 \otimes \bm{\delta}\mathbf{x}^4 \right] \otimes E \left[\bm{\delta}\mathbf{x}\right]} +12 \cancel{E\left[\bm{\delta}\mathbf{x} \otimes \bm{\delta}\mathbf{x}^4 \right]\otimes E\left[\bm{\delta}\mathbf{x}^2\right]}  \nonumber \\
        \,& \qquad \qquad \left. +12 \cancel{E\left[\bm{\delta}\mathbf{x}^2\right] \otimes E\left[\bm{\delta}\mathbf{x}\right] \otimes E\left[\bm{\delta}\mathbf{x}^4\right]}-12 \cancel{E\left[\bm{\delta}\mathbf{x}\right]\otimes E\left[\bm{\delta}\mathbf{x}\right]^2 \otimes E\left[\bm{\delta}\mathbf{x}^4\right]} \right \}\nonumber\\
        \,& \,\,\,\,+\frac{1}{2!3!}\left[ {{\left( \nabla \mathbf{M} \right)}^{2}}\otimes {{\nabla }^{2}}\mathbf{M} \otimes  {{\nabla }^{3}}\mathbf{M} \right]\left\{12\cancel{E\left[\bm{\delta }{{\mathbf{x}}}^2 \otimes \bm{\delta }{{\mathbf{x}}}^2 \otimes \bm{\delta }{{\mathbf{x}}}^3\right]}   \right.  \nonumber\\ 
    \,& \qquad \qquad -12 \cancel{E\left[\bm{\delta }{{\mathbf{x}}}^2 \otimes \bm{\delta }{{\mathbf{x}}}^3 \right]\otimes E\left[\bm{\delta }{{\mathbf{x}}}^2\right]}-12 \cancel{E\left[\bm{\delta }{{\mathbf{x}}}^2 \otimes \bm{\delta }{{\mathbf{x}}}^2\right] \otimes E\left[\bm{\delta }{{\mathbf{x}}}^3\right]} \nonumber \\
    \,& \qquad \qquad -24 \cancel{E\left[{\bm{\delta}\mathbf{x}} \otimes {\bm{\delta}\mathbf{x}}^2 \otimes {\bm{\delta}\mathbf{x}}^3\right] \otimes E\left[{\bm{\delta}\mathbf{x}}\right]} +24 \cancel{E\left[{\bm{\delta}\mathbf{x}}\otimes {\bm{\delta}\mathbf{x}}^3 \right]\otimes E\left[{\bm{\delta}\mathbf{x}}\right] \otimes E\left[{\bm{\delta}\mathbf{x}}^2\right]} \nonumber \\
    \,& \qquad \qquad +24\cancel{E\left[{\bm{\delta}\mathbf{x}}\otimes {\bm{\delta}\mathbf{x}}^2\right] \otimes E\left[{\bm{\delta}\mathbf{x}}\right] \otimes E\left[{\bm{\delta}\mathbf{x}}^3\right]} +12 \cancel{E\left[{\bm{\delta}\mathbf{x}}^2\right]^2 \otimes E\left[{\bm{\delta}\mathbf{x}}^3 \right]} \nonumber \\
    \,& \qquad \qquad \left. -24 \cancel{E\left[{\bm{\delta}\mathbf{x}}\right]^2 \otimes E\left[{\bm{\delta}\mathbf{x}}^2 \right] \otimes E\left[{\bm{\delta}\mathbf{x}}^3 \right]}  \right\} \nonumber
\end{align}

\begin{align}
	\label{eq:A.9}    
    \,& \,\,\,\,+\frac{1}{2^4}\left[ {{\left( \nabla^2 \mathbf{M} \right)}^{4}}\right] \left\{ E\left[{\bm{\delta}\mathbf{x}}^2\otimes {\bm{\delta}\mathbf{x}}^2 \otimes {\bm{\delta}\mathbf{x}}^2 \otimes {\bm{\delta}\mathbf{x}}^2\right] \right. \nonumber \\
    \,& \qquad  \qquad +6 E\left[{\bm{\delta}\mathbf{x}}^2 \otimes {\bm{\delta}\mathbf{x}}^2 \right]\otimes E\left[{\bm{\delta}\mathbf{x}}^2\right]^2 -4 E\left[{\bm{\delta}\mathbf{x}}^2 \otimes {\bm{\delta}\mathbf{x}}^2 \otimes{\bm{\delta}\mathbf{x}}^2\right] \otimes E \left[{\bm{\delta}\mathbf{x}}^2\right] \nonumber \\
    \,& \qquad \qquad \left. -3 E\left[{\bm{\delta}\mathbf{x}}^2\right]^4 \right\} +\cdots
\end{align}

Along the same lines, contributions of central point $\mathbf{X}_0$, through $\mathbf{Y}_0=\mathbf{M}(\mathbf{X}_0)$, on the vectorized first four output moments can be expressed as:

\begin{align}
	\mathbf{Y}_0-\bar{\mathbf{y}} &= - \dfrac{1}{2!}\left[\nabla^{2} \mathbf{M}\right]\left\{E\left[\bm{\delta}\mathbf{x}^2\right]\right\} - \cancel{\dfrac{1}{3!}\left[\nabla^{3} \mathbf{M}\right]\left\{E\left[\bm{\delta}\mathbf{x}^3\right]\right\}} - \dfrac{1}{4!}\left[\nabla^{4} \mathbf{M}\right]\left\{E\left[\bm{\delta}\mathbf{x}^4\right]\right\} \nonumber \\
                                    & -\cancel{\dfrac{1}{5!}\left[\nabla^{5} \mathbf{M}\right]\left\{E\left[\bm{\delta}\mathbf{x}^5\right]\right\}} - \dfrac{1}{6!}\left[\nabla^{6} \mathbf{M}\right]\left\{E\left[\bm{\delta}\mathbf{x}^6\right]\right\} - \cancel{\dfrac{1}{7!}\left[\nabla^{7} \mathbf{M}\right]\left\{E\left[\bm{\delta}\mathbf{x}^7\right]\right\}} -  \cdots
\end{align}

\begin{align}
	{{\left( \mathbf{Y}_0-\mathbf{\bar{y}} \right)}^{2}} & =\frac{1}{2!2!}\left[ \left({{\nabla}^{2}} \mathbf{M}\right)^{2} \right]\left\{E\left[\bm{\delta}\mathbf{x}^2 \right]^{2}\right\} + \frac{2}{2!3!}\left[ {{\nabla}^{2}} \mathbf{M}\otimes {{\nabla}^{3}}\mathbf{M} \right]\left\{\cancel{E\left[\bm{\delta}\mathbf{x}^2 \right] \otimes E\left[\bm{\delta}\mathbf{x}^3 \right]} \right\} \nonumber\\
	{}& +\frac{2}{2!4!}\left[ {{\nabla}^{2}} \mathbf{M}\otimes {{\nabla}^{4}}\mathbf{M} \right]\left\{E\left[\bm{\delta}\mathbf{x}^2 \right] \otimes E\left[\bm{\delta}\mathbf{x}^4 \right] \right\} + \frac{1}{3!3!}\left[ \left({{\nabla}^{3}} \mathbf{M}\right)^{2} \right]\left\{\cancel{E\left[\bm{\delta}\mathbf{x}^3 \right]^{2}} \right\} \nonumber\\ 
	{}& +\frac{2}{2!5!}\left[{{\nabla}^{2}} \mathbf{M}\otimes {{\nabla}^{5}}\mathbf{M} \right]\left\{\cancel{E\left[\bm{\delta}\mathbf{x}^2 \right] \otimes E\left[\bm{\delta}\mathbf{x}^5 \right]} \right\} +\frac{2}{3!4!}\left[{{\nabla}^{3}} \mathbf{M}\otimes {{\nabla}^{4}}\mathbf{M} \right]\left\{\cancel{E\left[\bm{\delta}\mathbf{x}^3 \right] \otimes E\left[\bm{\delta}\mathbf{x}^4 \right]} \right\} \nonumber\\ 
	{} & +\frac{2}{2!6!}\left[{{\nabla}^{2}} \mathbf{M}\otimes {{\nabla}^{6}}\mathbf{M} \right]\left\{E\left[\bm{\delta}\mathbf{x}^2 \right] \otimes E\left[\bm{\delta}\mathbf{x}^6 \right] \right\} +\frac{2}{3!5!}\left[{{\nabla}^{3}} \mathbf{M}\otimes {{\nabla}^{5}}\mathbf{M} \right]\left\{\cancel{E\left[\bm{\delta}\mathbf{x}^3 \right] \otimes E\left[\bm{\delta}\mathbf{x}^5 \right]} \right\} \nonumber\\
        {} & +\frac{1}{4!4!}\left[ \left({{\nabla}^{4}} \mathbf{M}\right)^{2} \right]\left\{E\left[\bm{\delta}\mathbf{x}^4 \right]^{2}\right\} + \cdots
\end{align}

\begin{align}
	{{\left( \mathbf{Y}_0-\mathbf{\bar{y}} \right)}^{3}} & =-\frac{1}{2!2!2!}\left[ \left({{\nabla}^{2}} \mathbf{M}\right)^{3} \right]\left\{E\left[\bm{\delta}\mathbf{x}^2 \right]^{3}\right\} - \frac{3}{2!2!3!}\left[ \left({{\nabla}^{2}} \mathbf{M}\right)^{2}\otimes {{\nabla}^{3}}\mathbf{M} \right]\left\{\cancel{E\left[\bm{\delta}\mathbf{x}^2 \right]^2 \otimes E\left[\bm{\delta}\mathbf{x}^3 \right]} \right\} \nonumber\\
	{}& -\frac{3}{2!2!4!}\left[ \left({{\nabla}^{2}} \mathbf{M}\right)^{2}\otimes {{\nabla}^{4}}\mathbf{M} \right]\left\{E\left[\bm{\delta}\mathbf{x}^2 \right]^2 \otimes E\left[\bm{\delta}\mathbf{x}^4 \right] \right\} \nonumber \\
        {}& -\frac{3}{2!3!3!}\left[ {{\nabla}^{2}} \mathbf{M}\otimes \left({{\nabla}^{3}}\mathbf{M}\right)^{2} \right]\left\{\cancel{E\left[\bm{\delta}\mathbf{x}^2 \right] \otimes E\left[\bm{\delta}\mathbf{x}^3 \right]^2} \right\} \nonumber \\
        {}& -\frac{1}{3!3!3!}\left[\left({{\nabla}^{3}} \mathbf{M}\right)^{3} \right]\left\{\cancel{E\left[\bm{\delta}\mathbf{x}^3 \right]^{3}}\right\} -\frac{3}{2!2!5!}\left[\left({{\nabla}^{2}} \mathbf{M}\right)^2\otimes {{\nabla}^{5}}\mathbf{M} \right]\left\{\cancel{E\left[\bm{\delta}\mathbf{x}^2 \right]^2 \otimes E\left[\bm{\delta}\mathbf{x}^5 \right]} \right\}  \nonumber\\ 
	{}& -\frac{6}{2!3!4!}\left[{{\nabla}^{2}} \mathbf{M}\otimes{{\nabla}^{3}}\mathbf{M} \otimes {{\nabla}^{3}}\mathbf{M} \right]\left\{\cancel{E\left[\bm{\delta}\mathbf{x}^2 \right] \otimes E\left[\bm{\delta}\mathbf{x}^3 \right] \otimes E\left[\bm{\delta}\mathbf{x}^4 \right]}\right\} - \cdots 
\end{align}

\begin{align}
	{{\left( \mathbf{Y}_0-\mathbf{\bar{y}} \right)}^{4}} & =\frac{1}{2!2!2!2!}\left[ \left({{\nabla}^{2}} \mathbf{M}\right)^{4} \right]\left\{E\left[\bm{\delta}\mathbf{x}^2 \right]^{4}\right\} + \frac{4}{2!2!2!3!}\left[ \left({{\nabla}^{2}} \mathbf{M}\right)^{3}\otimes {{\nabla}^{3}}\mathbf{M} \right]\left\{\cancel{E\left[\bm{\delta}\mathbf{x}^2 \right]^2 \otimes E\left[\bm{\delta}\mathbf{x}^3 \right]} \right\} \nonumber\\
	{}& +\frac{4}{2!2!2!4!}\left[ \left({{\nabla}^{2}} \mathbf{M}\right)^{3}\otimes {{\nabla}^{4}}\mathbf{M} \right]\left\{E\left[\bm{\delta}\mathbf{x}^2 \right]^3 \otimes E\left[\bm{\delta}\mathbf{x}^4 \right] \right\} \nonumber \\
        {}& +\frac{6}{2!2!3!3!}\left[ \left({{\nabla}^{2}} \mathbf{M}\right)^{2}\otimes \left({{\nabla}^{3}}\mathbf{M}\right)^2 \right]\left\{\cancel{E\left[\bm{\delta}\mathbf{x}^2 \right]^2 \otimes E\left[\bm{\delta}\mathbf{x}^3 \right]^2} \right\} \nonumber \\
        {}& +\frac{4}{2!2!2!5!}\left[ \left({{\nabla}^{2}} \mathbf{M}\right)^{3}\otimes {{\nabla}^{5}}\mathbf{M} \right]\left\{\cancel{E\left[\bm{\delta}\mathbf{x}^2 \right]^3 \otimes E\left[\bm{\delta}\mathbf{x}^5 \right]} \right\} \nonumber \\
        {}& +\frac{12}{2!2!3!4!}\left[ \left({{\nabla}^{2}} \mathbf{M}\right)^{2}\otimes {{\nabla}^{3}}\mathbf{M}\otimes {{\nabla}^{4}}\mathbf{M} \right]\left\{\cancel{E\left[\bm{\delta}\mathbf{x}^2 \right]^2 \otimes E\left[\bm{\delta}\mathbf{x}^3 \right]\otimes E\left[\bm{\delta}\mathbf{x}^4 \right]} \right\} + \cdots
\end{align}

\bibliography{reference_qpem.bib}

\end{document}